\documentclass[english,reqno]{amsart}

\usepackage{diagbox}
\usepackage{makecell}
\usepackage{amsmath, appendix, ulem, amsthm}
\usepackage[nobysame]{amsrefs}
\usepackage{amssymb, color,multirow}
\usepackage{amsfonts}
\usepackage[margin=1in]{geometry}
\usepackage{mathrsfs}
\usepackage{graphicx}
\usepackage{float}
\usepackage{epsf}
\usepackage{titletoc}
\usepackage{cases}
\usepackage[linesnumbered,ruled]{algorithm2e}
\usepackage{subfigure}
\usepackage{bbold}
\usepackage{dsfont}
\graphicspath{{../Figures/}}
\usepackage{appendix}
\numberwithin{equation}{section}

\newtheorem{thm}{Main Theorem}
\theoremstyle{remark}

\newtheorem{remark}[thm]{Remark}

\renewcommand{\hat}{\widehat}

\newcommand{\f}{\frac}

\newcommand{\beq}{\begin{equation}}
\newcommand{\eeq}{\end{equation}}

\newcommand{\eps}{\varepsilon}

\newcommand{\rr}{\boldsymbol{R}}

\newcommand{\Ps}{{\mathbb P}}

\newcommand{\PP}{\mathcal{P}}

\newcommand{\rd}{\mathrm{d}}

\newcommand{\Dx}{\Delta x}

\newcommand{\order}{\mathcal{O}}
\newcommand{\average}[1]{\left\langle#1\right\rangle}
\newcommand{\dc}{D_d}

\begin{document}

\title{A Spatial-Temporal asymptotic preserving scheme  for   radiation magnetohydrodynamics in the equilibrium and non-equilibrium diffusion limit} 
\author{Shi Jin}
\address{School of Mathematics, Institute of Natural Sciences and MOE-LSC, Shanghai Jiao Tong University, China}
\email{shijin-m@sjtu.edu.cn}
\author{Min Tang}
\address{School of Mathematics, Institute of Natural Sciences and MOE-LSC, Shanghai Jiao Tong University, China}
\email{tangmin@sjtu.edu.cn}
\author{Xiaojiang Zhang}
\address{School of Mathematics, Institute of Natural Sciences and MOE-LSC, Shanghai Jiao Tong University, China}
\email{xjzhang123@sjtu.edu.cn}
%
%
\date{\today}
\begin{abstract}
The radiation magnetohydrodynamics (RMHD) system couples the ideal magnetohydrodynamics equations with a gray radiation transfer equation. The main challenge is that the radiation travels at the speed of light while the magnetohydrodynamics changes with the time scale of the fluid. The time scales of these two processes can vary dramatically. In order to use mesh sizes and time steps that are independent of the speed of light, asymptotic preserving (AP) schemes in both space and time are desired. In this paper, we develop an AP scheme in both space and time for the RMHD system. Two different scalings are considered. One results in an  equilibrium diffusion limit system, while the other results in a non-equilibrium system. The main idea is to decompose the radiative intensity into three parts, each part is treated differently with suitable combinations of explicit and implicit discretizations guaranteeing the favorable stability conditionand computational efficiency. The performance of the AP method is  presented, for both optically thin and thick regions, as well as for the radiative shock problem. 
\end{abstract}

\maketitle


\section{Introduction}

Radiation magnetohydrodynamics (RMHD) is concerned with the dynamical behavior of magnetized  fluids that have nonnegligible exchange of energy and momentum  with radiation, which is important in high temperature flow systems, solar and space physics and astrophysics.
The radiation transfer equation (RTE) for the radiative intensity $I $ and the fluid temperature $T$ in the mixed frame (the radiative intensity $I$ is described in an Eulerian frame while the material coupling terms are described in the comoving frame) is governed by 
\begin{equation}\label{rtmhd}
\begin{aligned}
\frac{\partial I}{\partial t}+C\boldsymbol{n}\cdot\nabla I=C\sigma_a\left( \f{a_rT^4}{4\pi}-I\right)&+C\sigma_s(J-I)+3\boldsymbol{n}\cdot \boldsymbol{v}\sigma_a\left(\f{a_rT^4}{4\pi}-J\right)+\boldsymbol{n}\cdot \boldsymbol{v}(\sigma_a+\sigma_s)(I+3J)\\
&-2\sigma_s \boldsymbol{v}\cdot \boldsymbol{H}-(\sigma_a-\sigma_s)\f{\boldsymbol{v}\cdot \boldsymbol{v}}{C}J-(\sigma_a-\sigma_s)\f{\boldsymbol{v}\cdot (\boldsymbol{v}\cdot {\rm K})}{C},
\end{aligned}
\end{equation}
where $\sigma_a$ and $\sigma_s$ represent absorption and scattering opacities respectively, $\boldsymbol{n}\in V$  is the angular variable, $a_r$ is the radiation constant and $C$ is the speed of  light \cite{pub.1041264084}. $J$, $\boldsymbol{H}, {\rm K}$ are the zeroth, first, second angular moments of $I$, respectively:
$$
J=\f{1}{|V|}\int_{V} I  \text{d}\boldsymbol{n},\qquad
\boldsymbol{H}=\f{4\pi C}{|V|}\int_{V} \boldsymbol{n}I \text{d}\boldsymbol{n},\qquad
{\rm K}=\f{1}{|V|}\int_{V} \boldsymbol{n}\boldsymbol{n}I \text{d}\boldsymbol{n}.
$$
The ideal MHD equations with radiation energy and momentum source terms are
\begin{equation}
\left\{
\begin{aligned}
& \frac{\partial\rho}{\partial t}+\nabla\cdot(\rho \boldsymbol{v})=0,\\
& \frac{\partial(\rho \boldsymbol{v})}{\partial t}+\nabla\cdot(\rho \boldsymbol{v}\boldsymbol{v}-\boldsymbol{B}\boldsymbol{B}+\Ps^*)=-\boldsymbol{S_{rp}},\\
& \frac{\partial E}{\partial t}+\nabla\cdot[(E+P^*)\boldsymbol{v}-\boldsymbol{B}(\boldsymbol{B}\cdot \boldsymbol{v})]=-CS_{re},\\
&   \frac{\partial \boldsymbol{B}}{\partial t}+\nabla\times(\boldsymbol{v}\times \boldsymbol{B})=0\,,
\end{aligned}
\right.
\label{eqn:003}
\end{equation}
where
\begin{equation*}
S_{re}=\sigma_a( a_rT^4-E_r)+(\sigma_a-\sigma_s)\f{\boldsymbol{v}}{C^2}\cdot\left[\boldsymbol{F_r}-(\boldsymbol{v}E_r+\boldsymbol{v}\cdot {\rm P_r})\right],
\end{equation*}
\begin{equation*}
\boldsymbol{S_{rp}}=-\frac{\sigma_s+\sigma_a}{C}\left[\boldsymbol{F_r}-(\boldsymbol{v}E_r+\boldsymbol{v}\cdot {\rm P_r})\right]+\f{\boldsymbol{v}}{C}\sigma_a(a_rT^4-E_r),
\end{equation*}
with the energy density
\begin{equation*}
E_r= 4\pi J= \f{4\pi}{|V|}\int_{V} I  \text{d}\boldsymbol{n}=4\pi\average{I},
\end{equation*}
 the radiation flux
\begin{equation*}
\boldsymbol{F_r}= 4\pi C\boldsymbol{H}= \f{4\pi C}{|V|}\int_{V} \boldsymbol{n}I \text{d}\boldsymbol{n}=4\pi C\average{\boldsymbol{n}I},
\end{equation*}
and the radiation pressure
\begin{equation*}
{\rm P_r}= 4\pi {\rm K}=\f{4\pi}{|V|}\int_{V} \boldsymbol{n}\boldsymbol{n}I \text{d}\boldsymbol{n}=4\pi\average{\boldsymbol{n}\boldsymbol{n}I}.
\end{equation*}
In particular, the operators $\average{\cdot},\average{\boldsymbol{n}\cdot}$ and  $\average{\boldsymbol{n}\boldsymbol{n}\cdot}$ are respectively $\f{1}{|V|}\int_{V} \cdot  \text{d}\boldsymbol{n}$,  $\f{1}{|V|}\int_{V}\boldsymbol{n} \cdot  \text{d}\boldsymbol{n}$ and $\f{1}{|V|}\int_{V}\boldsymbol{n}\boldsymbol{n} \cdot  \text{d}\boldsymbol{n}$ in this context. Moreover, $\rho$ is the fluid density; $\boldsymbol{v}$ is the fluid velocity; $\boldsymbol{B}$ is the magnetic flux density; $p$ is the static pressure; $\Ps^*\equiv (p+B^2/2)I_d$ (with $I_d$ the $d\times d$ identity matrix) is the full pressure tensor of the fluid; $P^*\equiv (p+B^2/2)$  is the full pressure  of the fluid;
 \begin{equation*}
E=E_g+\frac{1}{2}\rho v^2+\frac{B^2}{2},
\end{equation*}
where $E_g$ is the internal energy density; $v^2=\boldsymbol{v}\cdot\boldsymbol{v}$ and $B^2=\boldsymbol{B}\cdot\boldsymbol{B}$. The radiation MHD system is closed by the perfect gas equation of state:
$$
E_g=p/(\gamma-1),\ \qquad \ T=p/(R_{\text{ideal}}\rho),
$$
where $\gamma$ is adiabatic index for an ideal gas, and $R_{\text{ideal}}$ is the ideal gas constant.

Numerical approximations to  the RTE have been extensively studied in \cite{MCCLARREN20087561,MOREL1996445,sun2015asymptotic,jin2009uniformly,KFJ16,klar1998asymptotic} and  for ideal MHD see for example \cite{stone2008athena,brio1988upwind,li2011central}. The RMHD coupled system  has been investigated lately \cite{pub.1005319601,BOLDING2017511,SEKORA20106819,KADIOGLU20108313,pub.1041264084,sekora2009higher,sun2020multiscale}. Since the transport term in the RTE and the advection term in the MHD equations are at different time scales, the RMHD system can be stiff, as can been seen more clearly after nondimensionalization \cites{lowrie1999coupling}. To solve the RMHD system, classical numerical discretizations require the space and time
steps to resolve the speed of light, which is very expensive. Asymptotic-Preserving (AP) schemes provide a generic framework for  such multiscale problems \cite{jin1999efficient, JinReview}. AP schemes were first studied for steady neutron transport problems in diffusive regimes \cite{larsen1987asymptotic,larsen1989asymptotic} and then \cite{golse1999convergence,2018Tang} for boundary value problems, unsteady transport problems \cite{jin1998diffusive, klar1998asymptotic} and gray radiative transport equations \cite{2001Klar,2021Tang}. When the  scaling parameter $\eps$ in the multiscale system can not be resolved numerically,  AP schemes should automatically become  a good solver for the macroscopic models. 

  AP schemes for RMHD have been proposed in the literature, but since most of them are based on operator splitting, they are only AP in space.  In \cite{BOLDING2017511}, Simon et al. have developed a second-order  scheme in both space and time for  the Euler equations coupled with a gray radiation $S_2$ model. The scheme in \cite{BOLDING2017511} uses MUSCL-Hancock method  to solve the Euler equations and the lumped linear-discontinuous Galerkin method for the radiation $S_2$ model. For the time discretization,  the TR/BDF2 time integration method is employed. Jiang {\it et al.}  devised an Implicit/Explicit (IMEX) scheme for the splitting system  that treats the transport term in the RTE and  convective term in the fluid equations explicitly, and the source term implicitly \cite{pub.1041264084}.  Recently, another method by Sun {\it et al.} \cite{sun2020multiscale} solves a coupled system of RTE and Euler equations, the authors use the gas-kinetic scheme (GKS) to solve the Euler equations and disretize the RTE by the unified gas-kinetic scheme (UGKS). However, all the aforementioned methods either require a time step that satisfies a hyperbolic constraint with CFL number proportional to the speed of light \cite{sun2020multiscale}, or use nonlinear iterations involving \eqref{rtmhd} in order to get rid of the time step constraint \cite{BOLDING2017511, pub.1041264084}. Both approaches are expensive. It is desired to design a scheme that can use a time step that is independent of the speed of light and employs nonlinear iterations that involve only macroscopic quantities, i.e. the $\rho$, $\boldsymbol{v}$, $E$, $\boldsymbol{B}$ or the moments of radiative intensity. 
  
  In this paper, we aim at developing a scheme for RMHD that is AP in $\it{both}$ space $\it{and}$ time. Two different parameter regimes are considered. One is for $\sigma_a$ small and $\sigma_s$ large, which results in the non-equilibrium diffusion limit system. The other is for $\sigma_a$ large, which gives the equilibrium diffusion limit system. The proposed scheme preserves both limits. The main idea is to decompose the intensity into three parts, two parts correspond to the zeroth and first order moments, while the third part is the residual. The two macroscopic moments are treated implicitly while the residual is explicit. Thanks to the properties of the residual term, one can update the macroscopic quantities first. The residual term is then updated by solving implicitly a linear transport equation for each direction $\mathbf{n}$. For the space discretization, we use Roe's method in the Athena code to solve the convective part in the ideal MHD equations and the UGKS  for the RTE. AP property of both semi-discretized and fully discretizad systems are proved. 

The computational cost of our proposed scheme is lower due to the following reasons. Compared to the explicit scheme, much larger time steps are allowed. On the other hand, compared to the fully implicit scheme, our scheme requires nonlinear iterations to solve a system with only macroscopic quantities and then one linear transport equation for each direction $\mathbf{n}$. Since different directions are {\it decoupled}, the computational cost is much lower than solving implicitly a linear RTE where all velocities are coupled. The computational bottleneck due to the resolution of the angular dimension is avoided.
  
  This paper is organized as follows. Section 2 gives the asymptotic limit of the coupled system, and the equilibrium and  non-equilibrium diffusion limit systems under different scalings are obtained.
 The  semi-discretization and one dimensional fully discretized
scheme for RMHD are illustrated in Sections 3 and 4 respectively. Their capability of capturing both diffusion limits are proved. In Section 5, the performance of our AP method is presented in both optically thin and thick regions. The radiative shock problems are tested. By comparing with  the semi-analytic solutions from \cite{lowrie2008radiative}, we observe that the scheme is accurate and stable using large time step and meshes that are independent of the speed of light.  Finally,  we conclude in Section 6.
\section{The diffusion limit of RMHD}\label{dimension}
As presented in \cite{lowrie1999coupling}, the dimensionless form for the coupled system can give a better insight on the relative importance of  different terms.  Consider the following nondimensionalization:
$$
x=\hat{x}\ell_\infty,\quad t=\hat{t}\ell_\infty/a_\infty,\quad \rho=\hat{\rho}\rho_\infty,\quad \boldsymbol{v}=\hat{\boldsymbol{v}}a_\infty,\quad p=\hat{p}\rho_\infty a_\infty^2,\quad T=\hat{T}T_\infty, \quad I=a_rT_\infty^4\hat{I},
$$
$$
E_r=a_rT_\infty^4\hat{E_r},\quad \boldsymbol{F_r}=Ca_rT_\infty^4\hat{\boldsymbol{F_r}},\quad {\rm P_r}=a_rT_\infty^4\hat{{\rm P_r}},\quad\sigma_a=\lambda_a\hat{\sigma_a}, \quad\sigma_s=\lambda_s\hat{\sigma_s},
$$
where  variables with a hat denote nondimensional quantities and variables with $\infty- $subscript are the characteristic values with units.  More precisely, $\ell_\infty$, $a_\infty$, $\rho_\infty$ and $T_\infty$ are respectively the reference length, sound speed, density and temperature. Moreover, $\lambda_a$ and $\lambda_s$ are respectively the characteristic values of absorption and scattering coefficients. Then the full dimensionless radiation MHD system becomes (after dropping the hat):
\begin{equation}
\left\{
\begin{aligned}
&\frac{\partial I}{\partial t}+\mathcal{C}\boldsymbol{n}\cdot\nabla I=\mathscr{L}_a\mathcal{C}\sigma_a\left( \f{T^4}{4\pi}-I\right)+\mathscr{L}_s\mathcal{C}\sigma_s(J-I)+3\mathscr{L}_a \boldsymbol{n}\cdot \boldsymbol{v}\sigma_a\left(\f{T^4}{4\pi}-J\right)+\boldsymbol{n}\cdot \boldsymbol{v}\left(\mathscr{L}_a\sigma_a+\mathscr{L}_s\sigma_s\right)(I+3J)\\
&\hspace{3.5cm}-2\mathscr{L}_s\sigma_s\boldsymbol{v}\cdot \boldsymbol{H}-(\mathscr{L}_a\sigma_a-\mathscr{L}_s\sigma_s)\f{\boldsymbol{v}\cdot \boldsymbol{v}}{\mathcal{C}}J-(\mathscr{L}_a\sigma_a-\mathscr{L}_s\sigma_s)\f{\boldsymbol{v}\cdot (\boldsymbol{v}\cdot {\rm K})}{\mathcal{C}}\triangleq\mathcal{C}S,\\
& \frac{\partial\rho}{\partial t}+\nabla\cdot(\rho \boldsymbol{v})=0,\\
& \frac{\partial(\rho \boldsymbol{v})}{\partial t}+\nabla\cdot(\rho \boldsymbol{v}\boldsymbol{v}-\boldsymbol{B}\boldsymbol{B}+\Ps^*)=-\PP_0\boldsymbol{S_{rp}},\\
& \frac{\partial E}{\partial t}+\nabla\cdot[(E+P^*)\boldsymbol{v}-\boldsymbol{B}(\boldsymbol{B}\cdot \boldsymbol{v})]=-\mathcal{C}\PP_0S_{re},\\
&   \frac{\partial \boldsymbol{B}}{\partial t}+\nabla\times(\boldsymbol{v}\times \boldsymbol{B})=0\,,
\end{aligned}
\right.
\label{eqn:non}
\end{equation}
where $\mathcal{C}=\f{C}{a_\infty}$, $\PP_0=\f{a_rT_\infty^4}{\rho_\infty a^2_\infty}$, $\mathscr{L}_a=\ell_\infty\lambda_a$, $\mathscr{L}_s=\ell_\infty\lambda_s$,  $S_{re}=4\pi\average{S}$ and $\boldsymbol{S_{rp}}=4\pi\average{\boldsymbol{n}S}$.
In \eqref{eqn:non}, $\PP_0$ is a nondimensional constant, which  measures the influence  of  radiation on the flow dynamics. In  \cite{godillon2005coupled,BCDB1,BCDB2}, the authors derived the equilibrium diffusion and non-equilibrium diffusion limits for a system that is composed of a kinetic equation and a drift-diffusion equation or Euler equations. In \cite{lowrie1999coupling}, the authors considered four different regimes for the RMHD system: equilibrium diffusion regime, strong equilibrium regime, isothermal regime and streaming regime. In the subsequent part, we  will focus on two regimes of the RMHD system: equilibrium diffusion and non-equilibrium diffusion limit regimes, whose derivations are the same as in \cite{godillon2005coupled,BCDB1,BCDB2,lowrie1999coupling}. For the convenience of readers, the details of the derivation are in Appendix \ref{apen01}.
\subsection{The non-equilibrium diffusion limit.} \label{sec21}
In this regime the radiation intensity $I$ tends to a state which can be characterized by a temperature which is different from the material temperature $T$.
As in the reference, we consider the following scaling
$$
\mathscr{L}_a = \eps,\quad\mathscr{L}_s = 1/\eps,\quad\PP_0=\order(1),\quad\mathcal{C}=c/\eps.
$$
Here the introduction of $\eps$ is for the convenience of asymptotic analysis, and the value of $\eps$ depends on the speed of light.
The  RMHD system \eqref{eqn:non} becomes:
\begin{subequations}\label{eqn:004}
\begin{numcases}{}
\frac{\partial I}{\partial t}+\frac{c}{\varepsilon}\boldsymbol{n}\cdot\nabla I=c\sigma_a\left( \f{T^4}{4\pi}-I\right)+\frac{c\sigma_s}{\varepsilon^2}(J-I)+3\eps \boldsymbol{n}\cdot \boldsymbol{v}\sigma_a\left(\f{T^4}{4\pi}-J\right)+\boldsymbol{n}\cdot \boldsymbol{v}\left(\eps\sigma_a+\f{\sigma_s}{\eps}\right)(I+3J)\nonumber\\
 \hspace{5cm}-2\f{\sigma_s}{\eps} \boldsymbol{v}\cdot \boldsymbol{H}-(\eps^2\sigma_a-\sigma_s)\f{\boldsymbol{v}\cdot \boldsymbol{v}}{c}J-(\eps^2\sigma_a-\sigma_s)\f{\boldsymbol{v}\cdot (\boldsymbol{v}\cdot {\rm K})}{c}, \label{eqn:case11004}\\
 \frac{\partial\rho}{\partial t}+\nabla\cdot(\rho \boldsymbol{v})=0,\label{eqn:case12004}\\
 \frac{\partial(\rho \boldsymbol{v})}{\partial t}+\nabla\cdot(\rho \boldsymbol{v}\boldsymbol{v}-\boldsymbol{BB}+\Ps^*)=-\PP_0\boldsymbol{S_{rp}},\label{eqn:case13004}\\
 \frac{\partial E}{\partial t}+\nabla\cdot[(E+P^*)\boldsymbol{v}-\boldsymbol{B}(\boldsymbol{B}\cdot \boldsymbol{v})]=-\frac{c}{\varepsilon}\PP_0S_{re}, \label{eqn:case14004}\\
   \frac{\partial \boldsymbol{B}}{\partial t}+\nabla\times(\boldsymbol{v}\times \boldsymbol{B})=0 \label{eqn:case15004}\,,
 \end{numcases}
\end{subequations}
where
\begin{equation*}
S_{re}=\eps\sigma_a\left( T^4-4\pi J\right)+(\eps^3 \sigma_a-\eps\sigma_s)\f{4\pi \boldsymbol{v}}{c^2}\cdot\left[\f{c}{\eps}\boldsymbol{H}-\left(\boldsymbol{v}J+\boldsymbol{v} \cdot {\rm K}\right)\right]	, \
\end{equation*}
\begin{equation*}
\boldsymbol{S_{rp}}=-\frac{4\pi(\sigma_s+\eps^2\sigma_a)}{c}\left[\f{c}{\eps}\boldsymbol{H}-\left(\boldsymbol{v}J+\boldsymbol{v} \cdot {\rm K}\right)\right]+\f{\boldsymbol{v}}{c}\eps^2\sigma_a\left(T^4-4\pi J\right).
\end{equation*}
When $\eps\to 0$ in (\ref{eqn:004}), the solution can be approximated by the solution of the following non-equilibrium system:
\begin{subequations}\label{eqnlimit1}
\begin{numcases}{}
\partial_t\rho +\nabla\cdot(\rho \boldsymbol{v})=0, \label{eqn:eqnlimit11}\\
 \partial_t(\rho \boldsymbol{v})+\nabla\cdot\left(\rho \boldsymbol{v}\boldsymbol{v}-\boldsymbol{B}\boldsymbol{B}+\Ps^*\right)=-4\pi\PP_0D_d\nabla J,\label{eqn:eqnlimit12}\\
 \partial_t\left(E+4\pi\PP_0J\right)+\nabla\cdot\left[(E+P^*)\boldsymbol{v}-\boldsymbol{B}(\boldsymbol{B}\cdot \boldsymbol{v})+16\pi\PP_0D_d\boldsymbol{v}J\right]=\nabla\cdot\left(\f{4\pi c\PP_0D_d}{\sigma_s}\nabla J\right),\label{eqn:eqnlimit13}\\
 4\pi\partial_tJ+\nabla\cdot\left(16\pi D_d\boldsymbol{v}J-\f{ 4\pi cD_d}{\sigma_s}\nabla J\right)=c\sigma_a\left(T^4-4\pi J\right)+4\pi D_d\boldsymbol{v}\cdot\nabla J, \label{eqn:eqnlimit14}\\
   \partial_t \boldsymbol{B}+\nabla\times(\boldsymbol{v}\times \boldsymbol{B})=0. \label{eqn:eqnlimit15}\,
 \end{numcases}
\end{subequations}
Here, the non-equilibrium system indicates that $J$ is away from $\f{T^4}{4\pi}$, while we will see in section \ref{sec2} that the equilibrium diffusion limit indicates that $J\approx\f{T^4}{4\pi}$.

\subsection{The equilibrium diffusion limit.}\label{sec2}
In this regime  the radiation intensity $I$ adapt to the material temperature. As in the reference, we consider the following scaling
$$
\mathscr{L}_a = 1/\eps,\quad\mathscr{L}_s = \eps,\quad\PP_0=\order(1),\quad\mathcal{C}=c/\eps.
$$
This is called the equilibrium diffusion regime, which will give the classical equilibrium diffusion limit \cites{pomraning2005equations,mihalas1984foundations,lowrie1999coupling}. The  RMHD system \eqref{eqn:non} becomes:
\begin{subequations}\label{eqn:eqb}
\begin{numcases}{}
\frac{\partial I}{\partial t}+\frac{c}{\varepsilon}\boldsymbol{n}\cdot\nabla I=\frac{c\sigma_a}{\varepsilon^2}\left( \f{T^4}{4\pi}-I\right)+c\sigma_s(J-I)+3\boldsymbol{n}\cdot \boldsymbol{v}\f{\sigma_a}{\eps}\left(\f{T^4}{4\pi}-J\right)+\boldsymbol{n}\cdot \boldsymbol{v}\left(\f{\sigma_a}{\eps}+\eps\sigma_s\right)(I+3J)\nonumber\\
\hspace{5cm}-2\eps\sigma_s \boldsymbol{v}\cdot \boldsymbol{H}-(\sigma_a-\eps^2\sigma_s)\f{\boldsymbol{v}\cdot \boldsymbol{v}}{c}J-(\sigma_a-\eps^2\sigma_s)\f{\boldsymbol{v}\cdot (\boldsymbol{v}\cdot {\rm K})}{c}, \label{eqn:case11eqb}\\
 \frac{\partial\rho}{\partial t}+\nabla\cdot(\rho \boldsymbol{v})=0,\label{eqn:case12eqb}\\
 \frac{\partial(\rho \boldsymbol{v})}{\partial t}+\nabla\cdot(\rho \boldsymbol{v}\boldsymbol{v}-\boldsymbol{B}\boldsymbol{B}+\Ps^*)=-\PP_0\boldsymbol{S_{rp}},\label{eqn:case13eqb}\\
 \frac{\partial E}{\partial t}+\nabla\cdot[(E+P^*)\boldsymbol{v}-\boldsymbol{B}(\boldsymbol{B}\cdot \boldsymbol{v})]=-\frac{c}{\varepsilon}\PP_0S_{re}, \label{eqn:case14eqb}\\
   \frac{\partial \boldsymbol{B}}{\partial t}+\nabla\times(\boldsymbol{v}\times \boldsymbol{B})=0 \label{eqn:case15eqb}\,,
 \end{numcases}
\end{subequations}
where
\begin{equation*}
S_{re}=\frac{\sigma_a}{\varepsilon}\left( T^4-4\pi J\right)+(\eps\sigma_a-\eps^3\sigma_s)\f{4\pi \boldsymbol{v}}{c^2}\cdot\left[\f{c}{\eps}\boldsymbol{H}-\left(\boldsymbol{v}J+\boldsymbol{v} \cdot {\rm K}\right)\right],
\end{equation*}
\begin{equation*}
\boldsymbol{S_{rp}}=-\frac{4\pi(\eps^2\sigma_s+\sigma_a)}{c}\left[\f{c}{\eps}\boldsymbol{H}-\left(\boldsymbol{v}J+\boldsymbol{v} \cdot {\rm K}\right)\right]+\f{\boldsymbol{v}}{c}\sigma_a\left(T^4-4\pi J\right).
\end{equation*}
When $\eps\to 0$ in (\ref{eqn:eqb}), the solution can be approximated by the solution of the following equilibrium system:\begin{subequations}\label{eqnlimit2}
\begin{numcases}{}
\partial_t\rho +\nabla\cdot(\rho \boldsymbol{v})=0, \label{eqn:eqnlimit21}\\
\partial_t(\rho \boldsymbol{v})+\nabla\cdot\left(\rho \boldsymbol{v}\boldsymbol{v}-\boldsymbol{B}\boldsymbol{B}+\Ps^*\right)=-\PP_0D_d\nabla T^4,\label{eqn:eqnlimit22}\\
 \partial_t\left(E+\PP_0T^4\right)+\nabla\cdot\left[(E+P^*)\boldsymbol{v}-\boldsymbol{B}(\boldsymbol{B}\cdot \boldsymbol{v})+4\PP_0D_d\boldsymbol{v}T^4\right]=\nabla\cdot\left(\f{ c\PP_0D_d}{\sigma_a}\nabla T^4\right),\label{eqn:eqnlimit23}\\
   \partial_t \boldsymbol{B}+\nabla\times(\boldsymbol{v}\times \boldsymbol{B})=0. \label{eqn:eqnlimit24}\,
 \end{numcases}
\end{subequations}
As we can see, when $J=\f{T^4}{4\pi}$ in \eqref{eqnlimit1}, the non-equilibrium diffusion system \eqref{eqnlimit1} is the same as the equilibrium diffusion system \eqref{eqnlimit2}, which indicates $J$ is close to $\f{T^4}{4\pi}$ in  the equilibrium system \eqref{eqn:eqb}.
\begin{remark}
To understand the above scalings, we choose the light speed,  the characteristic material temperature and the radiation constant to be respectively $C=299.8 \ cm\  sh^{-1}$ ($1\ sh=10^{-8}\ s$),  $0.1 keV$ and
$a_r=0.0001372 \ Jk \ cm^{-3} \ keV^{-4}$. The typical computational domain is 0.04 $cm$, and the units of the length, time, temperature and energy are respectively $cm,\ sh, \ keV$ and $Jk$ ($1\ Jk=10^9\ J$). For the above parameters, we can choose $\eps=0.01$. Then, when $\sigma_a=0.01 \ cm^{-1}$, $\sigma_s=577.35\ cm^{-1}$, the parameters are in the non-equilibrium diffusion regime and when $\sigma_a=577.35 \ cm^{-1}$, $\sigma_s=0.01\ cm^{-1}$, in the equilibrium diffusion regime. 
	\end{remark}




\section{Time discretization for   the  RMHD }
In this section, based on a decomposition for radiation intensity $I$, we will present a semi-discretization for the RMHD and show its AP property for any dimensions in space and angular variables.
 \subsection{A decomposition}  To preserve the asymptotic limit of \eqref{eqn:non}, the main idea is to decompose the intensity $I$ into three  parts, and each of them is treated in a different way. More precisely, we decompose $I(t,x,n)$ as follows:
 \begin{subequations}\label{eqn:011}
 \begin{equation}\label{eqn:0111}
I(t,x,\boldsymbol{n})=\average{I}+3 \boldsymbol{n}\cdot\average{\boldsymbol{n}I}+ Q(t,x,\boldsymbol{n}),
 \end{equation}
  \begin{equation} \label{eqn:0112}
\hspace{2.2cm}:=J(t,x)+\boldsymbol{n}\cdot \rr(t,x)+ Q(t,x,\boldsymbol{n}),
 \end{equation}
 \end{subequations}
where 3$\average{\boldsymbol{n}I}\equiv\rr$. 
Then taking $\average{\boldsymbol{n}\cdot}$  on both sides of \eqref{eqn:0111}  yields
\begin{equation*}
\average{\boldsymbol{n}I}=\average{\boldsymbol{n}I}+\average{\boldsymbol{n}Q},
\end{equation*}
which gives
\begin{equation}\label{eqn:rr}
\average{\boldsymbol{n}Q}=0.
\end{equation}
Moreover,  taking $\average{\cdot}$  on \eqref{eqn:0112}, yields
$$
\average{I}=\average{J+ \boldsymbol{n}\cdot \rr+ Q}=\average{J+ Q}=\average{I}+\average{ Q},
$$
which implies
\begin{equation}\label{eqn:qq}
\average{ Q}=0.
\end{equation}
The idea is close to the micro-macro decomposition in \cite{klar1998asymptotic,LM}. However, we do not write down the macroscopic and microscopic equations explicitly, which is different from the classical micro-macro decomposition method as in \cite{klar1998asymptotic}.

Substituting  the decomposition (\ref{eqn:011}) into  system (\ref{eqn:non}) leads to
\begin{subequations}\label{eqn:0012}
\begin{numcases}{}
\partial_t\left(J+\boldsymbol{n}\cdot\rr+Q\right)+\mathcal{C}\boldsymbol{n}\cdot\nabla \left(J+\boldsymbol{n}\cdot\rr+Q\right)=\mathscr{L}_a\mathcal{C}\sigma_a\left(\f{T^4}{4\pi}-J- \boldsymbol{n}\cdot\rr - Q\right)-\mathscr{L}_s\mathcal{C}\sigma_s(\boldsymbol{n}\cdot\rr+Q)\nonumber\\
\hspace{3cm}-\f{2}{3}\mathscr{L}_s\sigma_s\boldsymbol{v}\cdot\rr+\boldsymbol{n}\cdot \boldsymbol{v}\left(\mathscr{L}_a\sigma_a+\mathscr{L}_s\sigma_s\right)(4J+\boldsymbol{n}\cdot\rr+Q)+3\mathscr{L}_a\sigma_a\boldsymbol{n}\cdot \boldsymbol{v}\left(\f{T^4}{4\pi}-J\right)\nonumber\\
\hspace{3cm}-\f{(\mathscr{L}_a\sigma_a-\mathscr{L}_s\sigma_s)}{\mathcal{C}}\left(\boldsymbol{v}\cdot\boldsymbol{v}J+\boldsymbol{v}\cdot\left(\boldsymbol{v}\cdot D_d\right)J+ \boldsymbol{v}\cdot\left(\boldsymbol{v}\cdot{\rm K_Q}\right)\right), \label{eqn:case110012}\\
 \frac{\partial\rho}{\partial t}+\nabla\cdot(\rho \boldsymbol{v})=0,\label{eqn:case120012}\\
 \frac{\partial(\rho \boldsymbol{v})}{\partial t}+\nabla\cdot(\rho \boldsymbol{vv}-\boldsymbol{BB}+\Ps^*)=-\PP_0\hat{\boldsymbol{S_{rp}}},\label{eqn:case130012}\\
 \frac{\partial E}{\partial t}+\nabla\cdot[(E+P^*)\boldsymbol{v}-\boldsymbol{B}(\boldsymbol{B}\cdot \boldsymbol{v})]=-\mathcal{C}\PP_0\hat{S_{re}}, \label{eqn:case140012}\\
   \frac{\partial \boldsymbol{B}}{\partial t}+\nabla\times(\boldsymbol{v}\times \boldsymbol{B})=0 \label{eqn:case150012}\,,
 \end{numcases}
\end{subequations}
where $D_d=\average{\boldsymbol{n}\boldsymbol{n}}=\f{1}{3}I_d$\ ($I_d$ denotes the 3 by 3 identity matrix), $ {\rm K_Q} =\average{\boldsymbol{n}\boldsymbol{n}Q}$  and 
\begin{equation*}
\hat{S_{re}}=\mathscr{L}_a\sigma_a(T^4-4\pi J)+4\pi(\mathscr{L}_a\sigma_a-\mathscr{L}_s\sigma_s)\f{\boldsymbol{v}}{\mathcal{C}^2}\cdot\left[\mathcal{C}D_d\rr-\left(\boldsymbol{v}J+D_d\boldsymbol{v}J+ {\rm K_Q}\boldsymbol{v}\right)\right], \
\end{equation*}
\begin{equation*}
\hat{\boldsymbol{S_{rp}}}=-\frac{4\pi(\mathscr{L}_s\sigma_s+\mathscr{L}_a\sigma_a)}{\mathcal{C}}\left[\mathcal{C}D_d\rr-\left(\boldsymbol{v}J+\dc\boldsymbol{v}J+ {\rm K_Q}\boldsymbol{v}\right)\right]+\f{ \boldsymbol{v}}{\mathcal{C}}\mathscr{L}_a\sigma_a(T^4-4\pi J).
\end{equation*}
\eqref{eqn:case110012} has three unknown functions, $ J,\rr$ and $Q$. Thanks to the properties of $Q$ in \eqref{eqn:rr},\eqref{eqn:qq}, the zeroth and first moment equations of \eqref{eqn:case110012} have no time derivatives with respect to $Q$. Therefore, instead of solving the RMHD system \eqref{eqn:0012} directly, we take the zeroth and first moments of \eqref{eqn:case110012}, and then couple them together with \eqref{eqn:case120012}--\eqref{eqn:case150012} to solve variables $\rho,J,\rr,\boldsymbol{v},T,\boldsymbol{B}$.
More precisely, we have 
\begin{equation}
\left\{
\begin{aligned}
&\partial_t J+\mathcal{C}\nabla\cdot (\dc\rr)=\f{\mathcal{C}}{4\pi}\hat{S_{re}},\\
&\partial_t (\dc\rr)+\mathcal{C}\dc\nabla J+\mathcal{C}\nabla\cdot {\rm K_Q}=\f{\mathcal{C}}{4\pi}\hat{\boldsymbol{S_{rp}}},
\end{aligned}
\right.
\label{eqn:solsys}
\end{equation}
 coupled with \eqref{eqn:case120012}--\eqref{eqn:case150012} and by treating ${\rm K_Q}$ explicitly, one can update the macroscopic quantities first. Then $Q$
 can be updated implicitly, and the new $I$ can be given  by \eqref{eqn:011}. The details are illustrated in the next two subsections.
\subsection{The time discretization for \eqref{eqn:non}} Let 
$$
t^s=s\Delta t, \qquad s=0,1,2,\cdots,
$$
and 
$$
I^s\approx I(x,t^s,n),\qquad \rho^s\approx \rho(x,t^s),\qquad \boldsymbol{v}^s\approx \boldsymbol{v}(x,t^s),
$$
$$
T^s\approx T(x,t^s),\qquad \boldsymbol{B}^s\approx \boldsymbol{B}(x,t^s).
$$
We use the following semi-discrete scheme in time for  \eqref{eqn:0012}, which reads:
\begin{subequations} \label{eqn:013}
\begin{numcases}{}
\f{J^{s+1}-J^s}{\Delta t}+ \boldsymbol{n}\cdot\frac{\rr^{s+1}-\rr^s}{\Delta t}+\frac{Q^{s+1}-Q^s}{\Delta t} +\mathcal{C}\boldsymbol{n}\cdot\nabla( J^{s+1}+\boldsymbol{n}\cdot\rr^{s+1}+Q^{s})\nonumber\\
\hspace{0.5cm}
 =\mathcal{C}\mathscr{L}_a\sigma_a\left(\f{(T^{s+1})^4}{4\pi}-J^{s+1}-\boldsymbol{n}\cdot\rr^{s+1}- Q^{s+1}\right)-\mathcal{C}\mathscr{L}_s \sigma_s( \boldsymbol{n}\cdot\rr^{s+1}+ Q^{s+1})-\f{2\mathscr{L}_s\sigma_s}{3}\boldsymbol{v}^{s+1}\cdot \rr^{s+1}\nonumber\\
\hspace{0.5cm}+3\mathscr{L}_a\sigma_a\boldsymbol{n}\cdot \boldsymbol{v}^{s+1}\left(\f{(T^{s+1})^4}{4\pi}-J^{s+1} \right)+\boldsymbol{n}\cdot \boldsymbol{v}^{s+1}\left(\mathscr{L}_a\sigma_a+\mathscr{L}_s\sigma_s\right)(4J^{s+1}+\boldsymbol{n}\cdot\rr^{s+1}+Q^{s})\nonumber\\
\hspace{2cm}-\f{\mathscr{L}_a\sigma_a-\mathscr{L}_s\sigma_s}{\mathcal{C}}\left(\boldsymbol{v}^{s+1}\cdot\boldsymbol{v}^{s+1}J^{s+1}+\boldsymbol{v}^{s+1}\cdot\left(\boldsymbol{v}^{s+1}\cdot \dc\right)+\boldsymbol{v}^{s+1}\cdot\left(\boldsymbol{v}^{s+1}\cdot{\rm K_Q^{\it s}}\right)\right),\label{eqna}\\
\frac{\rho^{s+1}-\rho^s}{\Delta t} +\nabla\cdot(\rho^{s} \boldsymbol{v}^{s})=0,\label{eqnb}\\
\frac{(\rho \boldsymbol{v})^{s+1}-(\rho \boldsymbol{v})^s}{\Delta t}+\nabla\cdot\left(\rho^s \boldsymbol{v}^{s}\boldsymbol{v}^{s}-\boldsymbol{B}^{s}\boldsymbol{B}^{s}+(\Ps^*)^{s}\right)=-\PP_0(\hat{\boldsymbol{S_{rp}}})^{s+1},\label{eqnc}\\
 a\frac{(\rho T)^{s+1}-(\rho T)^s}{\Delta t}+\frac{\rho^{s+1}(v^{s+1})^2-\rho^{s}(v^{s})^2}{2\Delta t}+\frac{(B^{s+1})^2-(B^{s})^2}{2\Delta t}\nonumber\\
\hspace{5cm}
+\nabla\cdot\left[(E^{s}+(P^*)^{s})\boldsymbol{v}^{s}-\boldsymbol{B}^{s}(\boldsymbol{B}^{s}\cdot \boldsymbol{v}^{s})\right]=-\mathcal{C}\PP_0(\hat{S_{re}})^{s+1},\label{eqnd}\\
\frac{\boldsymbol{B}^{s+1}-\boldsymbol{B}^s}{\Delta t} +\nabla\times(\boldsymbol{v}^{s}\times \boldsymbol{B}^{s})=0.\label{eqne}\
\end{numcases}
\end{subequations}
where ${\rm K_Q^{\it s}} =\average{\boldsymbol{n}\boldsymbol{n}Q^{s}}$  and 
$$
(\hat{S_{re}})^{s+1}=\mathscr{L}_a\sigma_a((T^{s+1})^4-4\pi J^{s+1})+4\pi(\mathscr{L}_a\sigma_a-\mathscr{L}_s\sigma_s)\f{\boldsymbol{v}^{s+1}}{\mathcal{C}^2}\cdot\left[\mathcal{C}\dc\rr^{s+1}-\left(\boldsymbol{v}^{s+1}J^{s+1}+\dc\boldsymbol{v}^{s+1}J^{s+1}+ {\rm K_Q^{\it s}} \boldsymbol{v}^{s+1}\right)\right]	, \
$$
$$
(\hat{\boldsymbol{S_{rp}}})^{s+1}=-\frac{4\pi(\mathscr{L}_s\sigma_s+\mathscr{L}_a\sigma_a)}{\mathcal{C}}\left[\mathcal{C}\dc\rr^{s+1}-\left(\boldsymbol{v}^{s+1}J^{s+1}+\dc\boldsymbol{v}^{s+1}J^{s+1}+ {\rm K_Q^{\it s}} \boldsymbol{v}^{s+1}\right)\right]+\f{\mathscr{L}_a \boldsymbol{v}^{s+1}}{\mathcal{C}}\sigma_a((T^{s+1})^4-4\pi J^{s+1}).
$$
The zeroth and first moment equations of \eqref{eqna} are

\begin{subequations} \label{eqn:015}
\begin{numcases}{}
\frac{J^{s+1}-J^s}{\Delta t}+\mathcal{C} \nabla\cdot (\dc\rr^{s+1})=\f{\mathcal{C}}{4\pi}(\hat{S_{re}})^{s+1},\label{eqnee1}\\
 \dc\frac{\rr^{s+1}-\rr^s}{\Delta t}+\mathcal{C}\dc\nabla J^{s+1}+\mathcal{C} \nabla\cdot {\rm K_Q^{\it s}}=\f{\mathcal{C}}{4\pi}(\hat{\boldsymbol{S_{rp}}})^{s+1}.\label{eqnfe}\
\end{numcases}
\end{subequations}
We first update $\rho^{s+1}$,$\boldsymbol{B}^{s+1}$ by \eqref{eqnb}, \eqref{eqne} and then solve \eqref{eqnc}, \eqref{eqnd} and \eqref{eqn:015} to get $J^{s+1}$, $\rr^{s+1}$, $\boldsymbol{v}^{s+1}$, $T^{s+1}$. $Q^{s+1}$
 can then be updated implicitly such that
\begin{equation}
\begin{aligned}\label{eqn:imq}
\f{J^{s+1}-J^s}{\Delta t}&+ \boldsymbol{n}\cdot\frac{\rr^{s+1}-\rr^s}{\Delta t}+\frac{Q^{s+1}-Q^s}{\Delta t} +\mathcal{C}\boldsymbol{n}\cdot\nabla( J^{s+1}+\boldsymbol{n}\cdot\rr^{s+1}+Q^{s+1})\\
& =\mathcal{C}\mathscr{L}_a\sigma_a\left(\f{(T^{s+1})^4}{4\pi}-J^{s+1}-\boldsymbol{n}\cdot\rr^{s+1}- Q^{s+1}\right)-\mathcal{C}\mathscr{L}_s \sigma_s( \boldsymbol{n}\cdot\rr^{s+1}+ Q^{s+1})-\f{2\mathscr{L}_s\sigma_s}{3}\boldsymbol{v}^{s+1}\cdot \rr^{s+1}\\
&+3\mathscr{L}_a\sigma_a\boldsymbol{n}\cdot \boldsymbol{v}^{s+1}\left(\f{(T^{s+1})^4}{4\pi}-J^{s+1} \right)+\boldsymbol{n}\cdot \boldsymbol{v}^{s+1}\left(\mathscr{L}_a\sigma_a+\mathscr{L}_s\sigma_s\right)(4J^{s+1}+\boldsymbol{n}\cdot\rr^{s+1}+Q^{s+1})\\
&-\f{\mathscr{L}_a\sigma_a-\mathscr{L}_s\sigma_s}{\mathcal{C}}\left(\boldsymbol{v}^{s+1}\cdot\boldsymbol{v}^{s+1}J^{s+1}+\boldsymbol{v}^{s+1}\cdot\left(\boldsymbol{v}^{s+1}\cdot \dc\right)+\boldsymbol{v}^{s+1}\cdot\left(\boldsymbol{v}^{s+1}\cdot{\rm K_Q^{s}}\right)\right).
\end{aligned}
\end{equation}

In \eqref{eqn:imq}, one has to solve implicitly a linear transport equation for each direction $\boldsymbol{n}$. However, different directions are {\it decoupled} and the computational cost is much lower than solving implicitly a linear RTE where all directions are coupled together.
Finally, we get $I^{s+1}$  by using  \eqref{eqn:011}.
 
\subsection*{The non-equilibrium diffusion limit of (\ref{eqn:013}):}In the non-equilibrium regime, the scalings of the coefficients are
$$
\mathscr{L}_a = \eps,\quad\mathscr{L}_s = 1/\eps,\quad\PP_0=\order(1),\quad\mathcal{C}=c/\eps.
$$ 
Then we will show the AP property of the semi-discretization (\ref{eqn:015}) in the non-equilibrium regime. By using Champan-Enskog expansion in equation \eqref{eqn:imq}, one can get
\begin{equation}
\begin{aligned}\label{eqn:imq-1}
&c\sigma_s( \boldsymbol{n}\cdot\rr^{s+1}+ Q^{s+1})=-\eps^2\left(\f{J^{s+1}-J^s}{\Delta t}+ \boldsymbol{n}\cdot\frac{\rr^{s+1}-\rr^s}{\Delta t}+\f{Q^{s+1}-Q^s}{\Delta t}\right)\\
&-\eps c\boldsymbol{n}\cdot\nabla( J^{s+1}+\boldsymbol{n}\cdot\rr^{s+1}+Q^{s+1}) +\eps^2 c\sigma_a\left(\f{(T^{s+1})^4}{4\pi}-J^{s+1}-\boldsymbol{n}\cdot\rr^{s+1}- Q^{s+1}\right)-\f{2\eps\sigma_s}{3}\boldsymbol{v}^{s+1}\cdot \rr^{s+1}\\
&+3\eps^3\sigma_a\boldsymbol{n}\cdot \boldsymbol{v}^{s+1}\left(\f{(T^{s+1})^4}{4\pi}-J^{s+1} \right)+\boldsymbol{n}\cdot \boldsymbol{v}^{s+1}\left(\eps^3\sigma_a+\eps\sigma_s\right)(4J^{s+1}+\boldsymbol{n}\cdot\rr^{s+1}+Q^{s+1})\\
&-\f{\eps^4\sigma_a-\eps^2\sigma_s}{c}\left(\boldsymbol{v}^{s+1}\cdot\boldsymbol{v}^{s+1}J^{s+1}+\boldsymbol{v}^{s+1}\cdot\left(\boldsymbol{v}^{s+1}\cdot \dc\right)+\boldsymbol{v}^{s+1}\cdot\left(\boldsymbol{v}^{s+1}\cdot{\rm K_Q^{s}}\right)\right),
\end{aligned}
\end{equation}
which indicates
$$
 Q^{s+1}=-\boldsymbol{n}\cdot\rr^{s+1}+\order(\eps).
$$
The residual term $Q$ at current step can be gotten by equation \eqref{eqn:imq}, then one can get
\begin{equation}\label{kq}
{\rm K_Q^{\it s}}=\average{\boldsymbol{n}\boldsymbol{n}Q^s}=\average{\boldsymbol{n}\boldsymbol{n}(-\boldsymbol{n}\cdot\rr^{s}+\order(\eps))}=\order(\eps).
\end{equation}
From \eqref{eqnfe},
\begin{align}\label{eq:rphi1}
\f{1}{\eps}\rr^{s+1}  =   \f{4\boldsymbol{v}^{s+1}J^{s+1}}{c}-\f{1}{\sigma_s}\left(\nabla J^{s+1}-\f{1}{\dc}\nabla\cdot{\rm K_Q^{\it s}}\right)+ \order(\eps) \,.
\end{align}  
by using equation \eqref{kq}, equation \eqref{eq:rphi1} reduces to
\begin{align}\label{eq:rphi}
\f{1}{\eps}\rr^{s+1}  =   \f{4\boldsymbol{v}^{s+1}J^{s+1}}{c}-\f{\nabla J^{s+1}}{\sigma_s}+ \order(\eps)  \,.
\end{align}
Multiplying equation \eqref{eqnee1} by $4\pi\PP_0$, and adding it up with equation \eqref{eqnd}, we get
\begin{equation}\label{limitcase11-1}
\begin{aligned}
& a\frac{(\rho T)^{s+1}-(\rho T)^s}{\Delta t}+\frac{\rho^{s+1}(v^{s+1})^2-\rho^{s}(v^{s})^2}{2\Delta t}+\frac{(B^{s+1})^2-(B^{s})^2}{2\Delta t}+4\pi\PP_0\frac{J^{s+1}-J^s}{\Delta t}\\&\hspace{1.5cm} 
+\nabla\cdot\left[(E^{s}+(P^{*})^s)\boldsymbol{v}^{s}-\boldsymbol{B}^{s}(\boldsymbol{B}^{s}\cdot \boldsymbol{v}^{s})\right]+\f{4\pi c\PP_0}{\eps}\nabla\cdot (\dc\rr^{s+1})=0.\\
\end{aligned}
\end{equation}
By using \eqref{eq:rphi}, and sending $\varepsilon\to 0$ in equation \eqref{limitcase11-1}  yields 
\begin{equation}\label{limitcase11}
\begin{aligned}
&a\frac{(\rho T)^{s+1}-(\rho T)^s}{\Delta t}+\frac{\rho^{s+1}(v^{s+1})^2-\rho^{s}(v^{s})^2}{2\Delta t}+\frac{(B^{s+1})^2-(B^{s})^2}{2\Delta t}+4\pi\PP_0\frac{J^{s+1}-J^s}{\Delta t}\\&\hspace{1.5cm} 
+\nabla\cdot\left[(E^{s}+(P^{*})^s)\boldsymbol{v}^{s}-\boldsymbol{B}^{s}(\boldsymbol{B}^{s}\cdot \boldsymbol{v}^{s})+16\pi\PP_0\dc\boldsymbol{v}^{s+1}J^{s+1}\right]=\nabla\cdot\left(\f{4\pi c\PP_0\dc}{\sigma_s}\nabla J^{s+1}\right).\\
\end{aligned}
\end{equation}
This is a semi-discretization for \eqref{eqn:eqnlimit13}. Then multiplying equation \eqref{eqnfe} by $4\pi\PP_0$, and adding it to equation \eqref{eqnc}, sending $\varepsilon\to 0$, formally one gets
\begin{align}\label{limitcase12}
\frac{(\rho \boldsymbol{v})^{s+1}-(\rho \boldsymbol{v})^s}{\Delta t} +\nabla\cdot\left(\rho^s \boldsymbol{v}^{s}\boldsymbol{v}^{s}-\boldsymbol{B}^{s}\boldsymbol{B}^{s}+(\Ps^*)^{s}\right)=-4\pi\PP_0\dc\nabla J^{s+1}\,.
\end{align}
This is a semi-discretization for \eqref{eqn:eqnlimit12}. Sending  $\varepsilon\to 0$ in   equation   \eqref{eqnee1}, yields 
\begin{align}\label{limitcase13}
4\pi\frac{J^{s+1}-J^s}{\Delta t}+4\pi\dc \nabla\cdot(\boldsymbol{v}^{s+1}J^{s+1})- \nabla\cdot\left(\f{4\pi c\dc}{\sigma_s} \nabla J^{s+1}\right) =c\sigma_a((T^{s+1})^4-4\pi J^{s+1})+4\pi \dc\boldsymbol{v}^{s+1}\cdot\nabla J^{s+1}\,,
\end{align}
which is a semi-discretization \eqref{eqn:eqnlimit14}.
\subsection*{The equilibrium diffusion limit of (\ref{eqn:013}):}In the non-equilibrium regime, the scalings of the coefficients are
$$
\mathscr{L}_a = 1/\eps,\quad\mathscr{L}_s = \eps,\quad\PP_0=\order(1),\quad\mathcal{C}=c/\eps.
$$ 
The derivation of  the equilibrium diffusion limit is similar to the non-equilibrium case, but requires additional  analysis for the zeroth moment of equation \eqref{eqna}.
Sending $\varepsilon\to 0$ in the zeroth moment  equation of \eqref{eqna}, one  gets 
$$
\sigma_a((T^{s+1})^4-4\pi J^{s+1})+\order(\eps)=0,
$$
which implies
$$
4\pi J^{s+1}=(T^{s+1})^4+\order(\eps).
$$
Therefore,  \eqref{limitcase11} and \eqref{limitcase12} in the equilibrium regime become

$$
\begin{aligned}
& a\frac{(\rho T)^{s+1}-(\rho T)^s}{\Delta t}+\frac{\rho^{s+1}(v^{s+1})^2-\rho^{s}(v^{s})^2}{2\Delta t}+\frac{(B^{s+1})^2-(B^{s})^2}{2\Delta t}+\PP_0\frac{(T^{s+1})^4-(T^s)^{4}}{\Delta t}\nonumber\\&\hspace{1.5cm} 
+\nabla\cdot\left[(E^{s}+(P^{*})^s)\boldsymbol{v}^{s}-\boldsymbol{B}^{s}(\boldsymbol{B}^{s}\cdot \boldsymbol{v}^{s})+4\PP_0\dc\boldsymbol{v}^{s+1}(T^{s+1})^4\right]=\nabla\cdot\left(\f{4\pi c\PP_0\dc}{\sigma_a}\nabla (T^{s+1})^4\right),
\end{aligned}
$$
$$
\frac{(\rho \boldsymbol{v})^{s+1}-(\rho \boldsymbol{v})^s}{\Delta t} +\nabla\cdot\left(\rho^s \boldsymbol{v}^{s}\boldsymbol{v}^{s}-\boldsymbol{B}^{s}\boldsymbol{B}^{s}+(\Ps^*)^{s}\right)=-\PP_0\dc\nabla(T^{s+1})^4\,,
$$
which gives a semi-discretization of system \eqref{eqnlimit2}.

\section{The full discretization for the non-equilibrium regime}
For the ease of exposition, we will explain our spatial discretion in 1D. That is, $x \in [0,L]$, $n \in [-1,1]$, and $\average{f(n)} = \frac{1}{2} \int_{-1}^1 f(n) \rd n$. Recall the RMHD equations  in  slab geometry:
\begin{subequations} \label{orig}
\begin{numcases}{}
\partial_t (J+nR+Q)+\mathcal{C}n\partial_x(J+nR+Q)=
\mathcal{C}\mathscr{L}_a\sigma_a\left(\f{T^4}{4\pi}-J-nR-Q \right)-\mathcal{C}\mathscr{L}_s\sigma_s(nR+Q)+G,\label{orig1}\\
\partial_t\rho +\partial_x(\rho v_x)=0,\label{orig2}\\
 \partial_t(\rho v_x)+\partial_x\left(\rho v_x^2+p+B^2/2+B_x^2\right)=-\PP_0S_{rp},\label{orig3}\\
  \partial_t(\rho v_y)+\partial_x\left(\rho v_xv_y-B_xB_y\right)=0,\label{orig4}\\
  \partial_t(\rho v_z)+\partial_x\left(\rho v_xv_z-B_xB_z\right)=0,\label{orig5}\\
\partial_t \left(a\rho T+\f{\rho v^2+B^2}{2}\right)+\partial_x[(E+P^*)v_x-(\boldsymbol{B}\cdot \boldsymbol{v})B_x]=-\mathcal{C}\PP_0S_{re},\label{orig6}\\
\partial_t B_y+\partial_x( B_yv_x-B_xv_y)=0,\label{orig7}\\
\partial_t B_z+\partial_x( B_zv_x-B_xv_z)=0,\label{orig8}
 \end{numcases}
\end{subequations}
where
$$
\begin{aligned}
G=3\mathscr{L}_a\sigma_a nv_{x}\left(\f{T^4}{4\pi}-J\right)+n v_{x}\left(\mathscr{L}_a\sigma_{a}+\mathscr{L}_s\sigma_{s}\right)\left(4J+nR+Q\right)-\f{2\mathscr{L}_s\sigma_{s} }{3}v_{x} R-\f{(\mathscr{L}_a\sigma_{a}-\mathscr{L}_s\sigma_{s})(v_{x})^2}{\mathcal{C}}\left(\f{4}{3}J+ K_{Q}\right),
\end{aligned}
$$
$$
S_{re}=4\pi\mathscr{L}_a\sigma_{a}\left(\f{T^4}{4\pi}-J\right)+\f{4\pi(\mathscr{L}_a\sigma_{a}-\mathscr{L}_s\sigma_{s})v_{x}}{\mathcal{C}^2}\left[\f{\mathcal{C}}{3}R-v_{x}\left(\f{4}{3}J+ K_{Q}\right)\right]	, \
$$
$$
S_{rp}=-\frac{4\pi(\mathscr{L}_s\sigma_{s}+\mathscr{L}_a\sigma_{a})}{\mathcal{C}}\left[\f{\mathcal{C}}{3}R-v_{x}\left(\f{4}{3}J+ K_{Q}\right)\right]+\f{4\pi\mathscr{L}_a\sigma_{a} v_{x}}{\mathcal{C}}\left(\f{T^4}{4\pi}-J\right),
$$ and $v_x, v_y$ and $v_z$ are fluid velocity in the $x-$, $y-$, and $z-$directions, respectively. $B_x, B_y$ and $B_z$ are magnetic field in the $x-$, $y-$, and $z-$directions, respectively. The boundary condition becomes
\begin{equation} \label{BC1d}
I(t,0,n) = b_\text{L}(t,n), ~ \text{ for } ~n>0; \qquad I(t,L,n) = b_\text{R}(t,n), ~ \text{ for } ~ n<0\,.
\end{equation}
Higher dimensions can be treated in the dimension by dimension manner. Let $[a,b] $ be the computational domain, $\Delta x=(b-a)/N_x$, and we consider the uniform mesh as follows
\[
x_{i-\f{1}{2}}=a+(i-1)\Delta x, \quad i = 1, 2 \cdots, N_{x}+1,
\]
and let
$$x_{\f{1}{2}}=a<x_1<x_{\f{3}{2}}<\cdots<x_{i-\f{1}{2}}<x_i<x_{i+\f{1}{2}}<\cdots<x_{N_x}<x_{N_x+\f{1}{2}}=b,$$
$$x_i=\left(x_{i-\f{1}{2}}+x_{i+\f{1}{2}}\right)/2,\qquad \mbox{for $i=1,\cdots, N_x.$}$$

To get a consistent stencil in spatial discretization, we use the unified gas kinetic scheme (UGKS) for spatial discretization \cite{mieussens2013asymptotic}. Other space discretization that is AP can be applied as well. The most crucial point is that the space discretization of the MHD system has to be consistent with that of the RTE. 
\subsection{A finite volume approach}

UGKS is a finite volume method. Integrating the RMHD system \eqref{orig}  over 
$[t_s,t_{s+1}]$ and  $[x_{i-\f{1}{2}},x_{i+\f{1}{2}}]$ gives
\begin{subequations} \label{fvm}
	\begin{numcases}{}
\f{J_i^{s+1}-J_i^s}{\Delta t}+n\frac{R_i^{s+1}-R_i^s}{\Delta t}+\frac{Q_i^{s+1}-Q_i^s}{\Delta t}+\f{1}{\Delta x}(\zeta_{i+\f{1}{2}}-\zeta_{i-\f{1}{2}})=\nonumber\\
\hspace{2.5cm}\mathcal{C}\mathscr{L}_a\sigma_a\left(\f{(T_i^{s+1})^4}{4\pi}-J_i^{s+1}-nR_i^{s+1}-Q_i^{s+1} \right)-\mathcal{C}\mathscr{L}_s\sigma_s(nR_i^{s+1}+Q_i^{s+1})+G_i,\label{fvm1}\\
	\f{\rho_i^{s+1}-\rho_i^{s}}{\Delta t} +\f{1}{\Delta x}(F^s_{1,i+\f{1}{2}}-F^s_{1,i-\f{1}{2}})=0,\label{fvm2}\\
	\f{(\rho v_x)_i^{s+1}-(\rho v_x)_i^{s}}{\Delta t}+\f{1}{\Delta x}(F^s_{2,i+\f{1}{2}}-F^s_{2,i-\f{1}{2}})=-\PP_0(\hat{S_{rp}})_i^{s+1},\label{fvm3}\\
		\f{(\rho v_y)_i^{s+1}-(\rho v_y)_i^{s}}{\Delta t}+\f{1}{\Delta x}(F^s_{3,i+\f{1}{2}}-F^s_{3,i-\f{1}{2}})=0,\label{fvm4}\\
		\f{(\rho v_z)_i^{s+1}-(\rho v_z)_i^{s}}{\Delta t}+\f{1}{\Delta x}(F^s_{4,i+\f{1}{2}}-F^s_{4,i-\f{1}{2}})=0,\label{fvm5}\\
	a\f{(\rho T)_i^{s+1}-(\rho T)_i^{s}}{\Delta t}+	\f{(\rho v^2+B^2)_i^{s+1}-(\rho v^2+B^2)_i^{s}}{2\Delta t}+\f{1}{\Delta x}(F^s_{5,i+\f{1}{2}}-F^s_{5,i-\f{1}{2}})
	=-\mathcal{C}\PP_0(\hat{S_{re}})_i^{s+1},\label{fvm6}\\
	\f{B_{y,i}^{s+1}-B_{y,i}^{s}}{\Delta t}+\f{1}{\Delta x}(F^s_{6,i+\f{1}{2}}-F^s_{6,i-\f{1}{2}})=0,\label{fvm7}\\
	\f{B_{z,i}^{s+1}-B_{z,i}^{s}}{\Delta t}+\f{1}{\Delta x}(F^s_{7,i+\f{1}{2}}-F^s_{7,i-\f{1}{2}})=0,\label{fvm8}
	\end{numcases}
\end{subequations}
where  $G_i$, $(\hat{S_{re}})_i^{s+1}$ and $(\hat{S_{rp}})_i^{s+1}$  are 
$$
\begin{aligned}
G_i=&3\mathscr{L}_a\sigma_{a,i}nv_{x,i}^{s+1}\left(\f{(T_i^{s+1})^4}{4\pi}-J_{i}^{s+1}\right)+n v_{x,i}^{s+1}\left(\mathscr{L}_a\sigma_{a,i}+\mathscr{L}_s\sigma_{s,i}\right)\left(4J^{s+1}_{i}+nR_{i}^{s+1}+Q_{i}^{s}\right)\\
&\hspace{2cm}-\f{2}{3}\mathscr{L}_s\sigma_{s,i} v_{x,i}^{s+1} R_{i}^{s+1}-\f{(\mathscr{L}_a\sigma_{a,i}-\mathscr{L}_s\sigma_{s,i})(v_{x,i}^{s+1})^2}{\mathcal{C}}\left(\f{4}{3}J^{s+1}_{i}+ K_{Q,i}^{s}\right),
\end{aligned}
$$
$$
(\hat{S_{re}})_i^{s+1}=4\pi\mathscr{L}_a\sigma_{a,i}\left(\f{(T_i^{s+1})^4}{4\pi}-J_{i}^{s+1}\right)+\f{4\pi(\mathscr{L}_a\sigma_{a,i}-\mathscr{L}_s\sigma_{s,i})v_{x,i}^{s+1}}{\mathcal{C}^2}\left[\f{\mathcal{C}}{3}R_i^{s+1}-v_{x,i}^{s+1}\left(\f{4}{3}J^{s+1}_i+K_{Q,i}^{s}\right)\right]	, \
$$
$$
(\hat{S_{rp}})_i^{s+1}=-\frac{4\pi(\mathscr{L}_s\sigma_{s,i}+\mathscr{L}_a\sigma_{a,i})}{\mathcal{C}}\left[\f{\mathcal{C}}{3}R_i^{s+1}-v_{x,i}^{s+1}\left(\f{4}{3}J^{s+1}_i+ K_{Q,i}^{s}\right)\right]+\f{4\pi\mathscr{L}_a\sigma_{a,i} v_{x,i}^{s+1}}{\mathcal{C}}\left(\f{(T_i^{s+1})^4}{4\pi}-J_{i}^{s+1}\right).
$$
Here  the numerical fluxes $F^s_{\cdot,i\pm\f{1}{2}}$ at the interfaces use Roe's  Riemann solver \cite{stone2008athena}, based on the piecewise linear reconstruction. At last, we discuss  the microscopic flux $\zeta$ defined at the interface $x_{i+\f{1}{2}}$ in \eqref{fvm1}. In  \cite{mieussens2013asymptotic}, all the terms of  the microscopic flux are treated explicitly, in this paper all the terms except the initial value and $K_Q$ are treated $\it{implicitly}$, and the processes of UGKS are detailed in Appendix \ref{apen1}. 


With the expressions \eqref{isol}, \eqref{equ:i0} and \eqref{equ:j}, the numerical flux $\zeta_{i+\f{1}{2}}$ is
\begin{equation}\label{zsol}
\begin{aligned}
\zeta_{i+\f{1}{2}}&=A_{i+\frac{1}{2}} n\left(I_{i}^{s} \mathds{1}_{n>0}+I_{i+1}^{s} \mathds{1}_{n<0}\right)+C^1_{i+\frac{1}{2}} n J_{i+\frac{1}{2}}^{s+1}+\f{C^2_{i+\frac{1}{2}}}{4\pi} n (T_{i+\frac{1}{2}}^{s+1})^4+F_{i+\frac{1}{2}} n \hat{G}_{i+\f{1}{2}}\\
&+D^1_{i+\frac{1}{2}} n^{2}\left(\delta_{x} J_{i+\frac{1}{2}}^{s+1,+} \mathds{1}_{n>0}+\delta_{x} J_{i+\frac{1}{2}}^{s+1,-} \mathds{1}_{n<0}\right)+\f{D^2_{i+\frac{1}{2}}}{4\pi} n^{2}\left(\delta_{x} (T_{i+\frac{1}{2}}^{s+1,+})^4 \mathds{1}_{n>0}+\delta_{x} (T_{i+\frac{1}{2}}^{s+1,-})^4 \mathds{1}_{n<0}\right),
\end{aligned}
\end{equation}
where $\mathds 1_{n\lessgtr0}$ is the indicator function. 
Moreover, the coefficients in the numerical flux are given by:
\begin{equation}\label{coeff}
\begin{aligned}
&A=\f{\mathcal{C}}{\Delta t\mu}(1-e^{-\mu\Delta t}),\\
&C^1=\f{\mathcal{C}^2\mathscr{L}_s\sigma_s}{\Delta t\mu}\left(\Delta t-\f{1}{\mu}(1-e^{-\mu\Delta t})\right),\qquad C^2=\f{\mathcal{C}^2\mathscr{L}_a\sigma_a}{\Delta t\mu}\left(\Delta t-\f{1}{\mu}(1-e^{-\mu\Delta t})\right),\\
&D^1=-\f{\mathcal{C}^3\mathscr{L}_s\sigma_s}{\Delta t\mu^2}\left(\Delta t\left(1+e^{-\mu \Delta t}\right)-\frac{2}{\mu}\left(1-e^{-\mu \Delta t}\right)\right),\\
& D^2=-\f{\mathcal{C}^3\mathscr{L}_a\sigma_a}{\Delta t\mu^2}\left(\Delta t\left(1+e^{-\mu \Delta t}\right)-\frac{2}{\mu}\left(1-e^{-\mu \Delta t}\right)\right),\\
&F=\frac{\mathcal{C}}{\Delta t  \mu}\left(\Delta t-\frac{1}{\mu}\left(1-e^{-\mu \Delta t}\right)\right),
\end{aligned}
\end{equation}
with $\mu=\mathcal{C}\mathscr{L}_a\sigma_a+ \mathcal{C}\mathscr{L}_s\sigma_s$.

To obtain a scheme that can update the quantities $J$,$R$, $T$ and $v_x$, the zeroth moment and first moment of velocity field $n$ for  \eqref{fvm1} are needed. Integrating  equation \eqref{fvm1} for $n$ from -1 to 1, and multiplying the both sides by $2\pi$, one can obtain
\begin{equation}\label{zerom}
4\pi\f{J_i^{s+1}-J_i^{s}}{\Delta t}+\f{2\pi}{\Delta x}\int_{-1}^1\hat{\zeta}_{i+\f{1}{2}}-\hat{\zeta}_{i-\f{1}{2}}dn=\mathcal{C}(\hat{S_{re}})_i^{s+1},
\end{equation}
where
$$
\begin{aligned}
&\int_{-1}^1\hat{\zeta}_{i+\f{1}{2}}dn=\int_{-1}^1A_{i+\frac{1}{2}} n\left(I_{i}^{s} \mathds{1}_{n>0}+I_{i+1}^{s} \mathds{1}_{n<0}\right)dn+\f{2D^1_{i+\frac{1}{2}}}{3}\f{J_{i+1}^{s+1}-J_{i}^{s+1}}{\Delta x}+\f{D^2_{i+\frac{1}{2}}}{6\pi}\f{(T_{i+1}^{s+1})^4-(T_{i}^{s+1})^4}{\Delta x} \\
&+F_{i+\frac{1}{2}} \left[ 2\mathscr{L}_a\sigma_{a,i+\f{1}{2}}v_{x,i+\f{1}{2}}^{s+1}\left(\f{(T_{i+\f{1}{2}}^{s+1})^4}{4\pi}-J_{i+\f{1}{2}}^{s+1}\right)+v_{x,i+\f{1}{2}}^{s+1}\left(\mathscr{L}_s\sigma_{s,i+\f{1}{2}}+\mathscr{L}_a\sigma_{a,i+\f{1}{2}}\right)\left(\f{8}{3}J^{s+1}_{i+\f{1}{2}}+2 K_{Q,i+\f{1}{2}}^{s}\right)\right].
\end{aligned}
$$
Multiplying equation \eqref{fvm1} by $n$, then integrating it for $n$  from -1 to 1, and multiplying the both sides by $2\pi$, one can obtain
\begin{equation}\label{firm}
\f{4\pi}{3}\f{R_i^{s+1}-R_i^{s}}{\Delta t}+\f{2\pi}{\Delta x}\int_{-1}^1n\left(\hat{\zeta}_{i+\f{1}{2}}-\hat{\zeta}_{i-\f{1}{2}}\right)dn=\mathcal{C}(\hat{S_{rp}})_i^{s+1},
\end{equation}
where
$$
\begin{aligned}
\int_{-1}^1n\hat{\zeta}_{i+\f{1}{2}}dn&=\int_{-1}^1A_{i+\frac{1}{2}} n^2\left(I_{i}^{s} \mathds{1}_{n>0}+I_{i+1}^{s} \mathds{1}_{n<0}\right)dn+\f{C^1_{i+\frac{1}{2}}}{3} (J_{i+1}^{s+1}+J_{i}^{s+1})+\f{C^2_{i+\frac{1}{2}}}{12\pi} ((T_{i+1}^{s+1})^4+(T_{i}^{s+1})^4)\\
&+F_{i+\frac{1}{2}}\left[v_{x,i+\f{1}{2}}^{s+1}\left(\mathscr{L}_a\sigma_{a,{i+\f{1}{2}}}+\mathscr{L}_s\sigma_{s,{i+\f{1}{2}}}\right)\left(\f{2}{5}R_{i+\f{1}{2}}^{s+1}+\int_{-1}^1n^3Q_{i+\f{1}{2}}^sdn\right)-\f{4}{9}\mathscr{L}_s\sigma_{s,{i+\f{1}{2}}}v_{x,i+\f{1}{2}}^{s+1}R_{i+\f{1}{2}}^{s+1} \right.\\
\\&\left.-\f{2(\mathscr{L}_a\sigma_{a,{i+\f{1}{2}}}-\mathscr{L}_s\sigma_{s,{i+\f{1}{2}}})(v_{x,i+\f{1}{2}}^{s+1})^2}{3\mathcal{C}}\left(\f{4}{3}J^{s+1}_{i+\f{1}{2}}+ K_{Q,{i+\f{1}{2}}}^{s}\right)\right].
\end{aligned}
$$
Coupling  equations \eqref{zerom}, \eqref{firm}, \eqref{fvm3} and \eqref{fvm6} as a system and solving the system,  quantities $J$, $T$, $R$ and $v_x$ can be updated. Finally,  quantity $Q$ can be updated implicitly by  equation \eqref{fvm1}, 
which reads
$$
\begin{aligned}
	\f{J_i^{s+1}-J_i^s}{\Delta t}+&n\frac{R_i^{s+1}-R_i^s}{\Delta t}+\frac{Q_i^{s+1}-Q_i^s}{\Delta t}+\f{1}{\Delta x}(\zeta_{i+\f{1}{2}}-\zeta_{i-\f{1}{2}})=\\
&\mathcal{C}\mathscr{L}_a\sigma_a\left(\f{(T_i^{s+1})^4}{4\pi}-J_i^{s+1}-nR_i^{s+1}-Q_i^{s+1} \right)-\mathcal{C}\mathscr{L}_s\sigma_s(nR_i^{s+1}+Q_i^{s+1})+G_i,
\end{aligned}
$$
where
\begin{equation*}
\begin{aligned}
\zeta_{i+\f{1}{2}}&=A_{i+\frac{1}{2}} n\left((J_{i}^{s+1}+nR_{i}^{s+1}+Q_{i}^{s+1}) \mathds{1}_{n>0}+(J_{i+1}^{s+1}+nR_{i+1}^{s+1}+Q_{i+1}^{s+1}) \mathds{1}_{n<0}\right)\\
&+C^1_{i+\frac{1}{2}} n J_{i+\frac{1}{2}}^{s+1}+\f{C^2_{i+\frac{1}{2}}}{4\pi} n (T_{i+\frac{1}{2}}^{s+1})^4+F_{i+\frac{1}{2}} n \hat{G}_{i+\f{1}{2}}\\
&+D^1_{i+\frac{1}{2}} n^{2}\left(\delta_{x} J_{i+\frac{1}{2}}^{s+1,+} \mathds{1}_{n>0}+\delta_{x} J_{i+\frac{1}{2}}^{s+1,-} \mathds{1}_{n<0}\right)+\f{D^2_{i+\frac{1}{2}}}{4\pi} n^{2}\left(\delta_{x} (T_{i+\frac{1}{2}}^{s+1,+})^4 \mathds{1}_{n>0}+\delta_{x} (T_{i+\frac{1}{2}}^{s+1,-})^4 \mathds{1}_{n<0}\right),
\end{aligned}
\end{equation*}
and 
$$
\begin{aligned}
G_i=&3\mathscr{L}_a\sigma_{a,i}nv_{x,i}^{s+1}\left(\f{(T_i^{s+1})^4}{4\pi}-J_{i}^{s+1}\right)+n v_{x,i}^{s+1}\left(\mathscr{L}_a\sigma_{a,i}+\mathscr{L}_s\sigma_{s,i}\right)\left(4J^{s+1}_{i}+nR_{i}^{s+1}+Q_{i}^{s+1}\right)\\
&\hspace{2cm}-\f{2}{3}\mathscr{L}_s\sigma_{s,i} v_{x,i}^{s+1} R_{i}^{s+1}-\f{(\mathscr{L}_a\sigma_{a,i}-\mathscr{L}_s\sigma_{s,i})(v_{x,i}^{s+1})^2}{\mathcal{C}}\left(\f{4}{3}J^{s+1}_{i}+ K_{Q,i}^{s}\right),
\end{aligned}
$$
and the new $I$ can be given by using   \eqref{eqn:011}.


	\subsection{The diffusion limit of \eqref{fvm3}--\eqref{firm}}Firstly, the non-equilibrium regime is considered, i.e.,
	$$
\mathscr{L}_a = \eps,\quad\mathscr{L}_s = 1/\eps,\quad\PP_0=\order(1),\quad\mathcal{C}=c/\eps.
$$ 
Assuming $\sigma_a$ and $\sigma_s$ are positive, as $\eps\to 0$, the leading order of the coefficients in the numerical flux \eqref{zsol}
are
$$
\begin{aligned}
&A(\Delta t,\eps, c, \sigma_a, \sigma_s)=0+\order(\eps),\\
&\eps C^1(\Delta t,\eps, c, \sigma_a, \sigma_s)=c+\order(\eps),\qquad \eps C^2=0+\order(\eps),\\
&D^1=-\f{c}{\sigma_s}+\order(\eps),\qquad D^2=0+\order(\eps),\\
&\f{1}{\eps}F=\f{1}{\sigma_s}+\order(\eps).
\end{aligned}
$$
This means the zero moment and the first moment of the numerical flux, as $\eps\to 0$,  have the following limit:
\begin{equation}\label{limz} 
\int_{-1}^1\hat{\zeta}_{i+\f{1}{2}}dn\to\frac{8}{3}v_{x,i+\f{1}{2}}^{s+1}J_{i+\f{1}{2}}^{s+1}-\f{2c}{3\sigma_{s,i+\frac{1}{2}}}\f{J_{i+1}^{s+1}-J_{i}^{s+1}}{\Delta x},
\end{equation}
\begin{equation}\label{limfi}
\eps\int_{-1}^1n\hat{\zeta}_{i+\f{1}{2}}dn\to\f{c}{3}  (J_{i+1}^{s+1}+J_{i}^{s+1}),
\end{equation}
in which the condition the residual term $Q$ is the same order with $\eps$ in the non-equilibrium regime has been used.

Multiplying equation \eqref{zerom} by $\PP_0$, and adding it up with equation \eqref{fvm6}, the terms on the right hand sides cancel.   Sending $\eps\to 0$ and using \eqref{limz} yields 
\begin{equation}\label{tolenl}
\begin{aligned}
&a\f{(\rho T)_i^{s+1}-(\rho T)_i^{s}}{\Delta t}+	\f{(\rho v^2+B^2)_i^{s+1}-(\rho v^2+B^2)_i^{s}}{2\Delta t}+4\pi\PP_0\f{J_i^{s+1}-J_i^{s}}{\Delta t}+\f{1}{\Delta x}(F^s_{5,i+\f{1}{2}}-F^s_{5,i-\f{1}{2}})\\
&\hspace{1cm}+\f{16\pi\PP_0 }{3}\f{v_{x,i+\f{1}{2}}^{s+1}J_{i+\f{1}{2}}^{s+1}-v_{x,i-\f{1}{2}}^{s+1}J_{i-\f{1}{2}}^{s+1}}{\Delta x}=\f{4\pi\PP_0c}{3\sigma_{s,i+\frac{1}{2}}\Delta x}\f{J_{i+1}^{s+1}-J_{i}^{s+1}}{\Delta x}-\f{4\pi\PP_0c}{3\sigma_{s,i-\frac{1}{2}}\Delta x}\f{J_{i}^{s+1}-J_{i-1}^{s+1}}{\Delta x}.
\end{aligned}
\end{equation}
This is a full discretization for \eqref{eqn:eqnlimit13}. Then multiplying equation \eqref{firm} by $\f{\eps}{c}\PP_0$, and adding it to equation \eqref{fvm3}, sending $\eps\to 0$ and using \eqref{limfi}, one  obtains
\begin{equation}\label{tolmol}
\begin{aligned}
&\f{(\rho v_x)_i^{s+1}-(\rho v_x)_i^{s}}{\Delta t}+\f{1}{\Delta x}(F^s_{2,i+\f{1}{2}}-F^s_{2,i-\f{1}{2}})+\f{2\pi\PP_0 }{3}\f{J_{i+1}^{s+1}-J_{i-1}^{s+1}}{\Delta x}=0.
\end{aligned}
\end{equation}
This is a full discretization for \eqref{eqn:eqnlimit12}. Sending $\eps\to 0$ in  equation \eqref{firm} and using  \eqref{limfi} gives
\begin{align}\label{eq:rj}
\f{1}{\eps}R_i^{s+1}=\f{4}{c}v_{x,i}^{s+1}J_i^{s+1}-\f{J_{i+1}^{s+1}-J_{i-1}^{s+1}}{2\Dx\sigma_{s,i}},
\end{align}
using \eqref{eq:rj} in \eqref{zerom} and sending $\eps\to 0$, one  gets 
\begin{align}
&4\pi\frac{J_i^{s+1}-J_i^{s}}{\Delta t}-\frac{4\pi}{\Dx}\left[ \frac{c}{3\sigma_{s,i+\f{1}{2}}} \frac{ J_{i+1}^{s+1} - J_{i}^{s+1}}{\Dx } - \frac{c}{3\sigma_{s,i-\f{1}{2}}} \frac{ J_{i}^{s+1} - J_{i-1}^{s+1}}{\Dx } \right]\nonumber\\
&\hspace{1.5cm} 
+\f{16\pi }{3}\f{v_{x,i+\f{1}{2}}^{s+1}J_{i+\f{1}{2}}^{s+1}-v_{x,i-\f{1}{2}}^{s+1}J_{i-\f{1}{2}}^{s+1}}{\Delta x}=c\sigma_{a,i}\left((T_i^{s+1})^4-4\pi J_i^{s+1}\right)-\f{2\pi}{3}v_{x,i}^{s+1} \f{J^{s+1}_{i+1}-J^{s+1}_{i-1}}{\Dx},
\end{align}
which is a full discretization for \eqref{eqn:eqnlimit14}.
Next, we will give the diffusion limit of the full discretization of the equilibrium case.
Firstly, assuming $\sigma_a$ and $\sigma_s$ are positive, as $\eps\to 0$, the leading order of the coefficients in the numerical flux \eqref{zsol}
are
$$
\begin{aligned}
&A(\Delta t,\eps, c, \sigma_a, \sigma_s)=0+\order(\eps),\\
&\eps C^1(\Delta t,\eps, c, \sigma_a, \sigma_s)=0+\order(\eps),\qquad \eps C^2=c+\order(\eps),\\
&D^1=0+\order(\eps),\qquad D^2=-\f{c}{\sigma_a}+\order(\eps),\\
&\f{1}{\eps}F=\f{1}{\sigma_a}+\order(\eps).
\end{aligned}
$$
which means the zeroth moment and the first moment of the numerical flux, as $\eps\to 0$,  have the following limit:
\begin{equation}\label{limzcase2} 
\int_{-1}^1\hat{\zeta}_{i+\f{1}{2}}dn\to\frac{8}{3}v_{x,i+\f{1}{2}}^{s+1}J_{i+\f{1}{2}}^{s+1}-\f{2c}{3\sigma_{a,i+\frac{1}{2}}}\f{J_{i+1}^{s+1}-J_{i}^{s+1}}{\Delta x},
\end{equation}
\begin{equation}\label{limficase2}
\eps\int_{-1}^1n\hat{\zeta}_{i+\f{1}{2}}dn\to\f{c}{3}  (J_{i+1}^{s+1}+J_{i}^{s+1}).
\end{equation}
in which the condition the residual term $Q$ is the same order as $\eps$ in the equilibrium regime has been used.
In equation \eqref{zerom}, letting $\eps\to 0$ , one can obtain
$$
4\pi\sigma_{a,i}\left(\f{(T_i^{s+1})^4}{4\pi}-J_i^{s+1}\right)=0.
$$
Multiplying equation \eqref{zerom} by $\PP_0$, adding it to the energy equation  and substituting  limit \eqref{limzcase2} into the equation, then
multiplying equation \eqref{firm} by $\f{\eps}{c}\PP_0$, adding it to the momentum equation and substituting  limit \eqref{limficase2} into the equation, at last letting $\eps\to 0$ in the two equations, one can obtain
\begin{equation}\label{tolenlcase1}
\begin{aligned}
&a\f{(\rho T)_i^{s+1}-(\rho T)_i^{s}}{\Delta t}+	\f{(\rho v^2+B^2)_i^{s+1}-(\rho v^2+B^2)_i^{s}}{2\Delta t}+\PP_0\f{(T_i^{s+1})^4-(T_i^{s})^4}{\Delta t}+\f{1}{\Delta x}(F^s_{3,i+\f{1}{2}}-F^s_{3,i-\f{1}{2}})\\
&\hspace{1cm}+\f{4\PP_0 }{3}\f{v_{x,i+\f{1}{2}}^{s+1}(T_{i+\f{1}{2}}^{s+1})^4-v_{x,i-\f{1}{2}}^{s+1}(T_{i-\f{1}{2}}^{s+1})^4}{\Delta x}=\f{\PP_0c}{3\sigma_{a,i+\frac{1}{2}}\Delta x}\f{(T_{i+1}^{s+1})^4-(T_{i}^{s+1})^4}{\Delta x}-\f{\PP_0c}{3\sigma_{a,i-\frac{1}{2}}\Delta x}\f{(T_{i}^{s+1})^4-(T_{i-1}^{s+1})^4}{\Delta x},
\end{aligned}
\end{equation}
\begin{equation}\label{tolmolcase1}
\begin{aligned}
&\f{(\rho v_x)_i^{s+1}-(\rho v_x)_i^{s}}{\Delta t}+\f{1}{\Delta x}(F^s_{2,i+\f{1}{2}}-F^s_{2,i-\f{1}{2}})+\f{\PP_0 }{6}\f{(T_{i+1}^{s+1})^4-(T_{i-1}^{s+1})^4}{\Delta x}=0,
\end{aligned}
\end{equation}
which are full discretization  for \eqref{eqn:eqnlimit23} and \eqref{eqn:eqnlimit22}.
	\subsection{Boundary conditions} In this section, the numerical fluxes in the system  \eqref{fvm3}--\eqref{firm}  under the boundary condition \eqref{BC1d} are considered. The construction of the flux on the right boundary  is  similar to the construction  on the left boundary, for which we only consider the left boundary case. At the left boundary, the integral  representation of  the radiative intensity $I$ \eqref{isol} reads:
\begin{equation}\label{equ:inte}
I_{\f{1}{2}}(t) =\begin{cases}
b_L,&{\rm if} \quad n>0,\\
e^{-\mu_{\f{1}{2}}(t-t_s)}I\left(t_s,x_{\f{1}{2}}+\mathcal{C}n(t-t_s)\right)+\f{1-e^{-\mu_{\f{1}{2}}(t-t_s)}}{\mu_{\f{1}{2}}}\hat{G}_{\f{1}{2}}\\
+\int\limits^t_{t_s}e^{-\mu_{\f{1}{2}}(t-t_s)}\left(\mathcal{C}\mathscr{L}_s\sigma_{s,\f{1}{2}}J\left(z,x_{\f{1}{2}}-\mathcal{C}n(t-z)\right)+\f{\mathcal{C}\mathscr{L}_a\sigma_{a,\f{1}{2}}}{4\pi}T^4\left(z,x_{\f{1}{2}}-\mathcal{C}n(t-z)\right)\right)dz, &{\rm if} \quad n<0.
\end{cases}
\end{equation} 	
According to the approximation \eqref{equ:i0} and \eqref{equ:j} , 	the reconstruction of $I$ for $n<0$ and $x>x_{\f{1}{2}}$ at the left boundary can be written as:
$$
I(t_s,x,n)=I^s_{1},
$$	
and the reconstruction of $J$ for $t$ in the interval $[t_s,t_{s+1}]$:
$$
J(t,x)=J^{s+1}_{\f{1}{2}}+\delta_xJ^{s+1,-}_{\f{1}{2}}(x-x_{\f{1}{2}}),
$$
where $J^{s+1}_{\f{1}{2}}=\f{1}{2}(\average{b_L}+J^{s+1}_{1})$ and $\delta_xJ^{s+1,-}_{\f{1}{2}}=\f{J^{s+1}_{1}-J^{s+1}_{\f{1}{2}}}{\Delta x/2}$. Then the numerical flux at the left boundary reads:
$$
\begin{aligned}
\zeta_{\f{1}{2}}&=\mathcal{C}nb_L\mathds{1}_{n>0}+A_{\frac{1}{2}} nI_{1}^{s} \mathds{1}_{n<0}+C^1_{\frac{1}{2}} n J_{\frac{1}{2}}^{s+1}\mathds{1}_{n<0}+\f{C^2_{\frac{1}{2}}}{4\pi} n (T_{\frac{1}{2}}^{s+1})^4\mathds{1}_{n<0}+F_{\frac{1}{2}} n \hat{G}_{\frac{1}{2}}\mathds{1}_{n<0}\\
&+D^1_{\frac{1}{2}} n^{2}\delta_{x} J_{\frac{1}{2}}^{s+1,-} \mathds{1}_{n<0}+\f{D^2_{\frac{1}{2}}}{4\pi} n^{2}\delta_{x} (T_{\frac{1}{2}}^{s+1,-})^4 \mathds{1}_{n<0},
\end{aligned}
$$
which indicates
$$
\begin{aligned}
\int_{-1}^1\hat{\zeta}_{\f{1}{2}}dn&=\int_{-1}^1\left(\mathcal{C}nb_L\mathds{1}_{n>0}+A_{\frac{1}{2}} nI_{1}^{s} \mathds{1}_{n<0}\right)dn
-\f{1}{2}\left( \f{C_{\f{1}{2}}^1}{2}\left(J_1^{s+1}+\average{b_L}\right)+\f{C^2_{\frac{1}{2}}} {8\pi} ((T _1^{s+1})^4+(T_L^{s+1})^4)\right)\\
&+\int_{-1}^1F_{\frac{1}{2}} n \hat{G}_{\frac{1}{2}}\mathds{1}_{n<0}dn+\f{1}{3\Delta x}\left(D^1_{\frac{1}{2}} (J_1^{s+1}-\average{b_L})+\f{D^2_{\frac{1}{2}}}{4\pi}( (T _1^{s+1})^4-(T_L^{s+1})^4)\right),
\end{aligned}
$$
here $T_L^{s+1}$ is the material temperature at the left boundary, and 
$$
\begin{aligned}
\int_{-1}^1n\hat{\zeta}_{\f{1}{2}}dn&=\int_{-1}^1\left(\mathcal{C}n^2b_L\mathds{1}_{n>0}+A_{\frac{1}{2}} n^2I_{1}^{s} \mathds{1}_{n<0}\right)dn+\f{1}{3}\left( \f{C_{\f{1}{2}}^1}{2}\left(J_1^{s+1}+\average{b_L}\right)+\f{C^2_{\frac{1}{2}}} {8\pi}((T _1^{s+1})^4+(T_L^{s+1})^4)\right)\\
&+\int_{-1}^1F_{\frac{1}{2}} n^2 \hat{G}_{\frac{1}{2}}\mathds{1}_{n<0}dn-\f{1}{4\Delta x}\left(D^1_{\frac{1}{2}} (J_1^{s+1}-\average{b_L})+\f{D^2_{\frac{1}{2}}}{4\pi}( (T _1^{s+1})^4-(T_L^{s+1})^4)\right).
\end{aligned}
$$
To summarize, we have the following one time step update of the fully discrete version of RMHD.
\begin{algorithm}[!h]
\caption{one step of fully discrete update for RMHD}
\SetAlgoLined
\KwIn{$\rho_{i}^s$, $v^s_{x,i}$, $v^s_{y,i}$, $v^s_{z,i}$, $T_{i}^s$, $B_{y,i}^s$, $B_{z,i}^s$, $I_{i}^s$,$J_{i}^s$ ($J_{i}^s=\average{I_{i}^s}$)}
\KwOut{$\rho_{i}^{s+1}$, $v^{s+1}_{x,i}$, $v^{s+1}_{y,i}$, $v^{s+1}_{z,i}$, $T_{i}^{s+1}$, $B_{y,i}^{s+1}$, $B_{z,i}^{s+1}$, $I_{i}^{s+1}$,$J_{i}^{s+1}$ ($J_{i}^{s+1}=\average{I_{i}^{s+1}}$)}

\BlankLine
     obtain $\rho_i^{s+1}$, $v^{s+1}_{y,i}$, $v^{s+1}_{z,i}$, $B_{y,i}^{s+1}$, $B_{z,i}^{s+1}$ from \eqref{fvm2}, \eqref{fvm4}, \eqref{fvm5}, \eqref{fvm7} and \eqref{fvm8}; 
    \\  obtain $R_{i}^{s}$, $Q_{ i}^{s}$ from  \eqref{eqn:rr} and  \eqref{eqn:011};
    \\  get $J_{i}^{s+1}$, $T_{i}^{s+1}$,$v_{x,i}^{s+1}$ and $R_{i}^{s+1}$ from the system coupled by \eqref{fvm3}, \eqref{fvm6}, \eqref{zerom} and \eqref{firm};
    \\    get $Q_{i}^{s+1}$  from  equation \eqref{fvm1} and $Q_{i}^{s}$ replaced by $Q_{i}^{s+1}$ ;
    \\  obtain $I_{i}^{s+1}$ from \eqref{eqn:011}.
\end{algorithm}

\section{Numerical Examples}
In this section, we conduct several numerical experiments to test the performance of our proposed method. The Dirichlet boundary conditions are adopted in all the following examples.
For velocity discretization, we choose the discrete-ordinate method \cite{SKTS}. In all numerical examples, $S_8$ discrete ordinate method has been used.

The first example is for the RMHD system where the radiation does not effect the fluid i.e. $\PP_0=0$, but the fluid provides a source for radiation. The second example is to illustrate the performance of the radiation solver when {\it both optical thick and thin regions coexist}, in which the fluid density $\rho$, velocity $\mathbf{v}$ and temperature $T$ are fixed. In example 3, a range of radiation-hydrodynamic shock problems are tested. In example 4, using different absorption and scattering coefficients, AP property of the proposed scheme for the full RMHD system is illustrated numerically. Finally, optical thick and thin regions {\it coexisting} case are tested in Example 5 and we can see that large time step is allowed.  
\subsection{Example 1}\label{section1}
\begin{figure}[htbp]
    \centering
    \subfigure[Density $\rho$]{
        \includegraphics[width=2.9in]{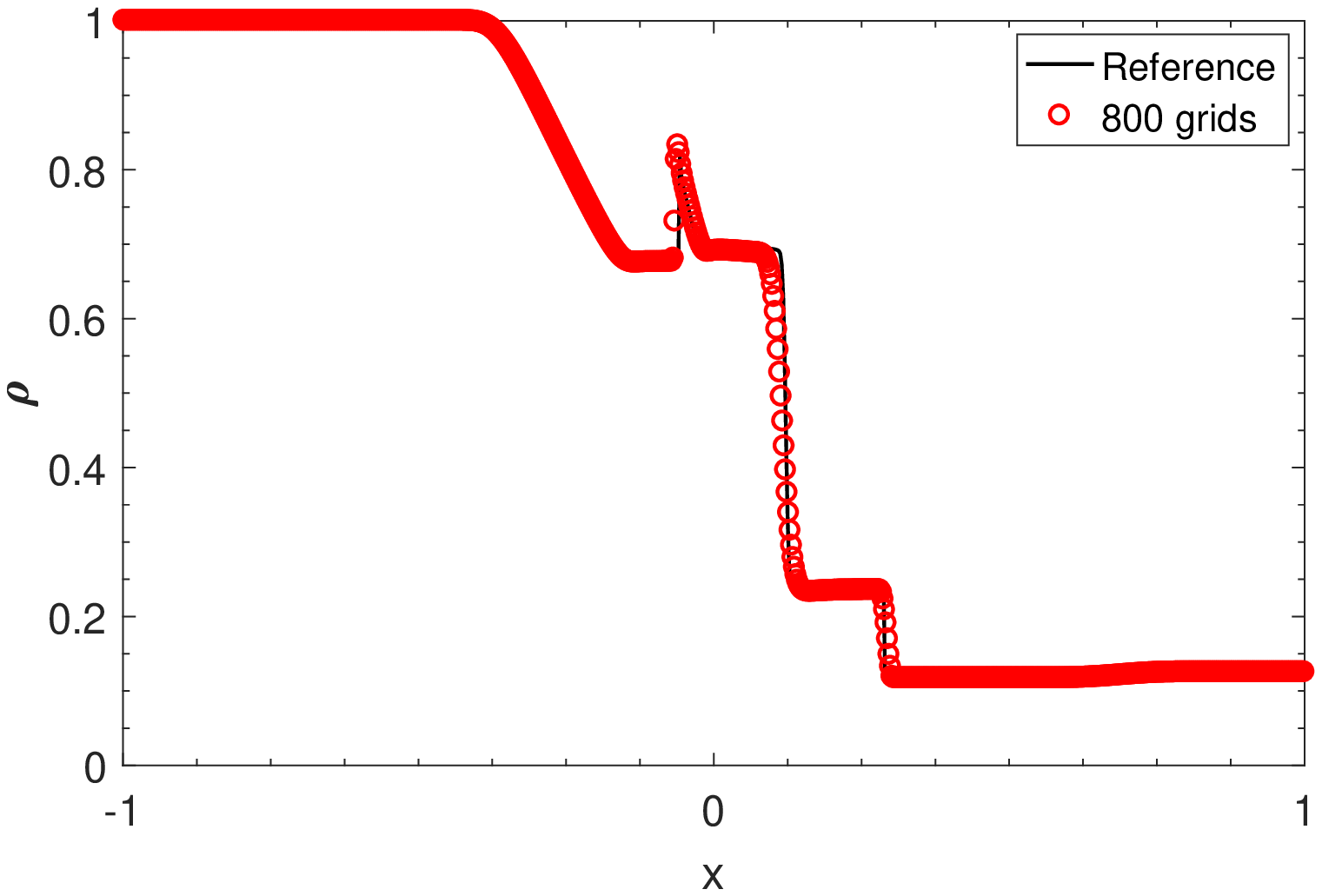}
    }
    \subfigure[Pressure p]{
	\includegraphics[width=2.9in]{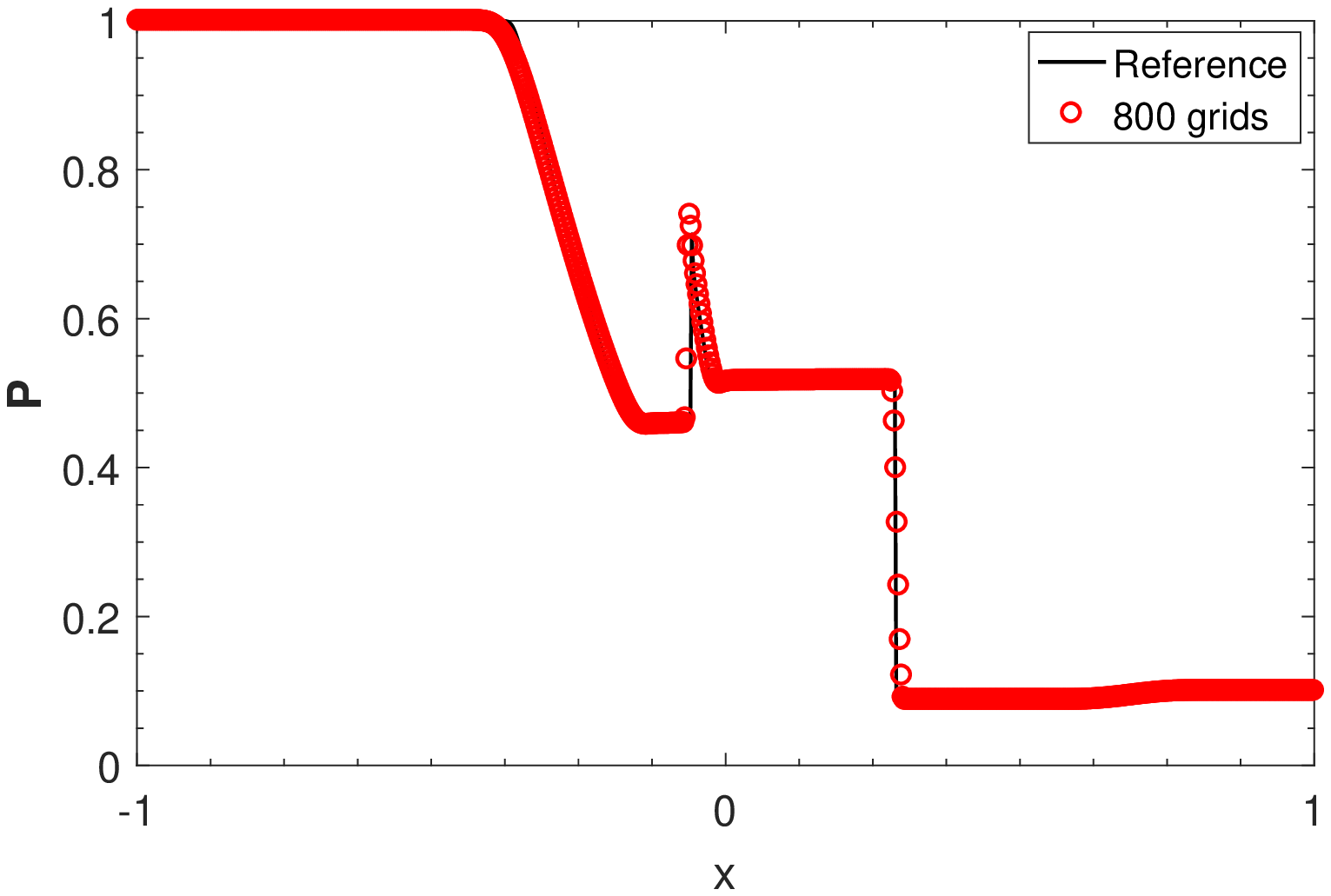}
    }
    \quad    
    \subfigure[Velocity in x-direction  $v_x$]{
    	\includegraphics[width=2.9in]{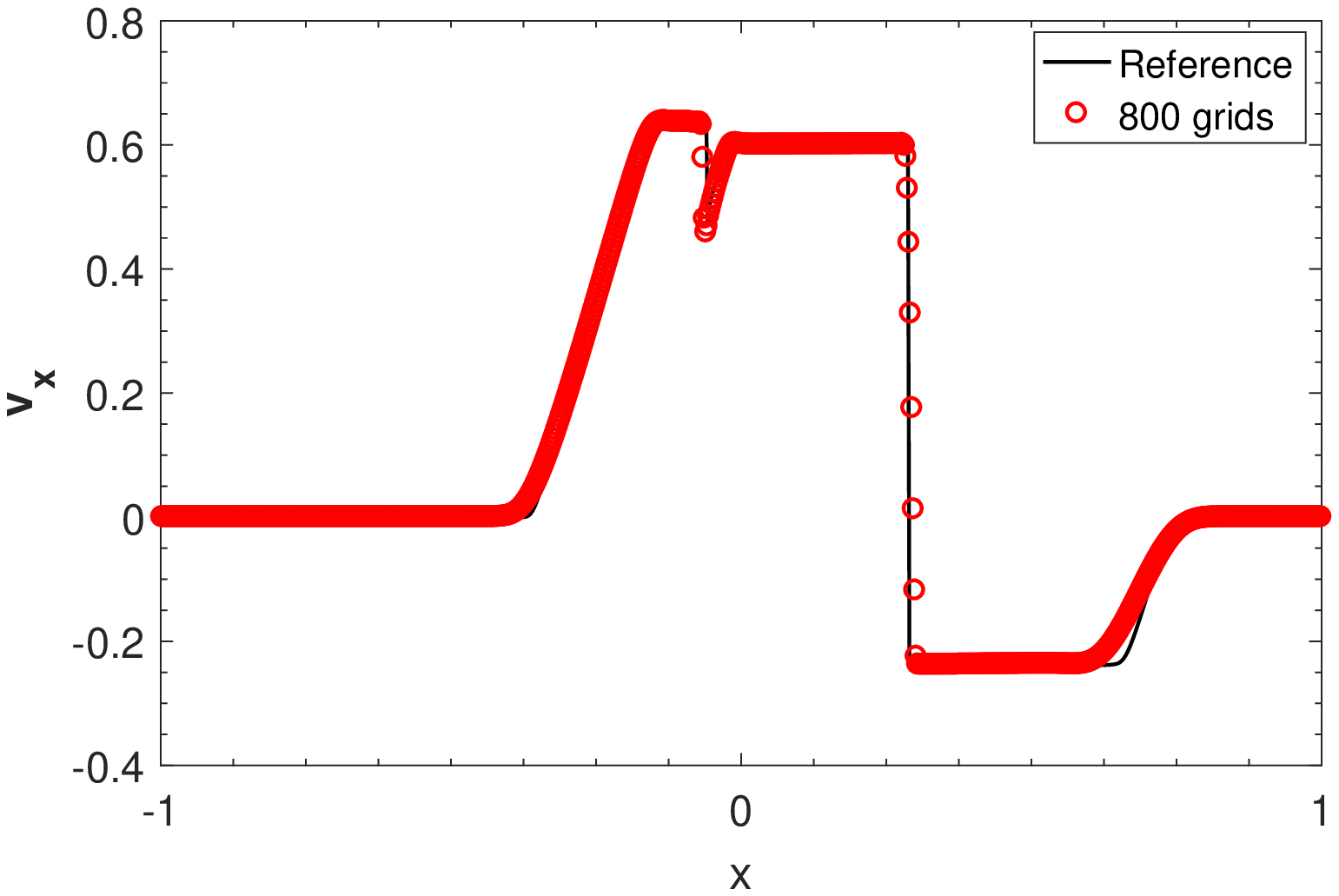}
    }
    \subfigure[Velocity in y-direction  $v_y$]{
	\includegraphics[width=2.9in]{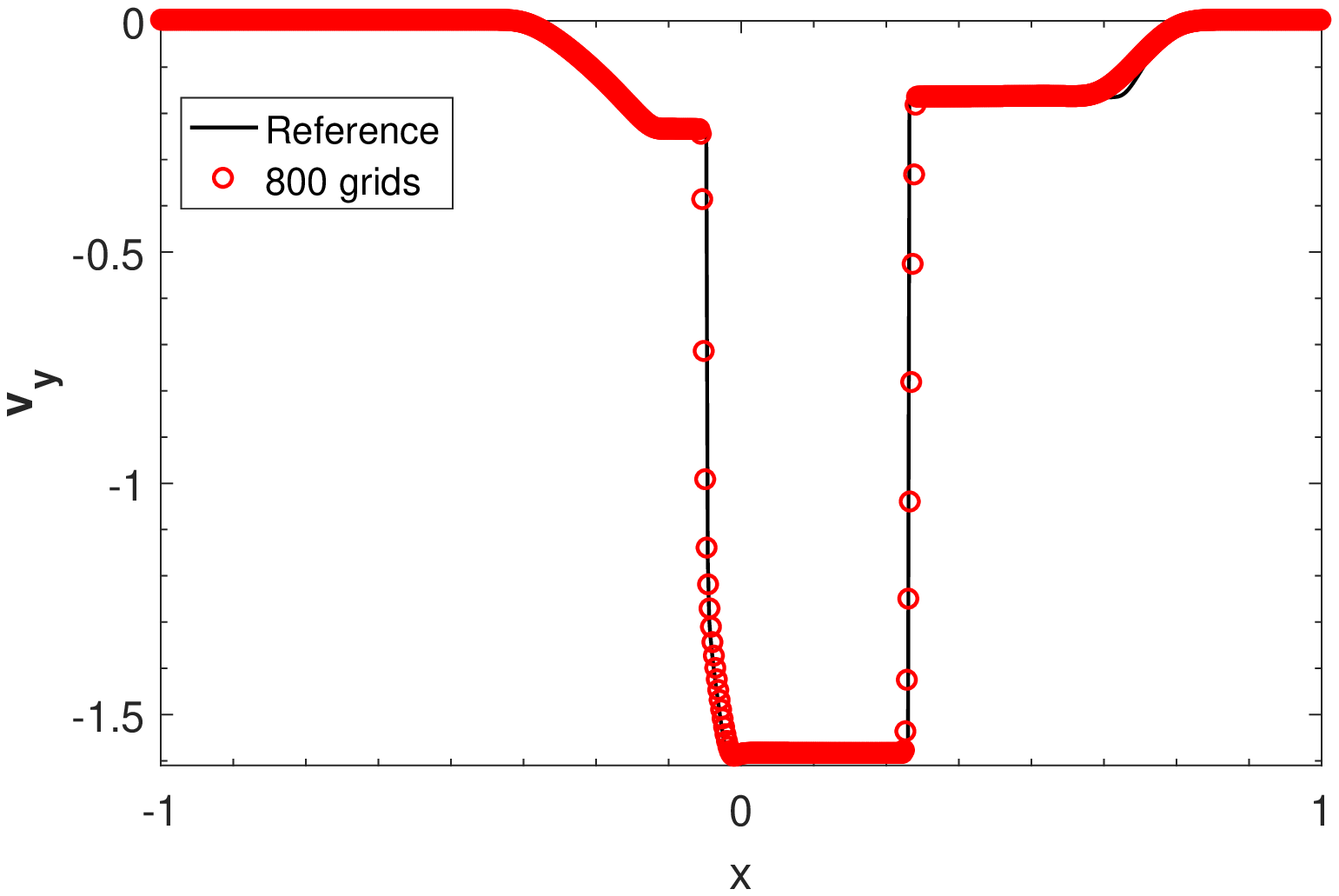}
    }
    \quad
    \subfigure[Magnetic field in y-direction  $B_y$]{
    	\includegraphics[width=2.9in]{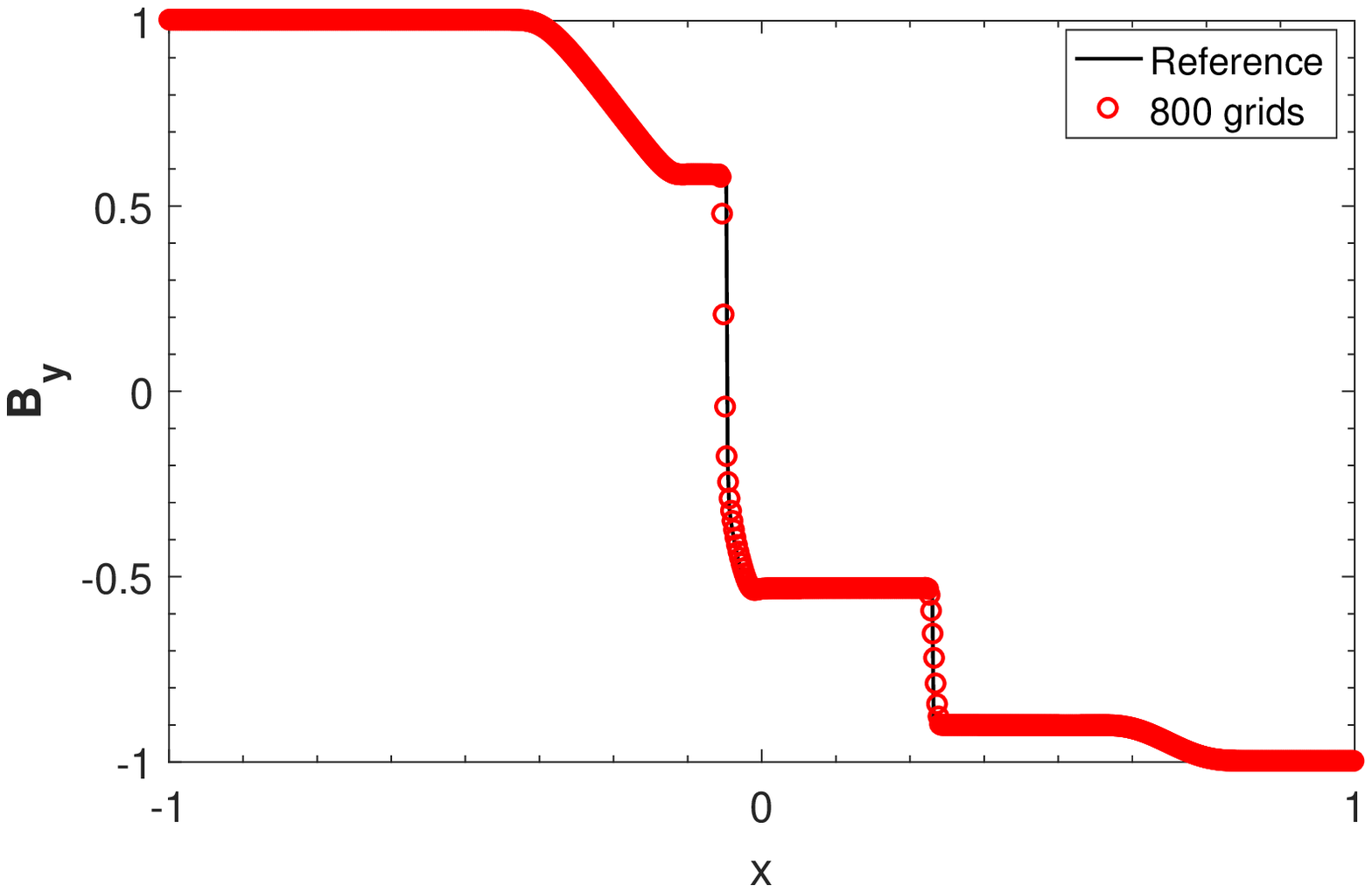}
    }
    \subfigure[The ratio of the pressure and the density]{
	\includegraphics[width=2.9in]{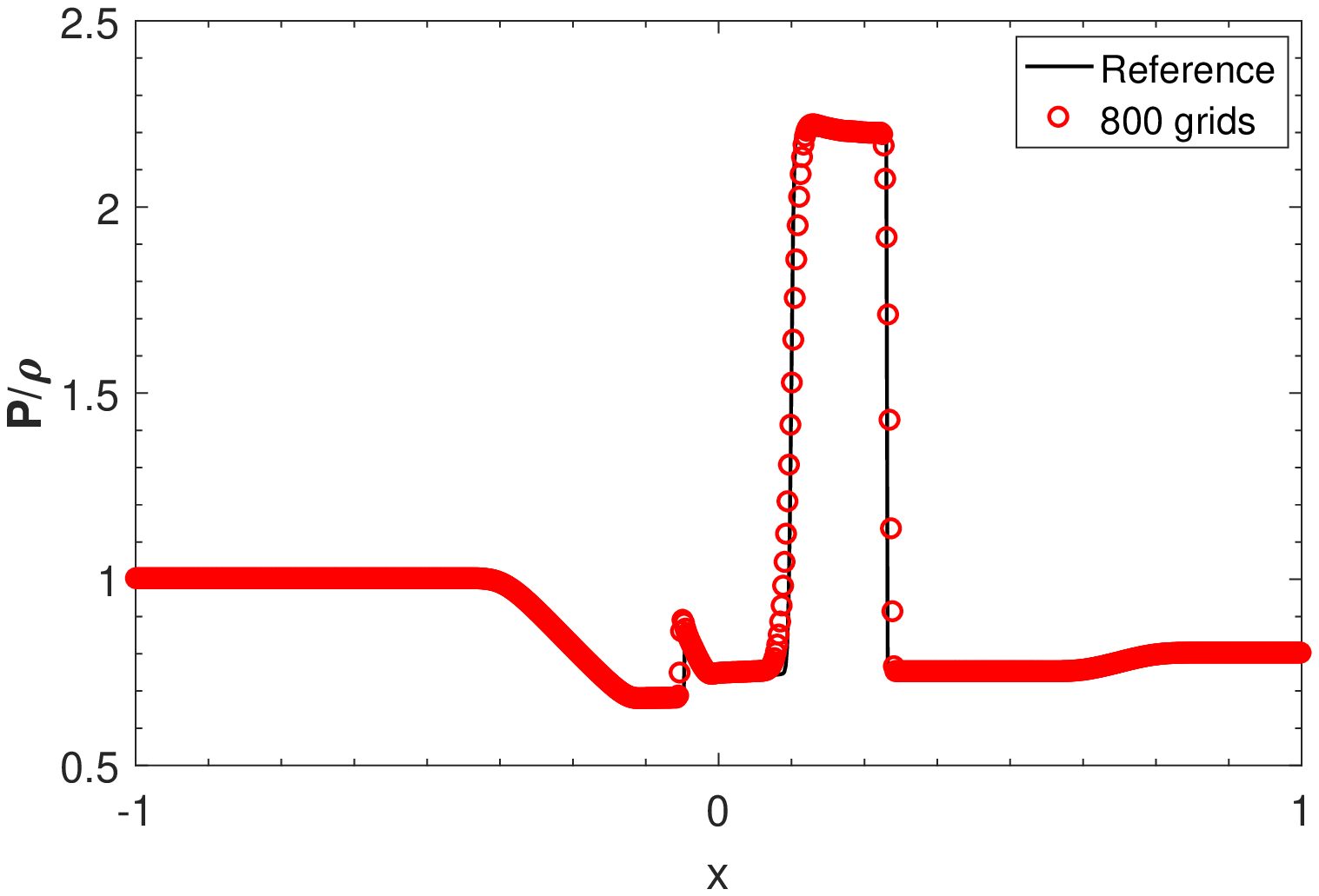}
    }
    \caption{Example 1. The numerical  results at time $t=0.2$  using our method  with $\Delta x=1/400$ (circles) and  the Roe solver with $\Delta x=1/2000$ (solid lines). For both schemes, $\Delta t = 0.2\Delta x$. }
    \label{fig.1}
\end{figure}

In this example, we will test an ideal MHD shock tube problem as in Section V in \cite{brio1988upwind}, where $\PP_0$=0 in \eqref{eqn:004}. The computational domain is $[-1,1]$ and the initial data  is given by
 \begin{equation}\label{equ:ini}
(\rho,v_x,v_y,v_z,B_y,B_z,p)(x,0)=\left\{
\begin{aligned}
 &(1,0,0,0,1,0,1),\qquad\qquad x<0,\\
&(0.125,0,0,0,-1,0,0.1),\quad x>0.
\end{aligned}
\right.
\end{equation}
The magnetic intensity in the $x$ direction  $B_x=0.75$, the ideal gas constant $ R_{\text{ideal}}=1$ and adiabatic index for an ideal gas $ \gamma=2$, the absorption collision cross-section $\sigma_a(x)=1/3$, the scatter collision cross-section  $\sigma_s(x)=1/3$, $c=0.1$ and $\eps=10^{-5}$ in the RMHD system \eqref{eqn:004}. The exact solutions at time $t>0$ involve two fast rarefaction waves, a slow compound wave, a contact discontinuity, and a slow shock.  Fig.~\ref{fig.1} displays the numerical solutions at $t=0.2$ obtained by our method   with $\Delta x=1/400$.  The reference solution is computed by  Roe's method shown in \cite{stone2008athena} with with $\Delta x=1/2000$. For all the  schemes, the time step is $\Delta t=0.2\Delta x.$  


\subsection{Example 2}
\begin{figure}[htbp]
	\centering
	\subfigure[absorption coefficient]{
		\includegraphics[width=2.9in]{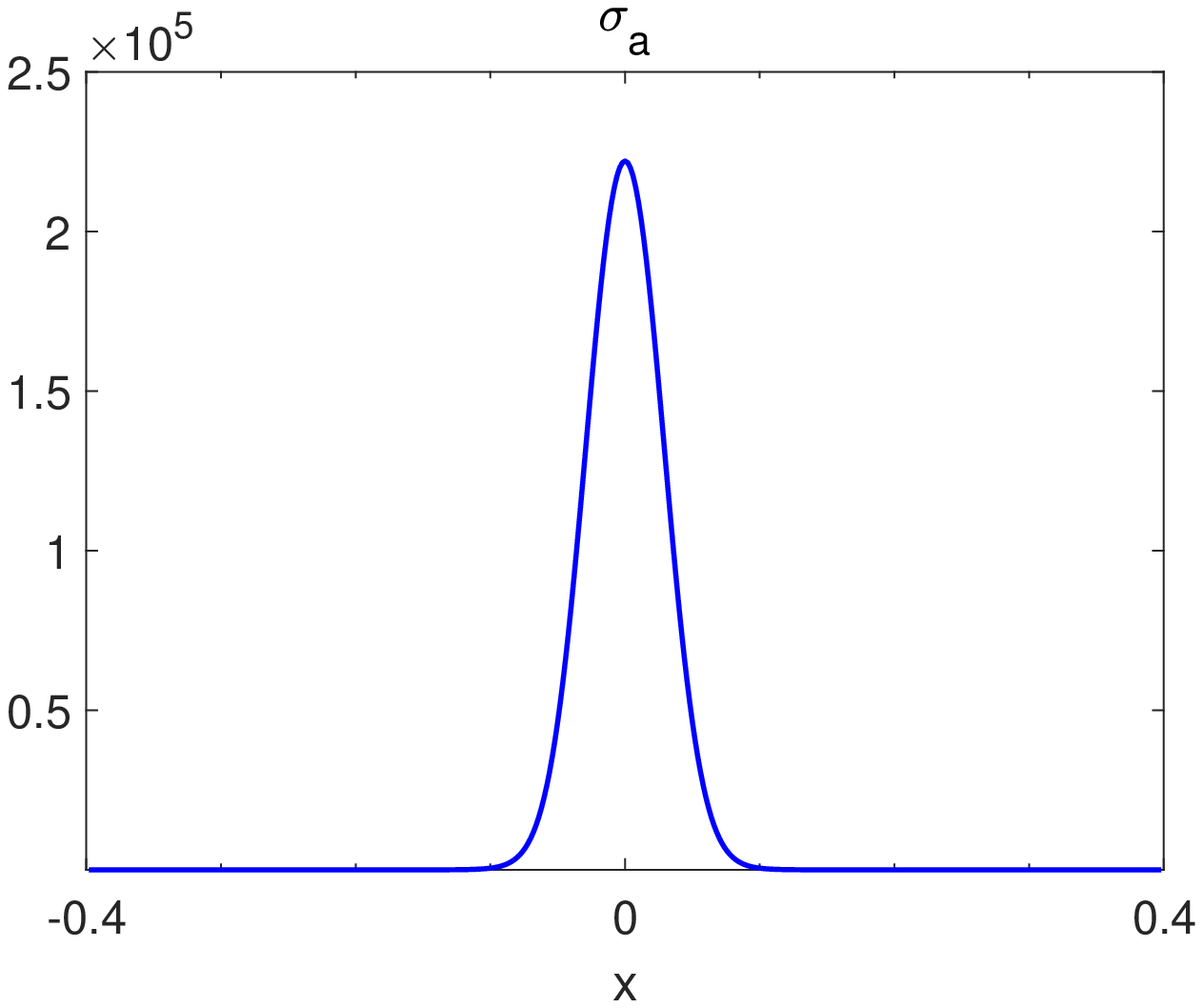}\label{fig.2_ab}
	}
	\subfigure[radiation temperature]{
		\includegraphics[width=2.9in]{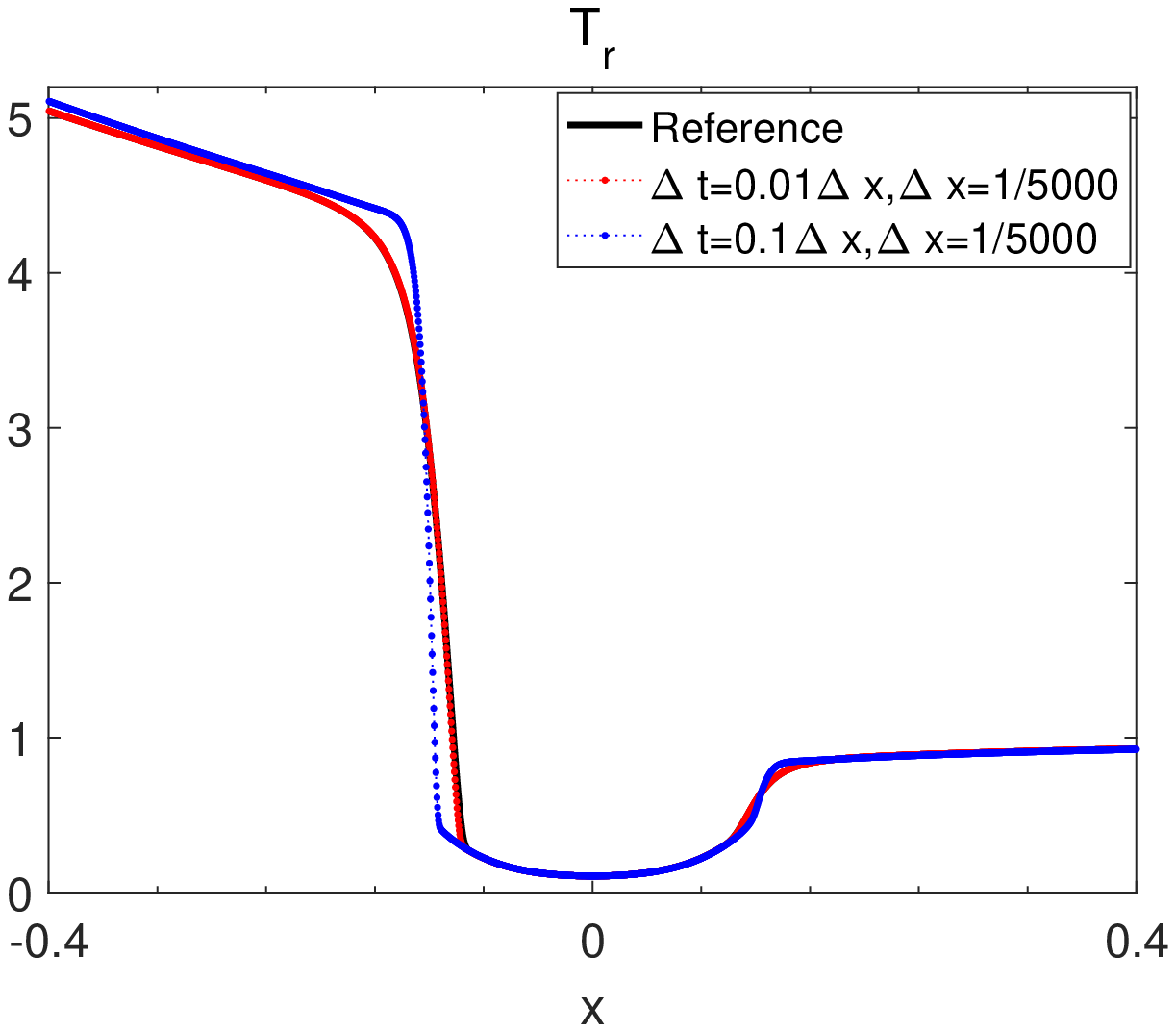}\label{fig.2_rt}
	}\\
	\subfigure[residual Q]{
		\includegraphics[width=2.9in]{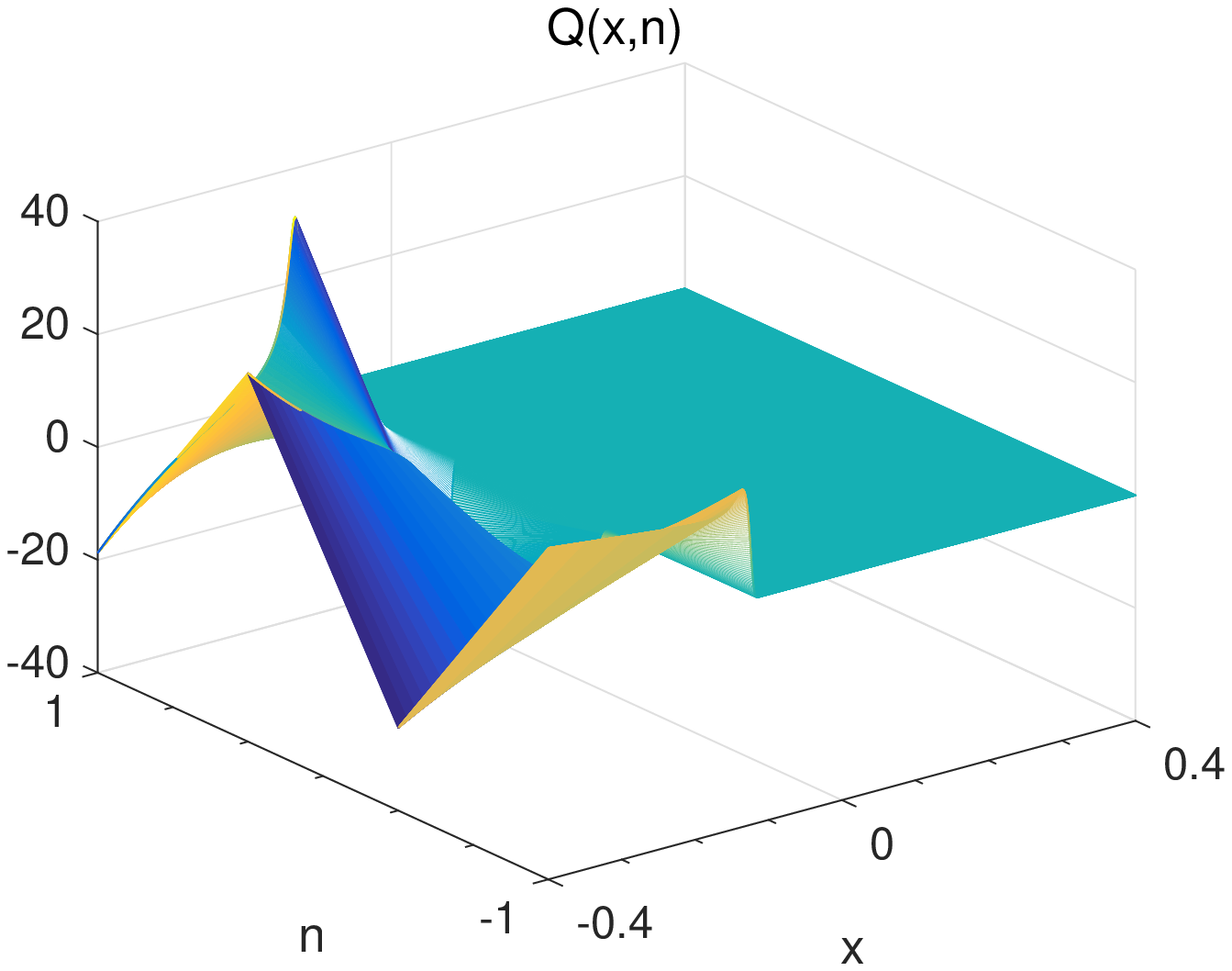}\label{fig.2_q1}
	}
	\subfigure[$K_Q$]{
		\includegraphics[width=2.9in]{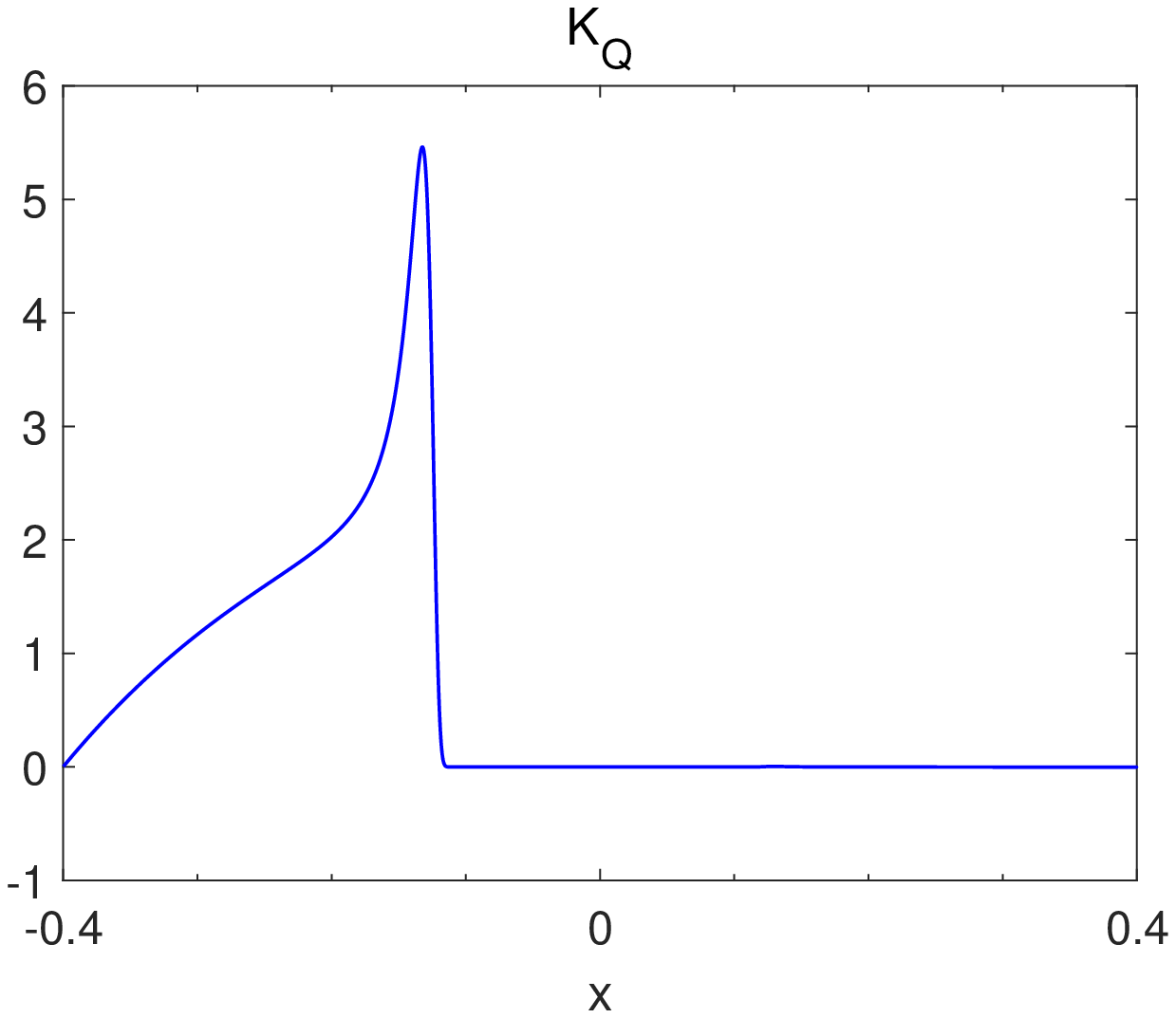}\label{fig.2_kq1}
	}
	\caption{Example 2. a) absorption coefficient $\sigma_a$; b) comparison of the numerical results at time $t =10$, obtained by our AP scheme and the explicit solver; \textcolor{black}{(The numerical result obtained by the AP solver using $\Delta t=0.01\Delta x$ and $\Delta x=1/5000$ overlaps with the reference solution.)} c) residual term $Q(x,n)$;\textcolor{black} {d) $K_Q$.} Here for AP scheme, we use $\Delta x=1/5000$, the time step $\Delta t =0.1 \Delta x$ \textcolor{black}{and $\Delta t =0.01 \Delta x$} For the explicit solver, we use  $\Delta x=1/5000$ and the time step $\Delta t =5*10^{-7}\Delta x$. }
\end{figure}

A problem with both optically thick and thin regions is tested for the radiation subsystem. This kind of problems has been studied in \cite{pub.1041264084,pub.1005319601,filbet2010class,hayes}, where the gas is held inactive, and only the radiation equation in system \eqref{eqn:non} is evolved, i.e., the evolution of the radiation does not change the gas density and momentum. In this test,  a blob of optically thick gas is located in the middle of the computational domain $[-0.4, 0.4]$, surrounded by an optically thin gas. The gas density is 
$$
\rho(x)=1+\f{11}{2}\left[{\rm tanh}(1-17x)+{\rm tanh}(1+17x)\right],
$$ 
and the gas temperature is adjusted to provide constant pressure in the domain, which is set to be $1/\rho$.  The initial radiative energy density $E_r$ is at local thermal equilibrium. The boundary condition is set to be $T(0,t)=6$ and $\mathcal{C}=5000$, $\mathscr{L}_a=1$, $\sigma_s=0$ and $\sigma_a=\rho^2T^{-3.5}$. See Fig. \ref{fig.2_ab} for the absorption coefficient, which varies from 1 to $10^5$. For this example, we compare the radiation temperature obtained by our AP scheme and an explicit solver on a finer mesh in Fig. \ref{fig.2_rt}. We use the space step $\Delta x=1/5000$ for both AP and the explicit solver, but the time step sizes are different such that \textcolor{black}{$\Delta t =0.1 \Delta x$ and $\Delta t =0.01 \Delta x$ for the AP solver and $\Delta t =5*10^{-7}\Delta x$ for the explicit solver. The good agreement between our and the reference solution indicates that the AP method works well in the coexistence case. Moreover, we plot  $K_Q$ in Fig. \ref{fig.2_kq1}, the value of $K_Q$ can be as large as $J$ in this example.} From Fig. \ref{fig.2_rt}, when at $t=10$, the radiation was stopped by the opaque material.
\subsection{Example 3}
In the simulation of RMHD, it is a challenging task to produce accurate radiative shocks, especially in the optically thick regime. We test a range of the radiation-hydrodynamic shock problems presented in \cite{lowrie2008radiative} in which several shocks are present in the solution. For each of these shocks, we set $\PP_0=10^{-4}$, adiabatic index for an ideal gas $\gamma=5/3$, absorption collision cross-section $\sigma_a=1/3$, scatter collision cross-section $\sigma_s=10^6$, $c=3\times 10^6$, $\eps=1$ in RMHD system \eqref{eqn:004}, and  the radiation temperature $T_r\equiv(4\pi J)^{0.25}$.  Moreover, the computational domain is $[-0.02,0.02]$ and the results at $t=0.04$ are displayed. In this example, $\Delta x=1/800,\ \Delta t=0.2\Delta x$. 
In \cite{lowrie2008radiative}, the authors have obtained some semi-analytic solutions for the non-equilibrium diffusion limit system \eqref{eqnlimit1} at different Mach numbers. We compare the numerical results with the semi-analytic solutions obtained \cite{lowrie2008radiative}. 
\begin{equation} \label{mach1.2}
(\rho,v_x,v_y,v_z,B_y,B_z,T,T_r)(x,0)=\left\{
\begin{aligned}
 &(1,1.2,0,0,0,0,1,1),\qquad\qquad x<0,\\
&(1.298088,0.9244363,0,0,0,0,1.194888,1.194888),\quad x>0,
\end{aligned}
\right.
\end{equation}
\begin{equation}  \label{mach2}
(\rho,v_x,v_y,v_z,B_y,B_z,T,T_r)(x,0)=\left\{
\begin{aligned}
 &(1,2,0,0,0,0,1,1),\qquad\qquad x<0,\\
&(2.287066,0.874482876,0,0,0,0,2.077223,2.077223),\quad x>0.
\end{aligned}
\right.
\end{equation}

\begin{figure}[htbp]
	\centering
	\subfigure[Density $\rho$]{
		\includegraphics[width=2.9in]{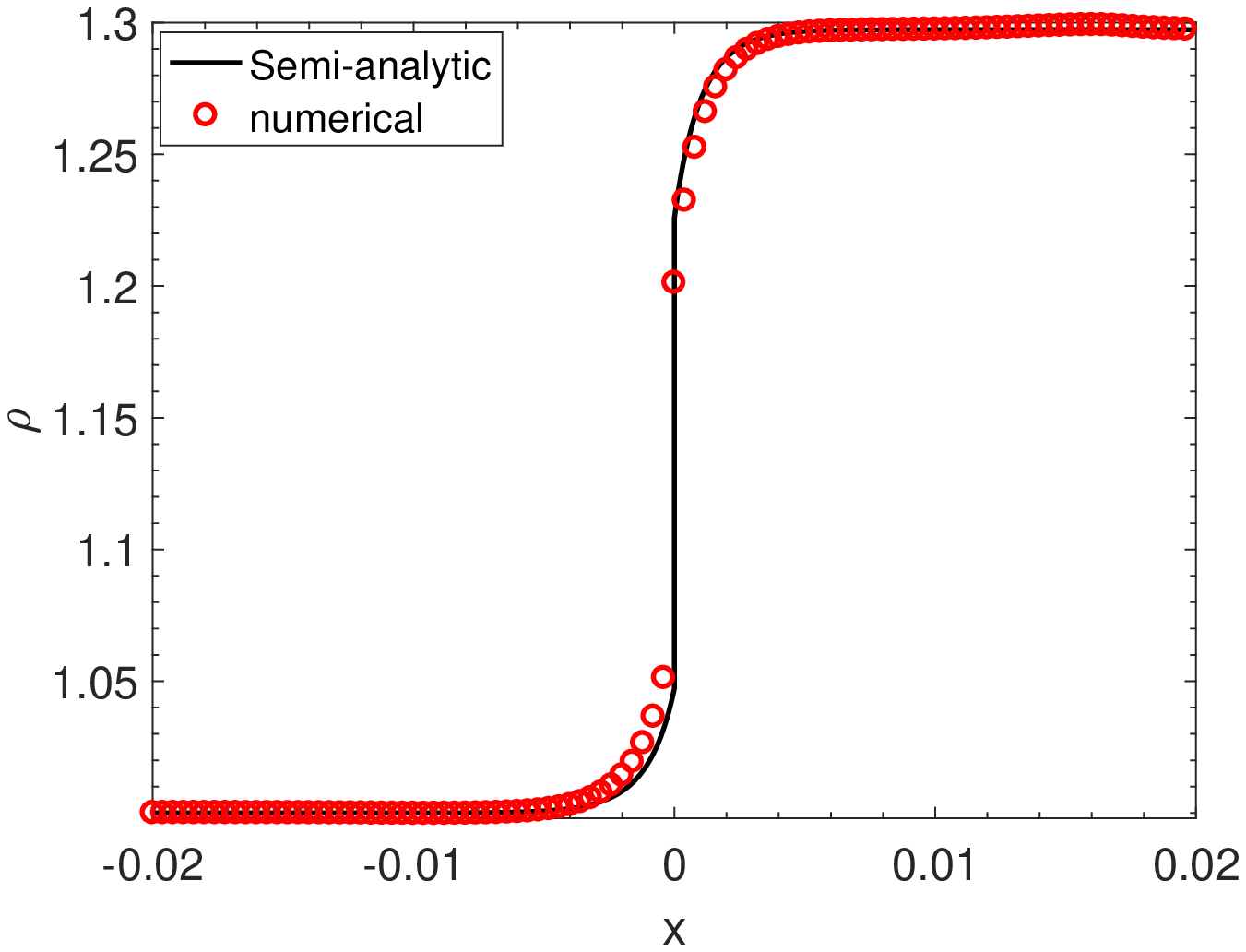}
	}
	\subfigure[Velocity $v$]{
		\includegraphics[width=2.9in]{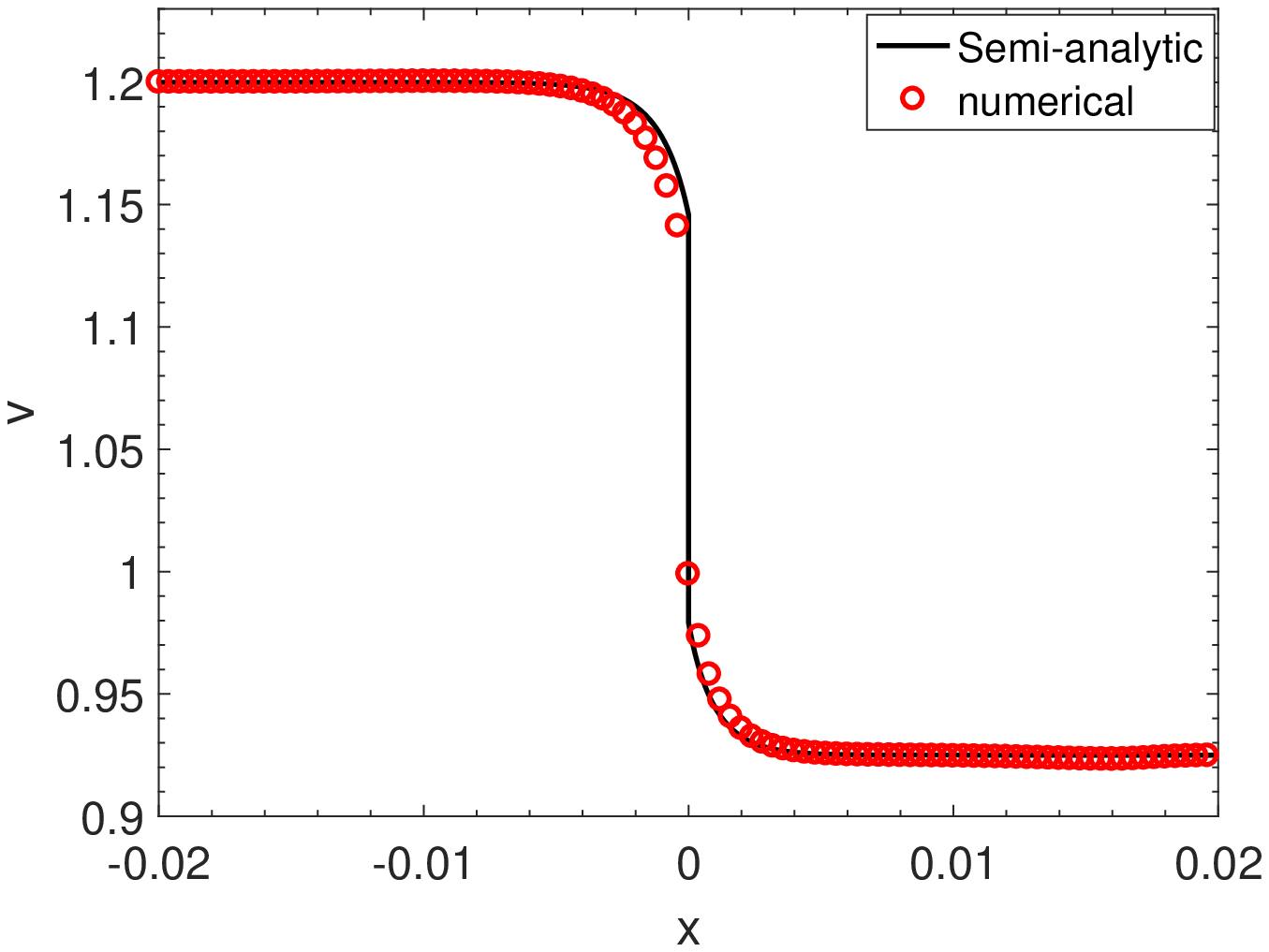}
	}
	\quad
	\subfigure[Temperature $T$]{
		\includegraphics[width=2.9in]{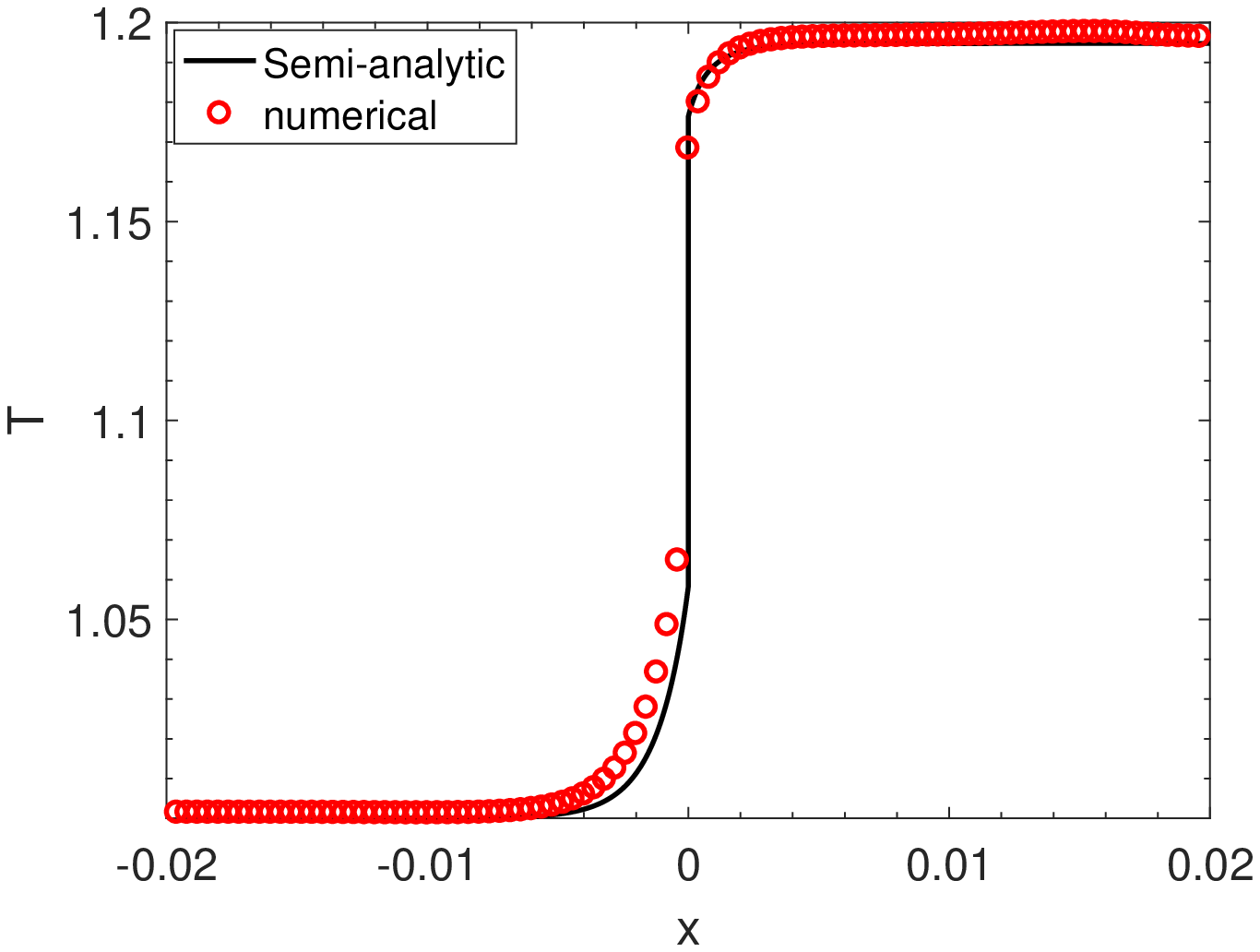}
	}
	\subfigure[Radiation temperature $T_r$]{
		\includegraphics[width=2.9in]{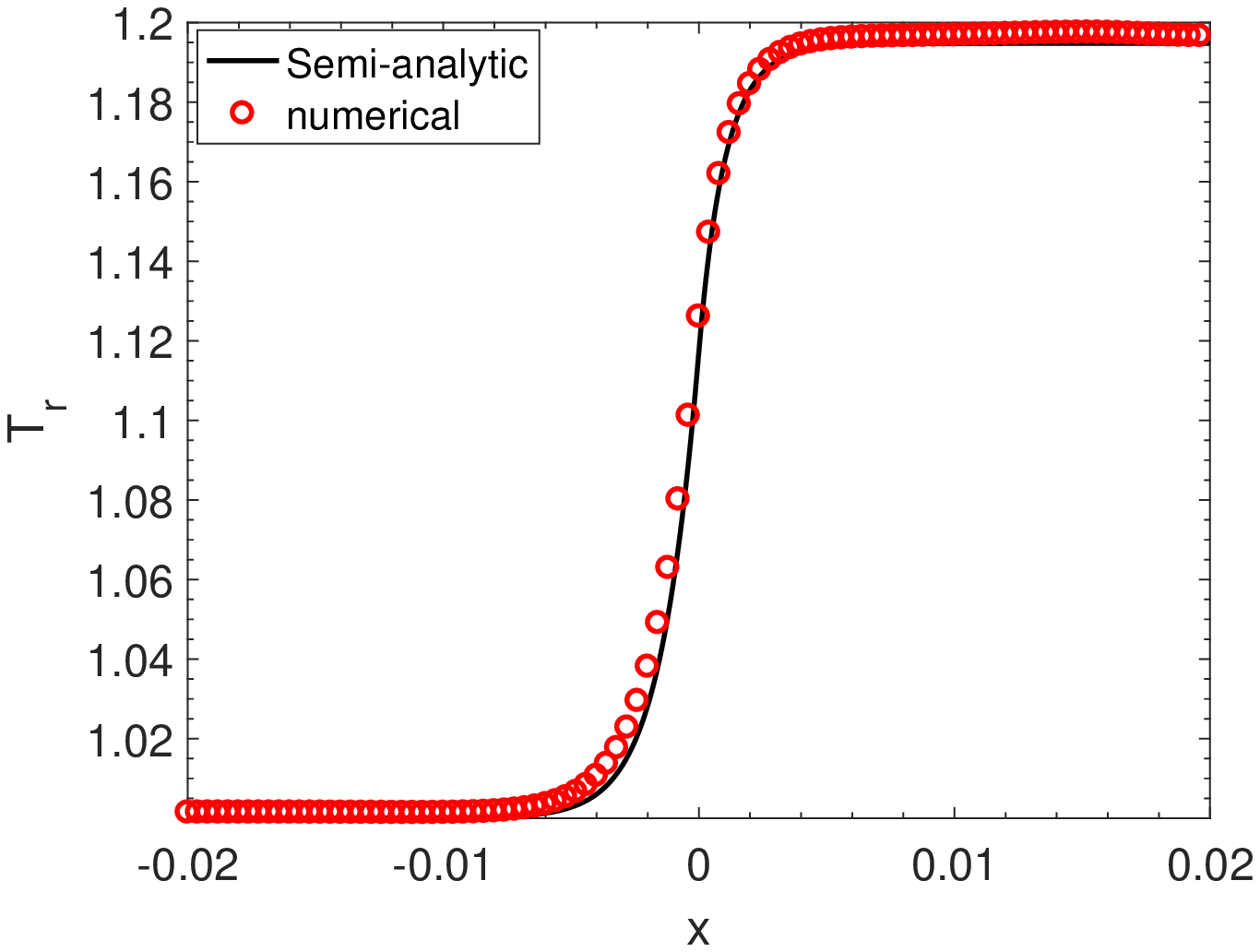}
	}
 \caption{Comparison of the results using our AP scheme at time $t = 0.04$ and the semi-analytic solution for Mach number $\mathcal{M}$ = 1.2. For AP scheme, we use  $\Delta x = 1/800$ and $\Delta t = 0.2\Delta x$.}
	\label{fig.1.2}
\end{figure}

First of all, a Mach 1.2 shock is considered, which has no isothermal sonic point (ISP) but a hydrodynamic shock.  The initial conditions are shown in \eqref{mach1.2}. Fig.~\ref{fig.1.2} compares our numerical results with the semi-analytic solutions, in which we can see good agreement for all quantities, including   density, velocity, material and radiation temperature. Due to the hydrodynamic shock, there are discontinuities in the solution profiles of density, velocity and material temperature, and the maximum material temperature is bounded, since there is no ISP to drive it further.

\begin{figure}[htbp]
	\centering
	\subfigure[Density $\rho$]{
		\includegraphics[width=2.9in]{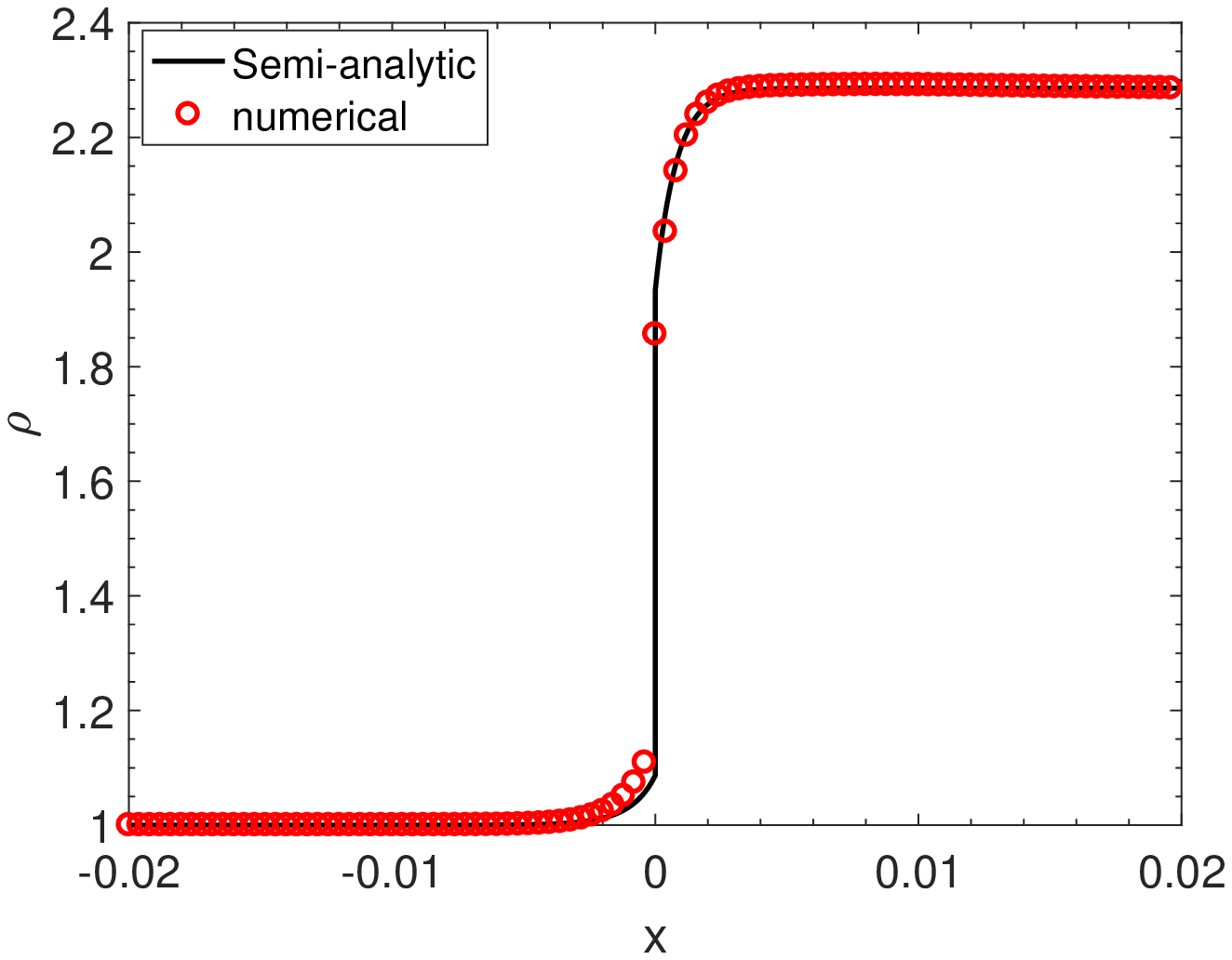}
	}
	\subfigure[Velocity $v$]{
		\includegraphics[width=2.9in]{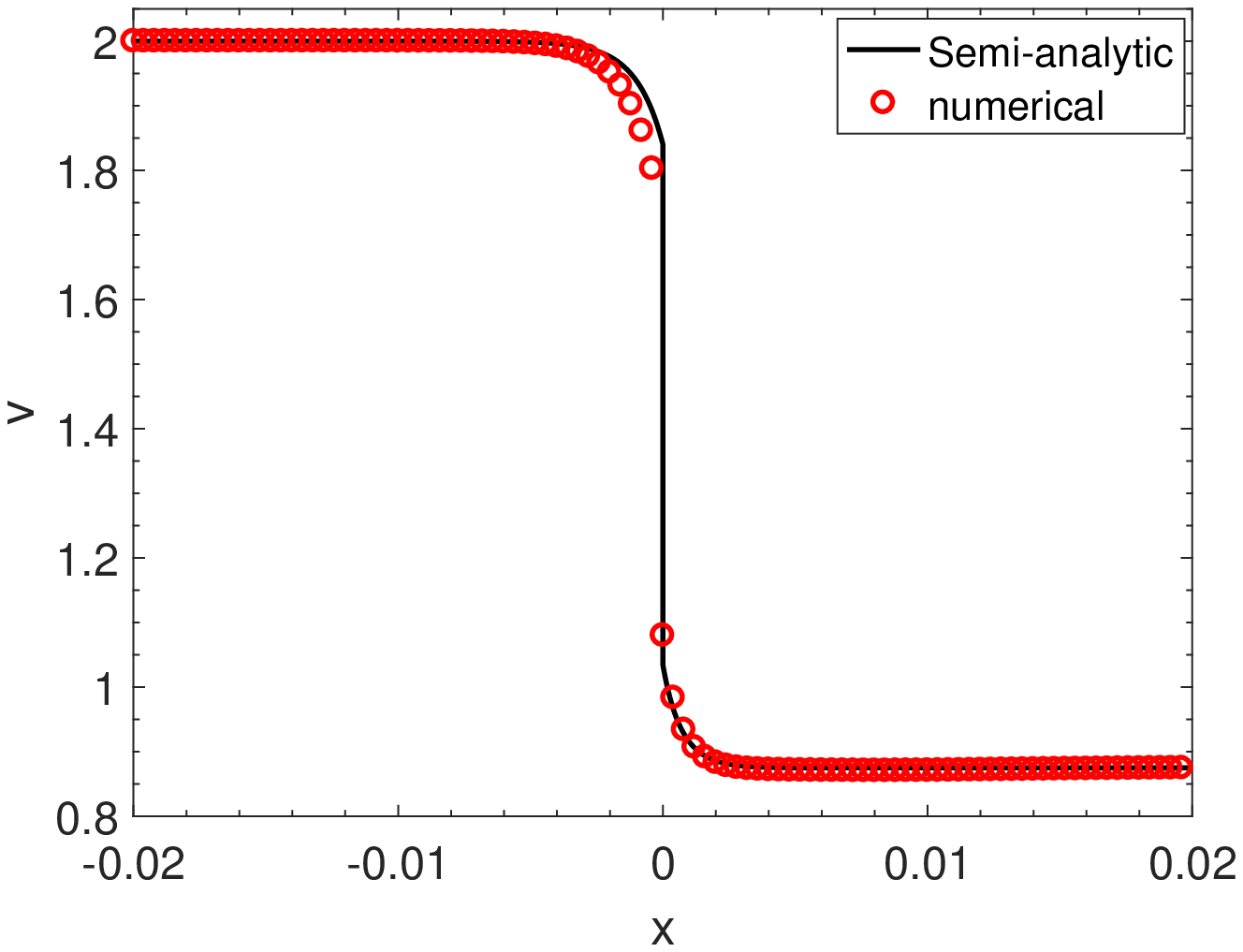}
	}
	\quad
	\subfigure[Temperature $T$]{
		\includegraphics[width=2.9in]{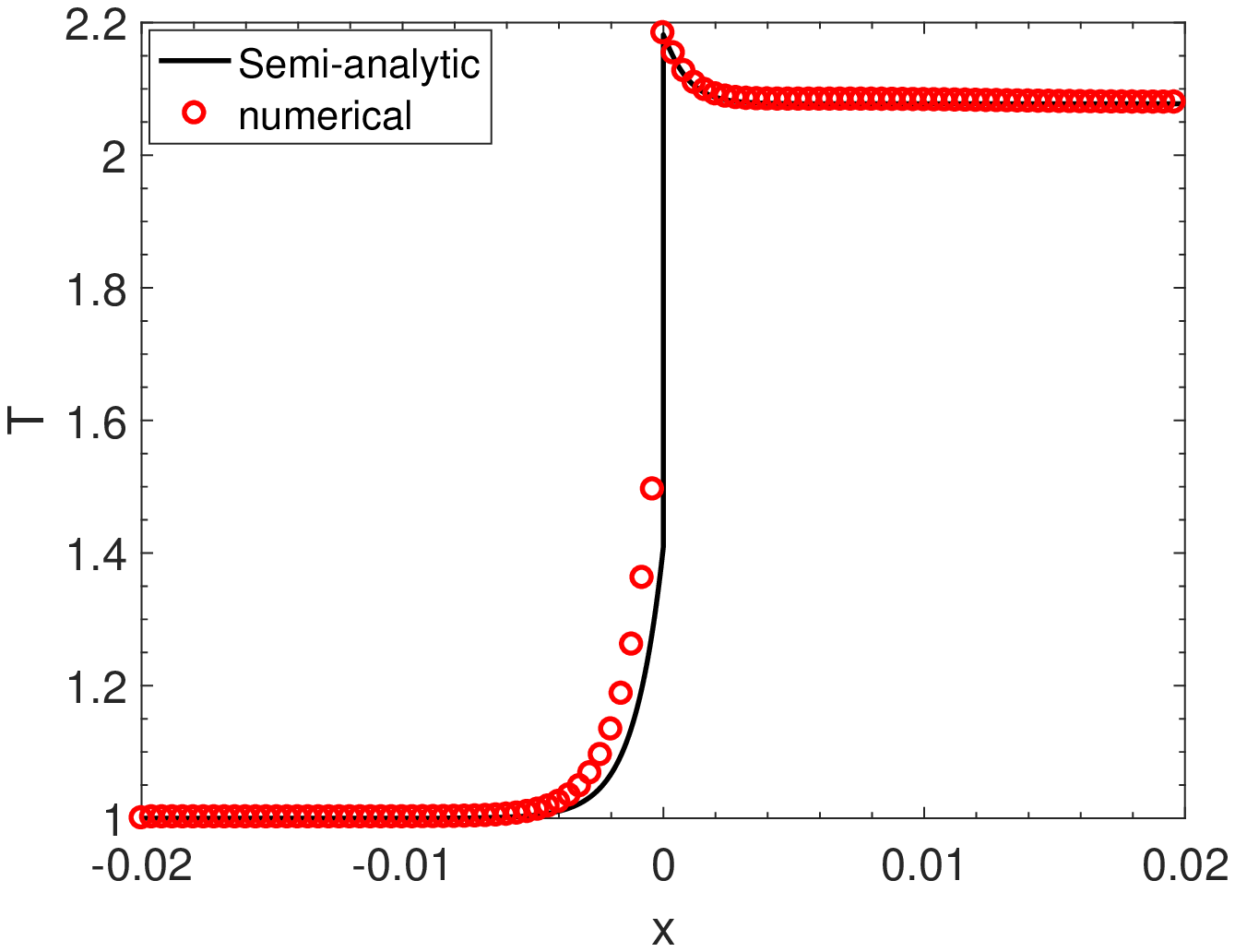}\label{ze}
	}
	\subfigure[Radiation temperature $T_r$]{
		\includegraphics[width=2.9in]{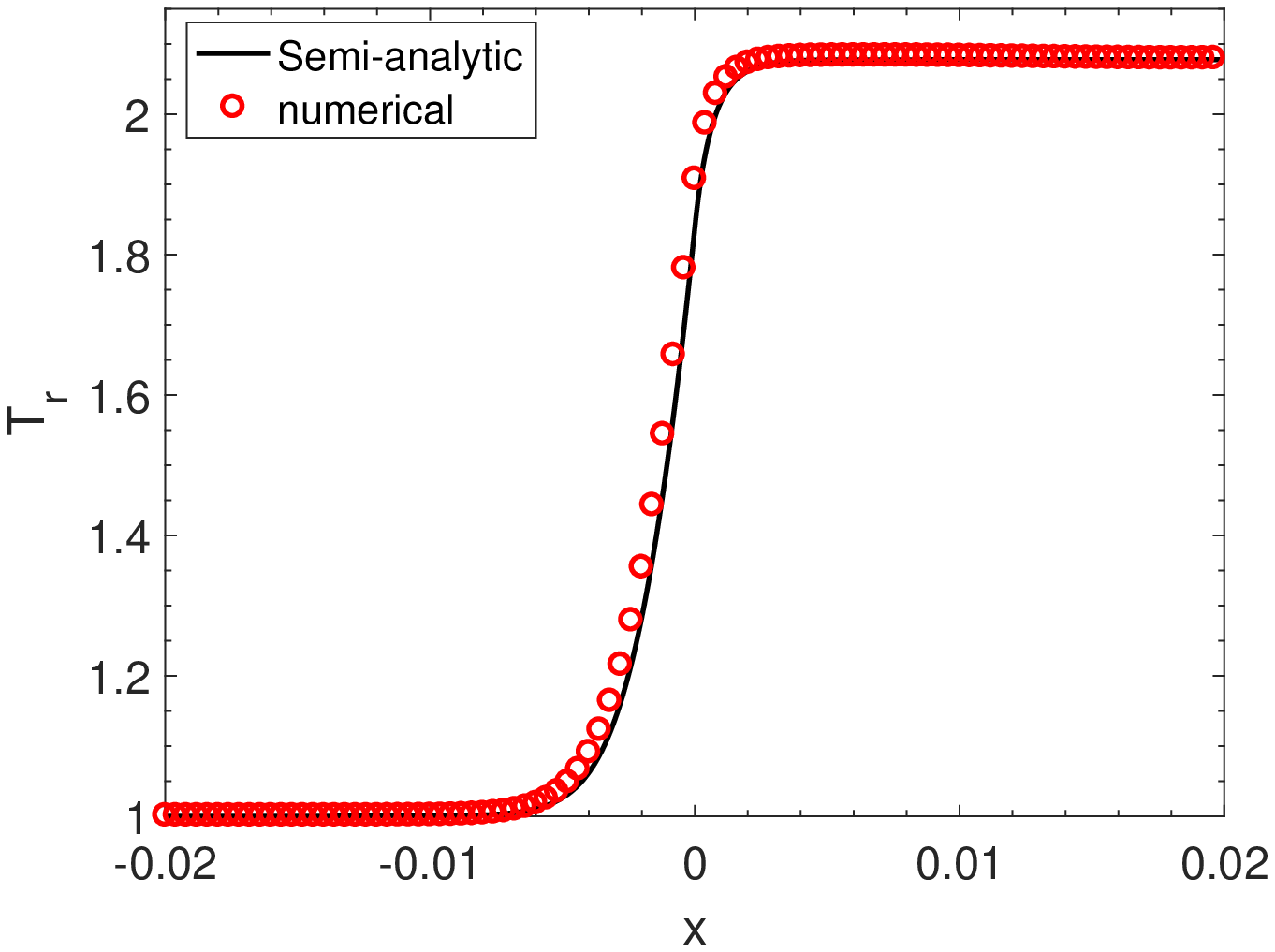}
	}
 \caption{Example 3. Comparison of the results using our AP scheme at time $t = 0.04$ and the semi-analytic solution for Mach number $\mathcal{M}$ = 2. For AP scheme, we use  $\Delta x = 1/800$ and $\Delta t = 0.2\Delta x$.}
	\label{fig.2.0}
\end{figure}

At Mach 2, there are both a hydrodynamic shock and an ISP. The initial conditions are given in \eqref{mach2}. Fig.~\ref{fig.2.0} compares our numerical results  with the semi-analytic solutions. One can see that our numerical results are  in good agreement with the semi-analytic solutions. In Fig.~\ref{fig.2.0}, discontinuities can   be seen in  material density, velocity and material temperature due to the hydrodynamic shock. Moreover, the Zel’dovich spike can be observed in material temperature as in Fig.~\ref{ze}. The Zel’dovich spike is caused by the ISP embedded within the hydrodynamic shock, which drives up the material temperature at the shock front.

\subsection{Example 4}\label{ex3}
\begin{figure}[ht]
	\centering
	\includegraphics[width=1\textwidth]{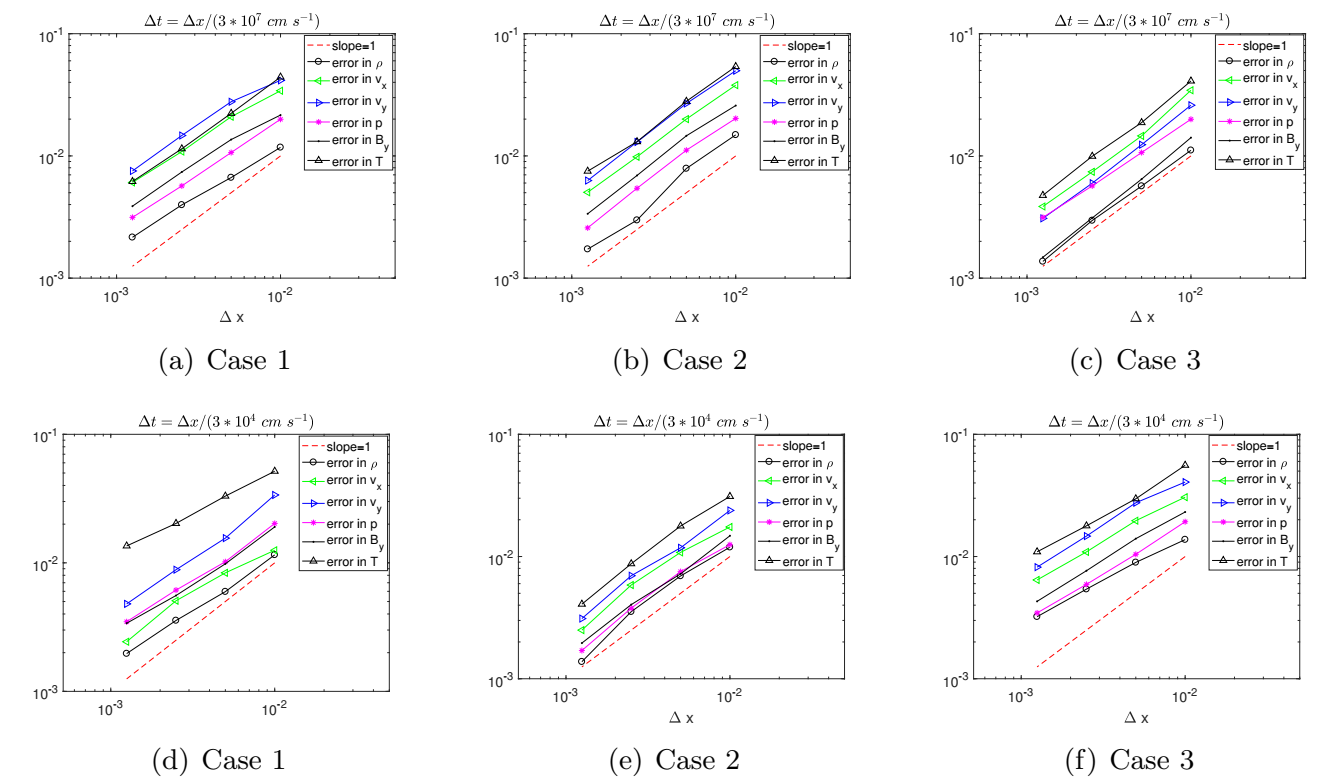}
 \caption{Example 4. Errors of different $\Delta x$, $\Delta t$ for the three cases. Here different $\Delta x=10^{-2}\ cm$, $5*10^{-3}\ cm$, $2.5*10^{-3}\ cm$, $1.25*10^{-3}\ cm$, $6.25*10^{-4}\ cm$ are tested and in (a), (b) and (c), $\Delta t$ are chosen to be $\Delta x/(3*10^7\ cm\ s^{-1})$. In (d), (e) and (f), $\Delta t$ are chosen to be $\Delta x/(3*10^4\ cm\ s^{-1})$. }
    \label{fig:er1-case}
\end{figure}

\begin{figure}[htbp]
	\centering
	\subfigure[$L^2$ norm of $\average{Q}$]{
		\includegraphics[width=2.9in]{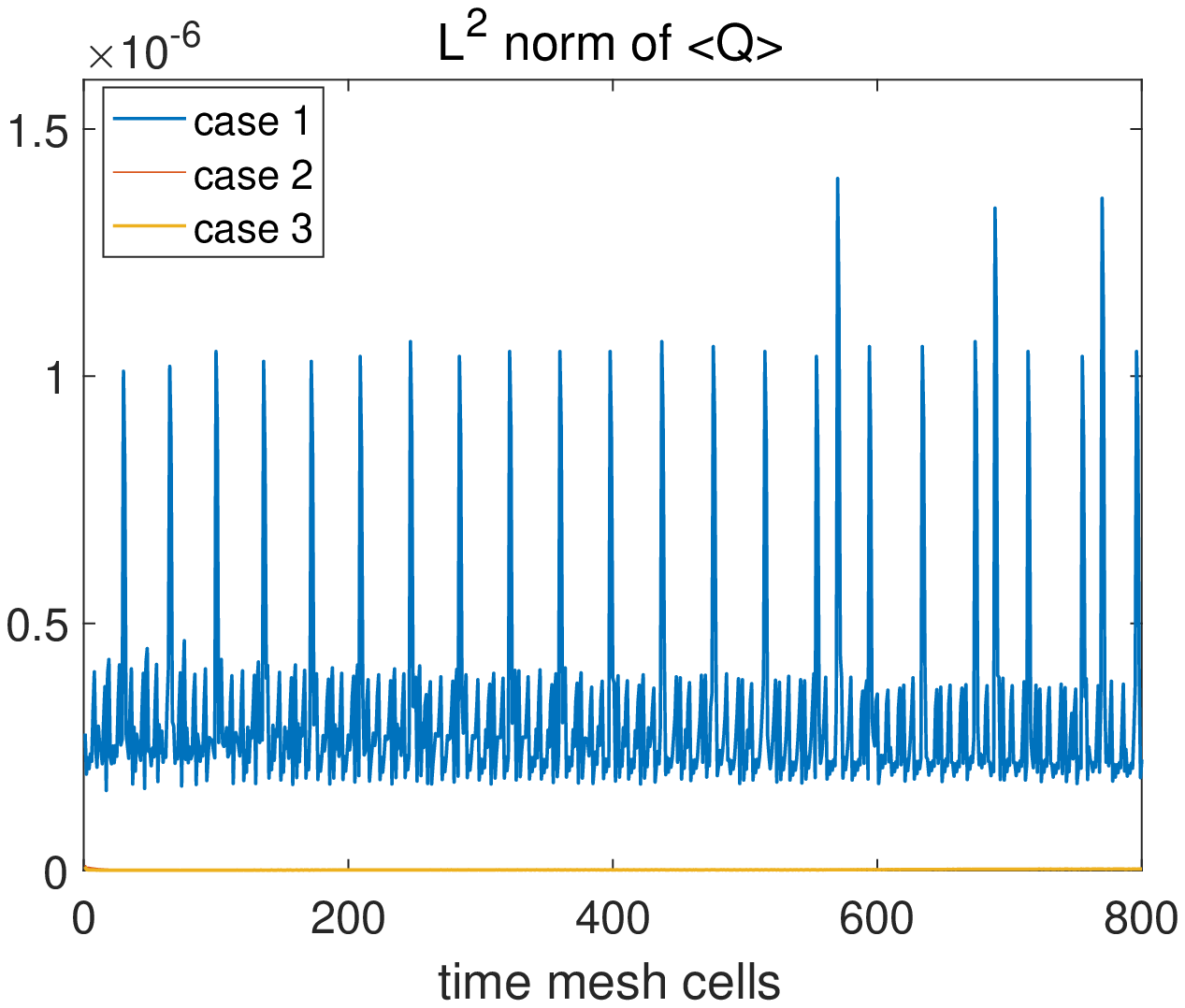}
	}
	\subfigure[$L^2$ norm of $\average{nQ}$]{
		\includegraphics[width=2.9in]{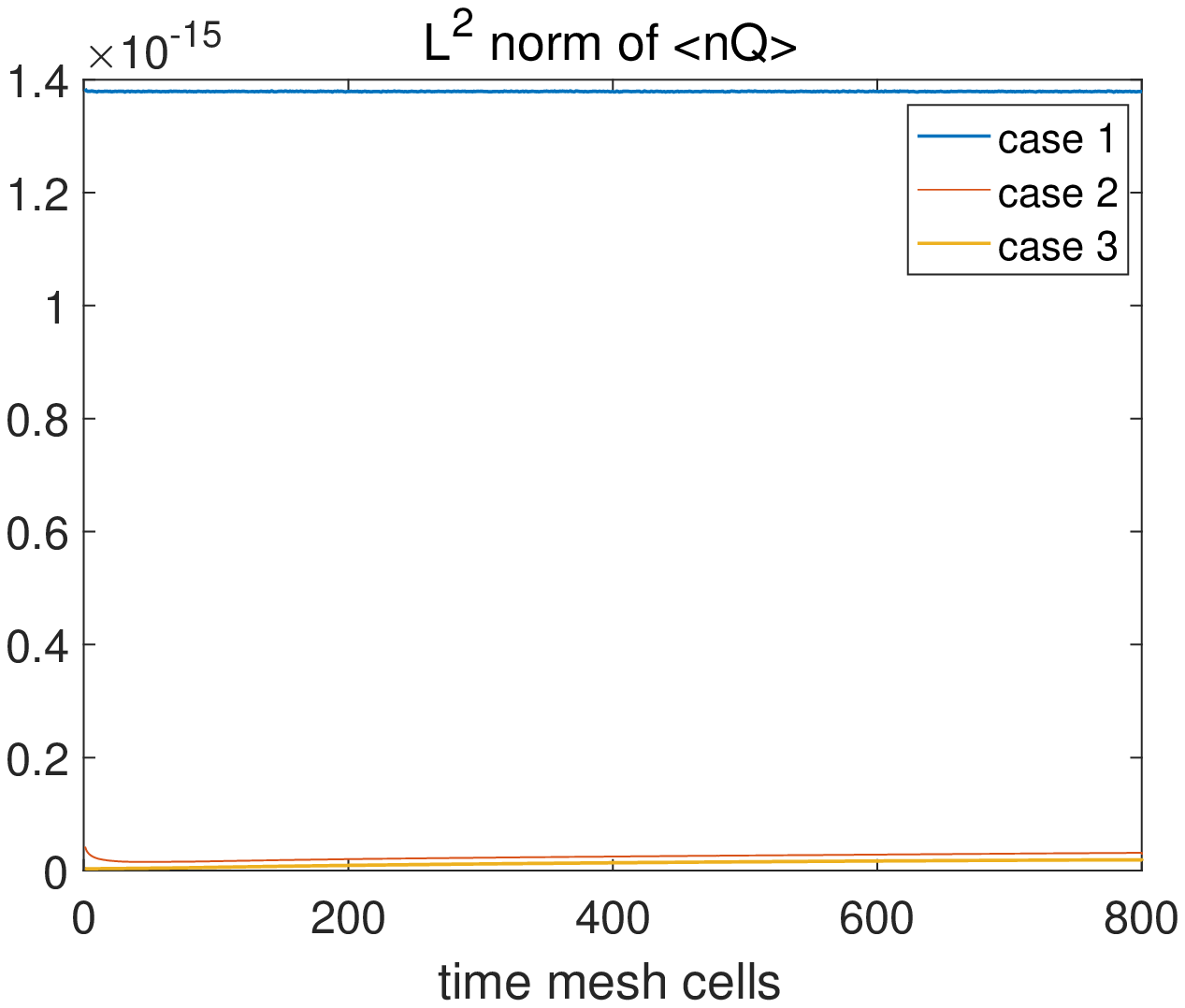}
	}
 \caption{\textcolor{black}{Example 4. the $L^2$ norm of $\average{Q}$ and $\average{nQ}$ in each time step with $\Delta x=2.5*10^{-3}\ cm$ and $\Delta t=\Delta x/(3*10^4\ cm\ s^{-1})$.}}
	\label{fig.4.2.0}
\end{figure}

\begin{figure}[htbp]
	\centering
	\includegraphics[width=1\textwidth]{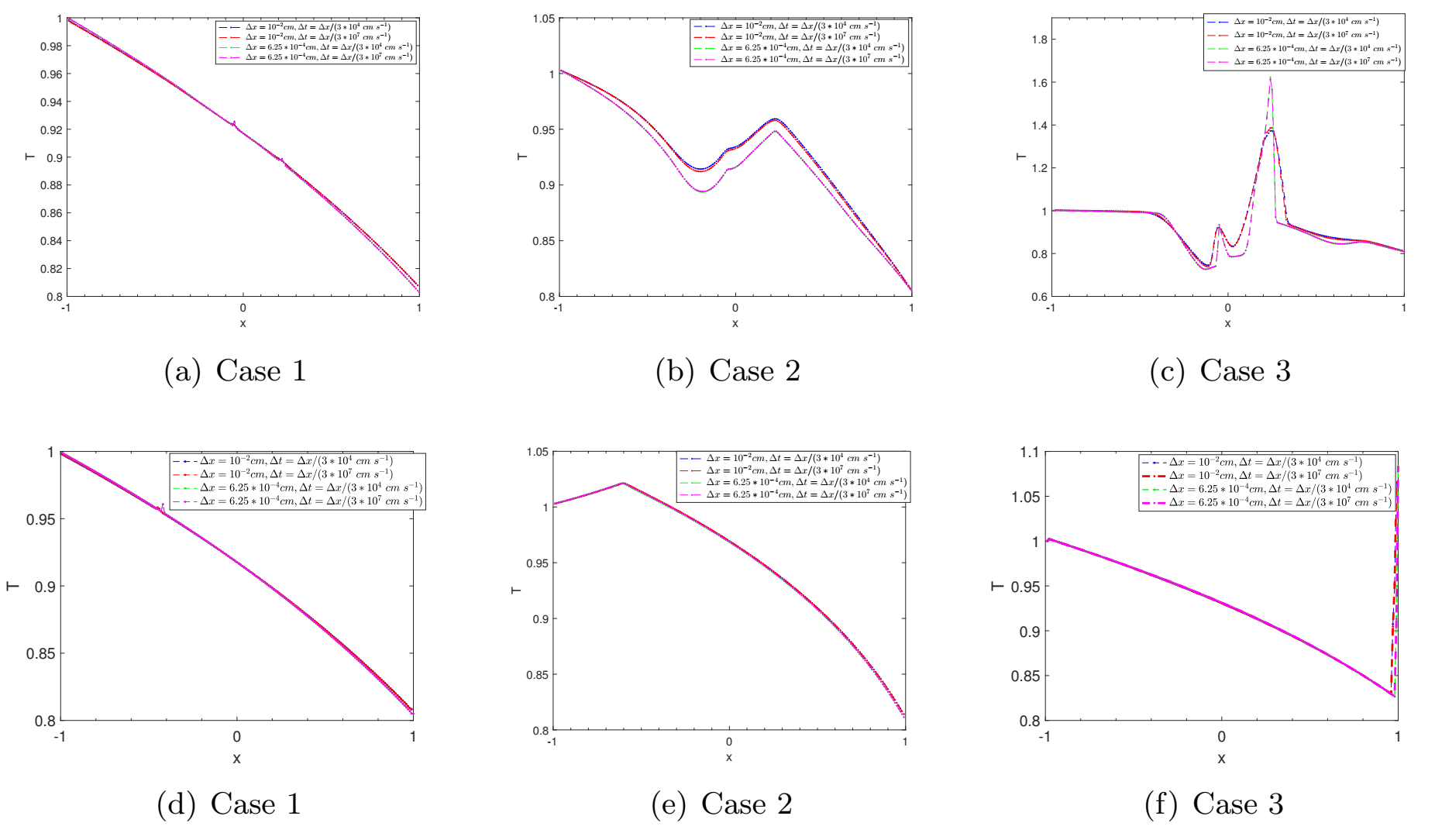}
	\caption{Comparison of numerical results of Case 1, 2 and 3 calculated by our AP solver with $\Delta x=10^{-2}\ cm$ and two different $\Delta t=\Delta x/(3*10^7\ cm\ s^{-1})$, $\Delta x/(3*10^4\ cm\ s^{-1})$, $\Delta x=6.25*10^{-4}\ cm$ and two different $\Delta t=\Delta x/(3*10^7\ cm\ s^{-1})$, $\Delta x/(3*10^4\ cm\ s^{-1})$. In (a),(b) and (c), the computed time $t=6.78*10^{-5}s$, in (d),(e) and (f), the computed time $t=3.33*10^{-3}s$. }
	\label{fig.3-11}
\end{figure}

We show that the convergence order and stability of our scheme for the coupled RMHD system \eqref{rtmhd}-\eqref{eqn:003} are independent of the light speed. In this example we assume that the speed of light $C=3\times 10^{10}\ cm\ s^{-1}$,  the radiation constant $a_r=0.0001\ Jk \ cm^{-3} \ keV^{-4}$. Three different sets of $\sigma_a$, $\sigma_s$ are tested: 1) $\sigma_a=1\ cm^{-1}$, $\sigma_s=1\ cm^{-1}$; 2) $\sigma_a=10^3\ cm^{-1}$, $\sigma_s=10^{-3}\ cm^{-1}$; 3) $\sigma_a=10^{-3}\ cm^{-1}$, $\sigma_s=10^{3}\ cm^{-1}$. Moreover, the initial data for the three cases are the same as in Example 1 and the units of density, magnetic flux density and pressure are respectively $g\ cm^{-3}$, ${\rm T}$ and $Pa$.
 Since the exact solutions are not known, the numerical errors are defined by
\begin{equation} \label{error}
\begin{aligned}
&\text{error}_\rho = \| \rho_{\Delta x}(\cdot, t_\text{max}) -\rho_{\Delta x/2}(\cdot, t_\text{max})\|_{l_1},\quad
\text{error}_p = \| p_{\Delta x}(\cdot, t_\text{max}) - p_{\Delta x/2}(\cdot, t_\text{max})\|_{l_1},\\
&\text{error}_{v_x} = \| (v_x)_{\Delta x}(\cdot, t_\text{max}) - (v_x)_{\Delta x/2}(\cdot, t_\text{max})\|_{l_1},\quad
\quad
\text{error}_{v_y} = \|  (v_y)_{\Delta x}(\cdot, t_\text{max}) -  (v_y)_{\Delta x/2}(\cdot, t_\text{max})\|_{l_1},\\
&\text{error}_{B_y} = \| (B_y)_{\Delta x}(\cdot, t_\text{max}) - (B_y)_{\Delta x/2}(\cdot, t_\text{max})\|_{l_1},\quad
\quad
\text{error}_{T} = \| T_{\Delta x}(\cdot, t_\text{max}) - T_{\Delta x/2}(\cdot, t_\text{max})\|_{l_1}.
\end{aligned}
\end{equation}
Different $\Delta x=10^{-2}\ cm$, $5*10^{-3}\ cm$, $2.5*10^{-3}\ cm$, $1.25*10^{-3}\ cm$, $6.25*10^{-4}\ cm$ are tested and $\Delta t$ are chosen to be $\Delta x/(3*10^7\ cm\ s^{-1})$ and $\Delta x/(3*10^4\ cm\ s^{-1})$. \textcolor{black}{The convergence orders are displayed in Figure. \ref{fig:er1-case} and we can see that the convergence orders are independent of the time step. And we plot the $L^2$ norm of $\average{Q}$ and $\average{nQ}$ in each time step with $\Delta x=2.5*10^{-3}\ cm$ and $\Delta t=\Delta x/(3*10^4\ cm\ s^{-1})$. We can see from Fig.\ref{fig.4.2.0}  that the values of 	$\average{Q}$	and	$\average{nQ}$ are close to zero. The two constraints $\average{Q}=0$ and	$\average{nQ}=0$ are not exactly satisfied due to numerical errors. However, we can consider of solving a new model system composed of $J$, $R$, $Q$ and the MHD system,  which can yields the solution to the original RMHD system.  The results are acceptable except some numerical errors.} Moreover, the CFL condition of explicit schemes requires that $\Delta t<\Delta x/C=\Delta x/(3*10^{10}\ cm\ s^{-1})$, but much larger time steps are allowed in our AP solver.  
To check the stability, we plot the numerical results of Case 1, 2 and 3 calculated by our AP solver with $\Delta x=10^{-2}\ cm$ and two different $\Delta t=\Delta x/(3*10^7\ cm\ s^{-1})$, $\Delta x/(3*10^4\ cm\ s^{-1})$, $\Delta x=6.25*10^{-4}\ cm$ and two different $\Delta t=\Delta x/(3*10^7\ cm\ s^{-1})$, $\Delta x/(3*10^4\ cm\ s^{-1})$. In Figure \ref{fig.3-11}, we can observe that the stability does not depend on the time step but the accuracy does.


\subsection{Example 5}
\begin{figure}[htbp]
	\centering
	\subfigure[Scattering coefficient $\sigma_s$]{
		\includegraphics[width=2.4in]{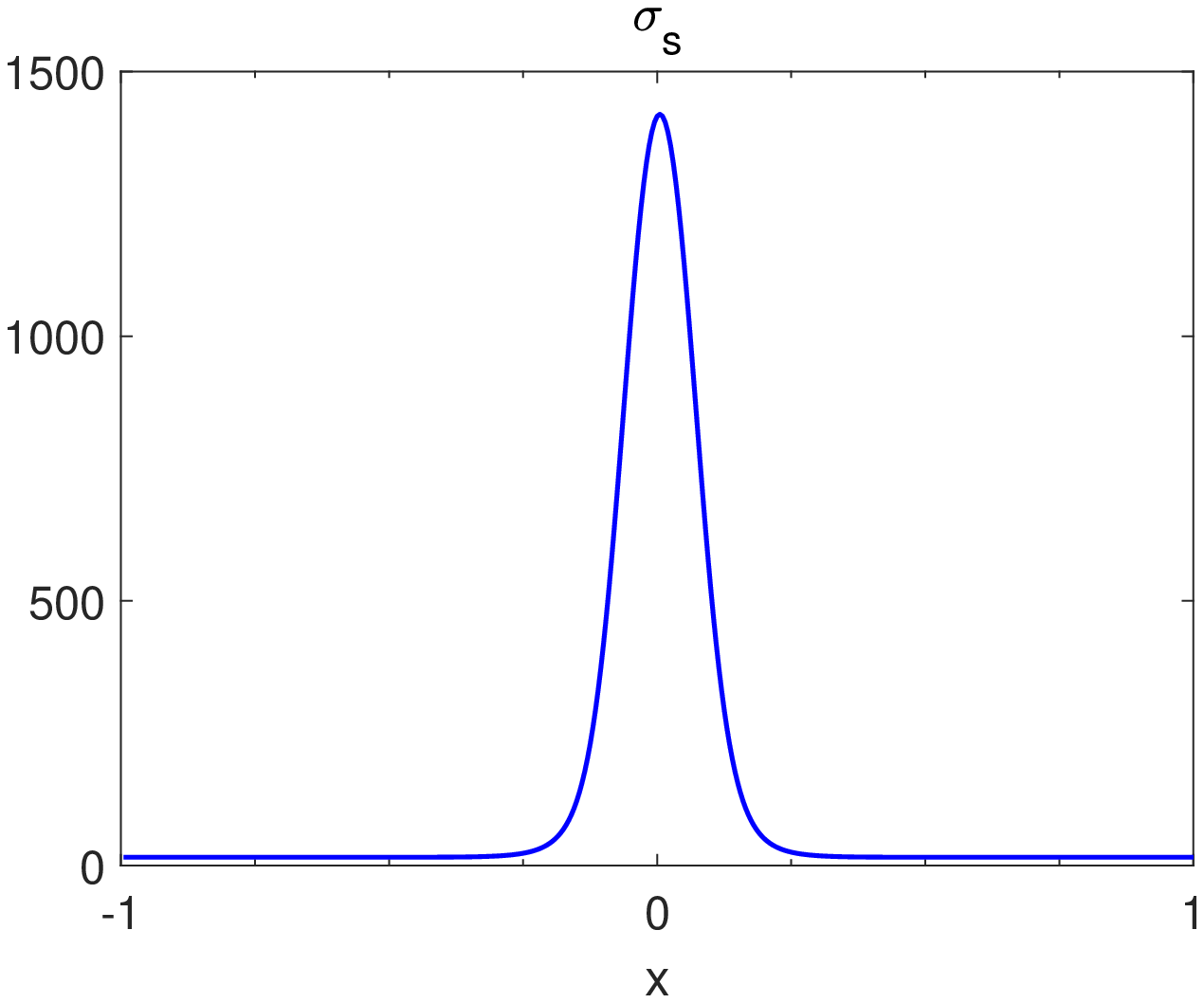}\label{fig.ex4-1}
	}
	\subfigure[Residual term $Q$]{
		\includegraphics[width=2.5in]{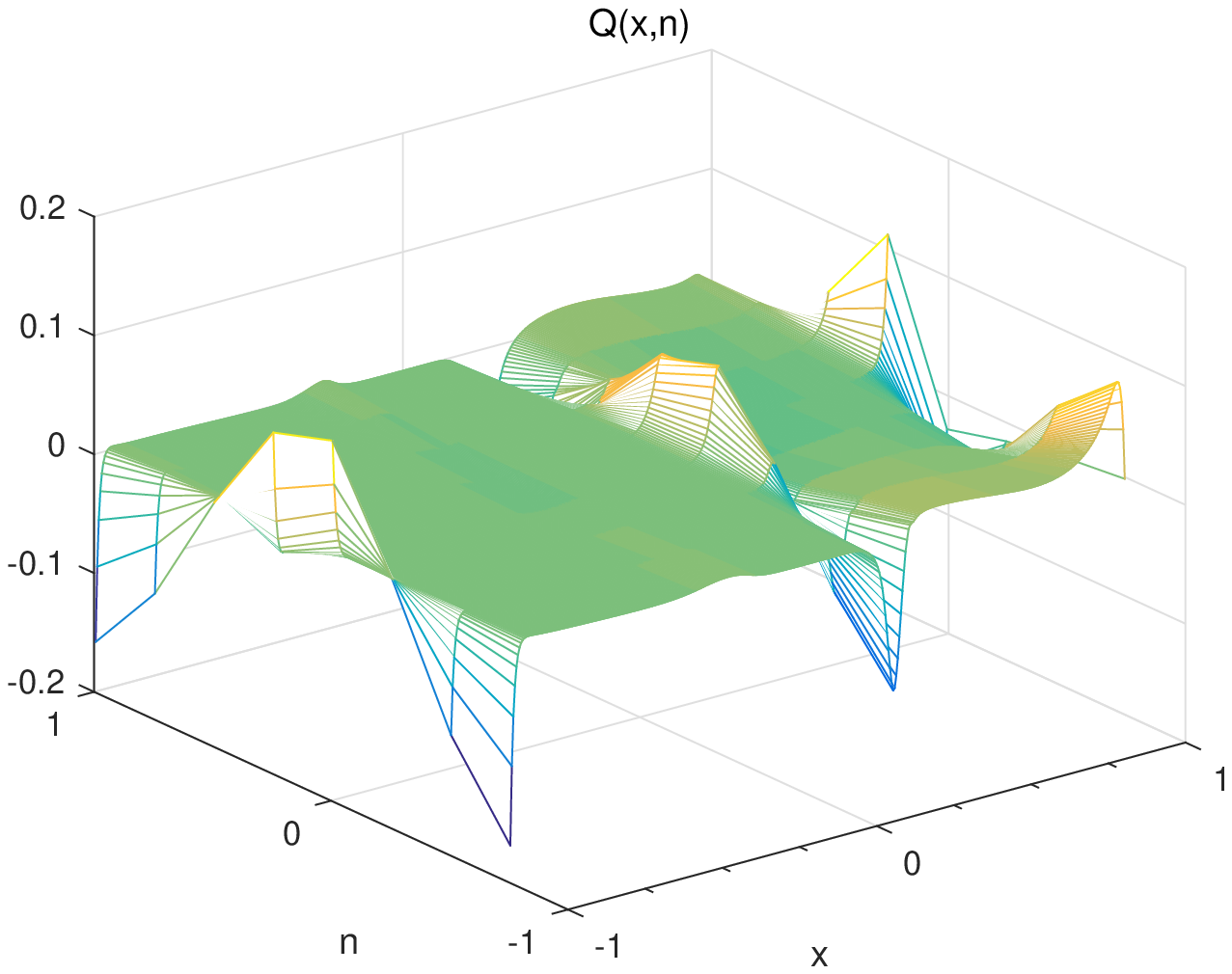}
	}
	\quad    
	\subfigure[Density $\rho$]{
		\includegraphics[width=2.5in]{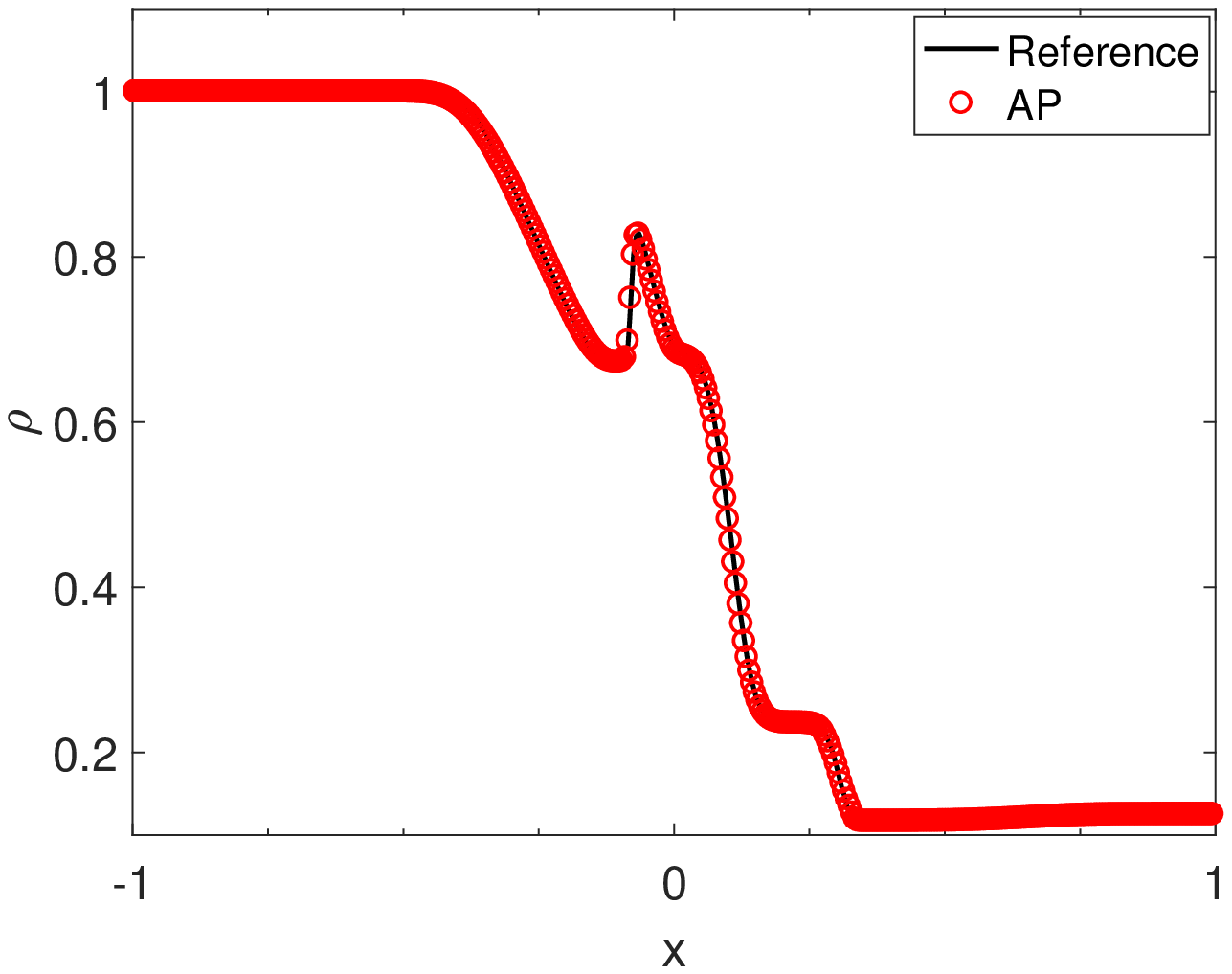}
	}
	\subfigure[Velocity in x-direction  $v_x$]{
		\includegraphics[width=2.5in]{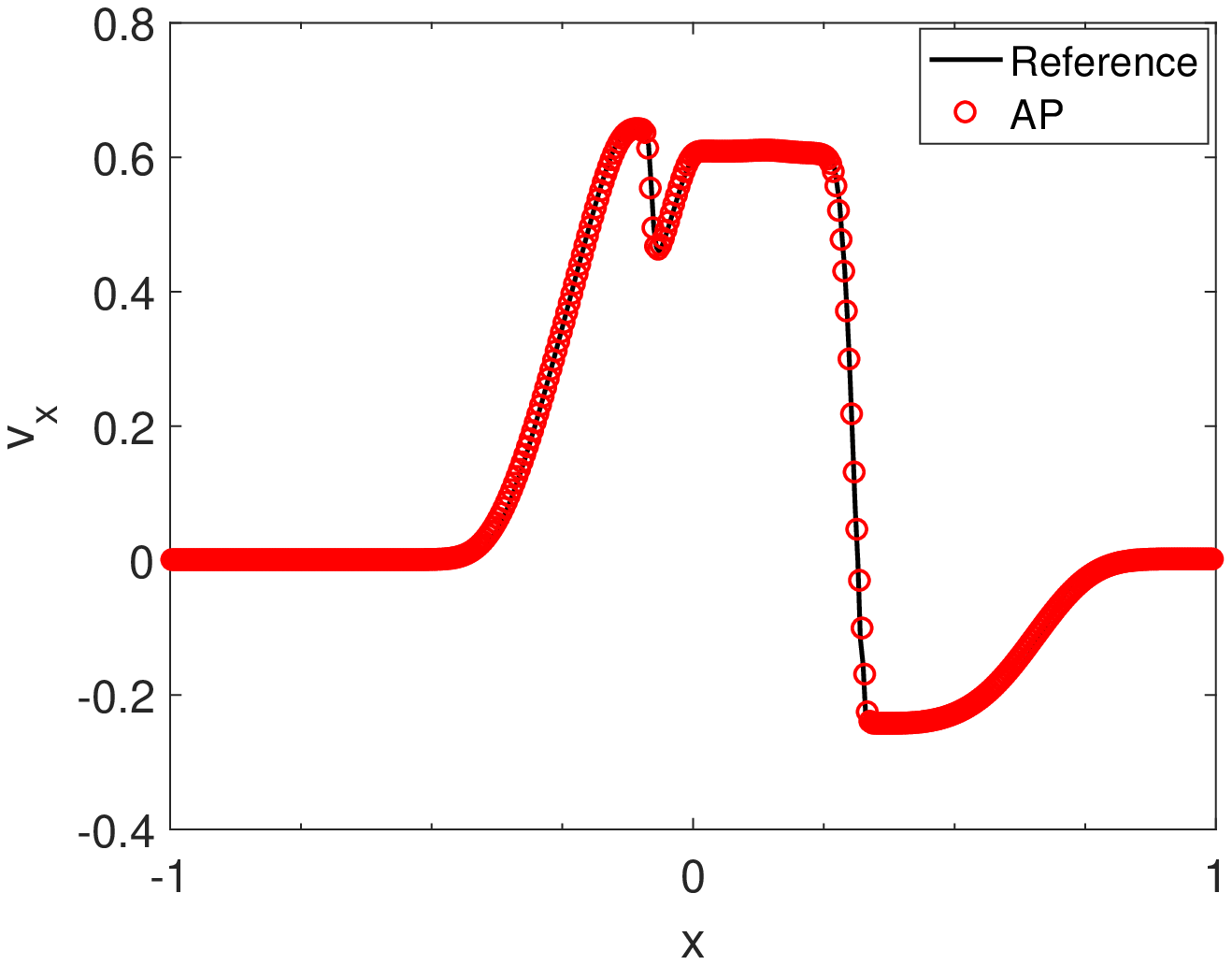}
	}
	\quad
	\subfigure[Velocity in x-direction  $v_x$]{
		\includegraphics[width=2.5in]{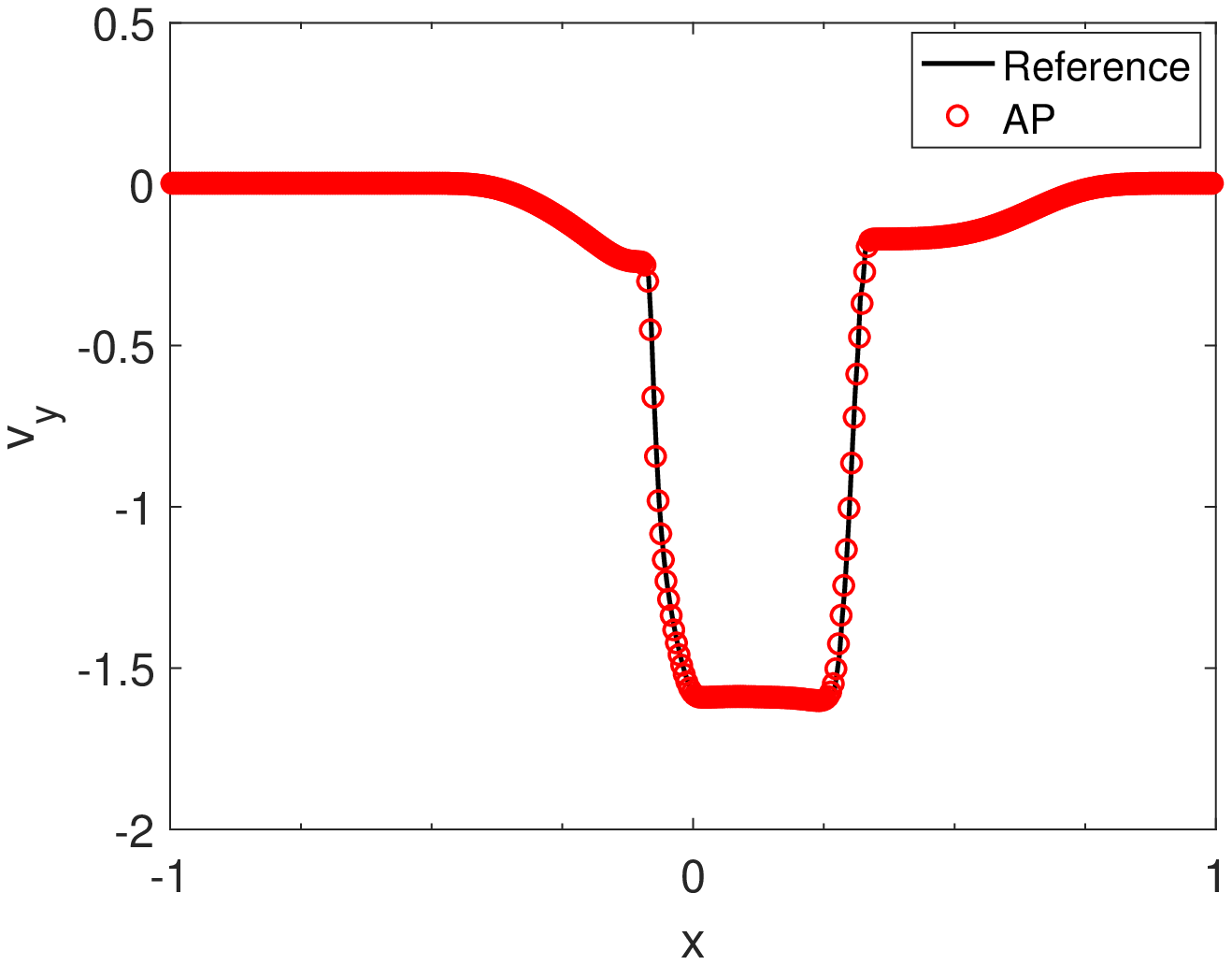}
	}
	\subfigure[Magnetic field in y-direction  $B_y$]{
		\includegraphics[width=2.5in]{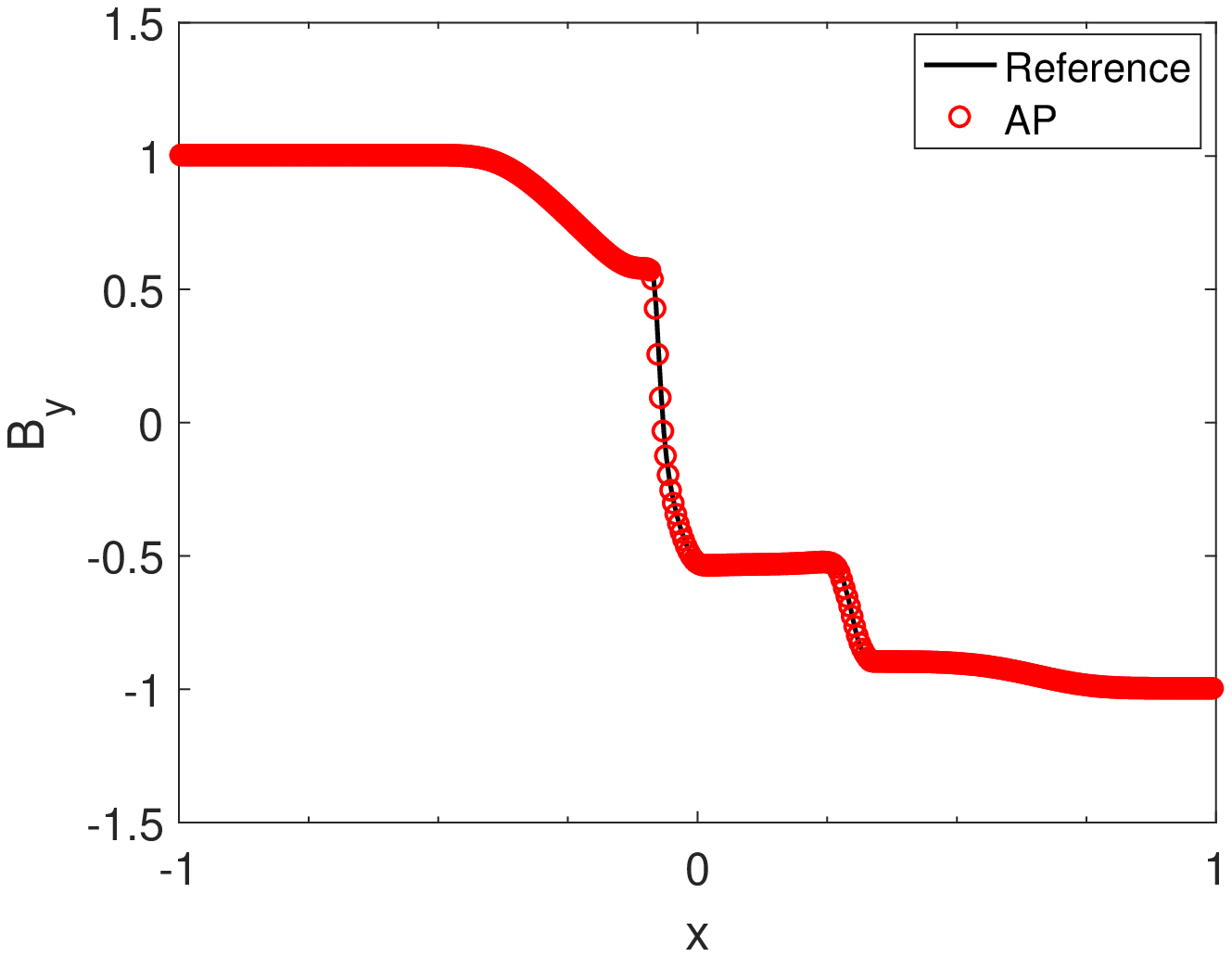}}
	\quad
	\subfigure[Fluid temperature $p/\rho$]{
		\includegraphics[width=2.5in]{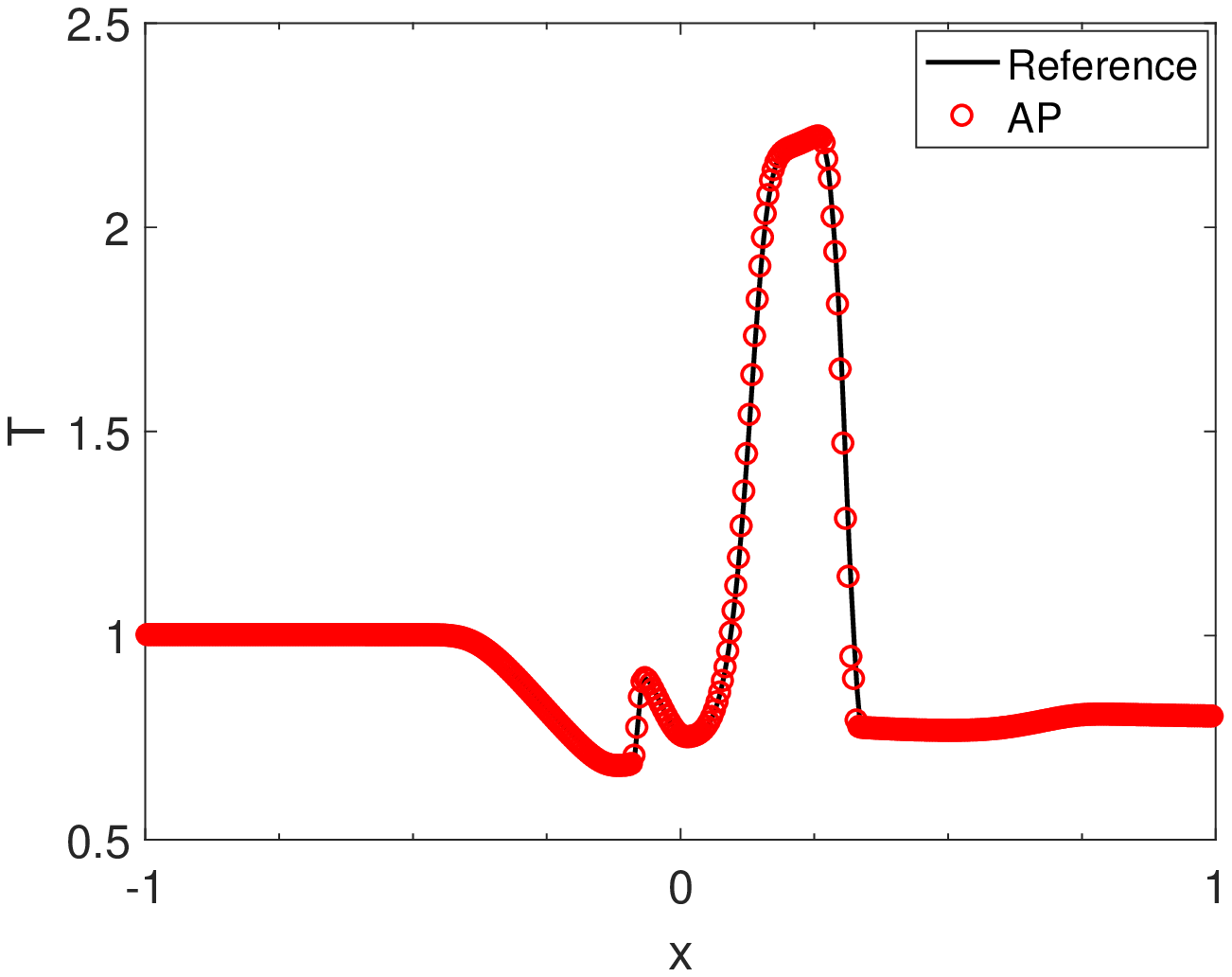}
	}
	\subfigure[Radiation temperature $T_r$]{
		\includegraphics[width=2.5in]{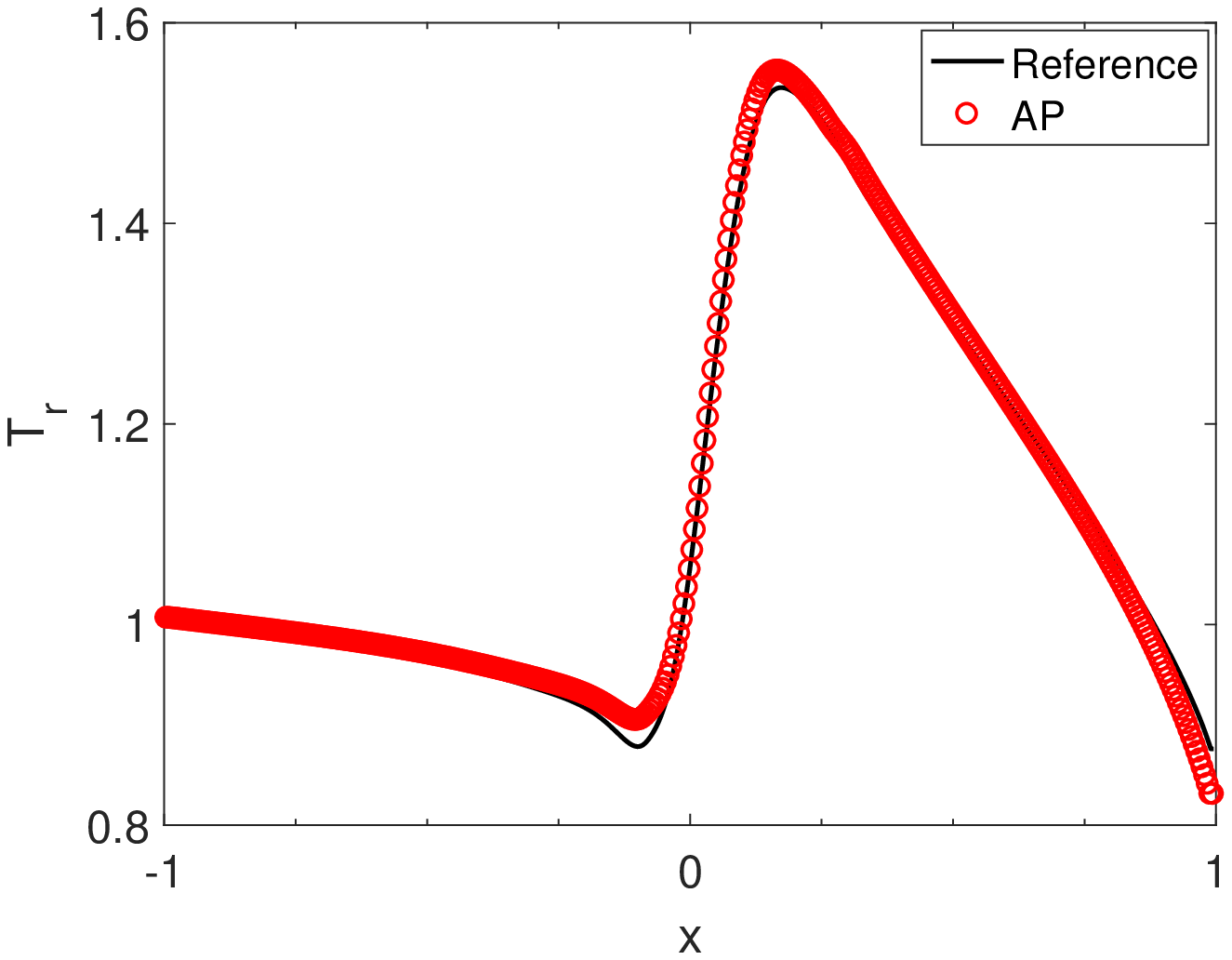}
	}
\caption{Example 5. Comparison of the results at time $t = 6.67*10^{-8} s$ of our AP scheme and  the reference solution.
Here $\Delta x = 1/400\ cm $, $\Delta t = 1.67* 10^{-10} s$  
for  AP schemes, $\Delta x = 1/400\ cm$, $\Delta t =3.33*10^{-15}s$ for reference solution.}
	\label{fig.ex4}
\end{figure}

The RMHD system with multi-scale absorption and scattering coefficients are tested in this example. 
The initial and boundary conditions are the same as in case 1 of Example 4. The absorption and scattering coefficients vary in space, such that
$$\sigma_a=1/3\quad cm^{-1},\qquad
\sigma_s=\left[3+10\left({\rm tanh}(1-11x)+{\rm tanh}(1+11x)\right)\right]^{2.5}\quad cm^{-1}.
$$ 
As we can see from the profile of $\sigma_s(x)$ in Fig. \ref{fig.ex4-1}, $\sigma_s$ varies from very small to very large, thus both optical thick and thin regimes coexist. To compute the reference solution, we use a   finite volume solver in \cite{sun2020multiscale} and a fine mesh  $\Delta x = 1/400\ cm$, $\Delta t =3.33*10^{-15}s$. We use $\Delta x = 1/400\ cm $, $\Delta t = 1.67\times 10^{-10} s$  in our AP scheme, which is much larger than the allowed time step for the explicit solver. The numerical results at $t = 6.67*10^{-8} s$ are displayed in Fig. \ref{fig.ex4}(b)-(h). As we can see from Fig. \ref{fig.ex4}(b) that $Q(x,n)$ is not small in some region, which indicate that one has to use RMHD system instead of the limit system.  
The good agreement between our solution and the reference solution indicates that our method works well for the coexisting case for the RMHD system. 

\section{Conclusion and discussions}
In this paper, we have introduced an AP scheme in both space and time for the RMHD system that couples the ideal MHD equations with a gray RTE. Two different scalings are considered, one results in an equilibrium diffusion limit system, while the other results in nonequilibrium. The main idea is to decompose the intensity into three
parts, two parts correspond to the zeroth and first order moments, while the third part is the residual.  The
two macroscopic moments are updated first by treating them implicitly while the residual explicitly. After the zeroth and first order moments are obtained, we update the residual term by solving implicitly a linear transport equation for each direction $\mathbf{n}$. For the space discretization, we use Roe’s
method in the Athena code to solve the convective part in the ideal MHD equations and the UGKS for the
RTE. Numerical results of optically thin and thick regions coexisting case are presented and the radiative shock problem are tested. The stability and convergence orders of the proposed scheme are independnent of the light speed. \textcolor{black}{ In each time step, our method requires nonlinear iterations to solve a system with only macroscopic quantities, which means the nonlinear iteration is independent of the number of angles. After the macroscopic quantities are obtained, decoupled linear transport equations are solved only once in each time step. The computational cost is much lower than a fully implicit solver. The scheme performance for real multi-dimensional anisotropic problem will be our future work.}

It would be worthwhile to construct a  second-order AP method for the radiation magnetohydrodynamics system. In the semi-implicit scheme, a simple time-integration method has been used. To obtain a second scheme, higher order IMEX time-integration method can be considered.
 
 \bigskip
\textbf{Acknowledgement:} The authors would like to thank Dr. Yanfei Jiang from  Computational Center for Astrophysics, the Flatiron Institute, for proposing this problem and useful discussions. S. Jin is partially supported by NSFC12031013 and the Strategic Priority Research Program of Chinese Academy of Sciences, XDA25010404; M. Tang and X. J. Zhang are partially supported by NSFC11871340 and the Strategic Priority Research Program of Chinese Academy of Sciences, XDA25010401.

 \appendix
  \renewcommand{\appendixname}{Appendix~\Alph{section}}
  \section{The derivation of asymptotic analysis for the two regimes}\label{apen01}
  In the non-equilibrium regime, the  RMHD system \eqref{eqn:non} becomes:
\begin{subequations}
\begin{numcases}{}
\frac{\partial I}{\partial t}+\frac{c}{\varepsilon}\boldsymbol{n}\cdot\nabla I=c\sigma_a\left( \f{T^4}{4\pi}-I\right)+\frac{c\sigma_s}{\varepsilon^3}(J-I)+3\eps \boldsymbol{n}\cdot \boldsymbol{v}\sigma_a\left(\f{T^4}{4\pi}-J\right)+\boldsymbol{n}\cdot \boldsymbol{v}\left(\eps\sigma_a+\f{\sigma_s}{\eps^2}\right)(I+3J)\nonumber\\
 \hspace{5cm}-2\f{\sigma_s}{\eps} \boldsymbol{v}\cdot \boldsymbol{H}-(\eps^2\sigma_a-\sigma_s)\f{\boldsymbol{v}\cdot \boldsymbol{v}}{c}J-(\eps^2\sigma_a-\sigma_s)\f{\boldsymbol{v}\cdot (\boldsymbol{v}\cdot {\rm K})}{c}, \\
 \frac{\partial\rho}{\partial t}+\nabla\cdot(\rho \boldsymbol{v})=0,\\
 \frac{\partial(\rho \boldsymbol{v})}{\partial t}+\nabla\cdot(\rho \boldsymbol{v}\boldsymbol{v}-\boldsymbol{BB}+{\rm P^*})=-\PP_0\boldsymbol{S_{rp}},\\
 \frac{\partial E}{\partial t}+\nabla\cdot[(E+P^*)\boldsymbol{v}-\boldsymbol{B}(\boldsymbol{B}\cdot \boldsymbol{v})]=-\frac{c}{\varepsilon}\PP_0S_{re}, \\
   \frac{\partial \boldsymbol{B}}{\partial t}+\nabla\times(\boldsymbol{v}\times \boldsymbol{B})=0 \,,
 \end{numcases}
\end{subequations}
where
\begin{equation*}
S_{re}=\eps\sigma_a\left( T^4-4\pi J\right)+(\eps^3 \sigma_a-\eps\sigma_s)\f{4\pi \boldsymbol{v}}{c^2}\cdot\left[\f{c}{\eps}\boldsymbol{H}-\left(\boldsymbol{v}J+\boldsymbol{v} \cdot {\rm K}\right)\right]	, \
\end{equation*}
\begin{equation*}
\boldsymbol{S_{rp}}=-\frac{4\pi(\sigma_s+\eps^2\sigma_a)}{c}\left[\f{c}{\eps}\boldsymbol{H}-\left(\boldsymbol{v}J+\boldsymbol{v} \cdot {\rm K}\right)\right]+\f{\boldsymbol{v}}{c}\eps^2\sigma_a\left(T^4-4\pi J\right).
\end{equation*}
When $\eps \to 0$,  we show that the solution to (\ref{eqn:004}) can be approximated by the solution to a nonlinear diffusion equation coupled with a MHD system. We assume that the radiation intensity $I$ and temperature $T$ have the following Chapman-Enskog expansion such that
\begin{equation}\label{expan}
\begin{aligned}
I= I^{(0)}+\eps I^{(1)}+\eps ^2 I^{(2)}+\cdots,\\
T= T^{(0)}+\eps T^{(1)}+\eps ^2 T^{(2)}+\cdots.
\end{aligned}
\end{equation}
By substituting  ansatz \eqref{expan} into  equation  (\ref{eqn:case11004}) and collecting the terms of the same order in $\eps$, we have
\begin{subequations} \label{limit-case1}
\begin{align}
 &\order \left(\f{1}{\eps^2}\right): I^{(0)}=J^{(0)},\label{limit-case11}\\
 & \order \left(\f{1}{\eps}\right): c\boldsymbol{n}\cdot\nabla I^{(0)}=c\sigma_s\left( J^{(1)} -I^{(1)}\right)+\boldsymbol{n}\cdot \boldsymbol{v}\sigma_s\left(I^{(0)}+3J^{(0)} \right).\label{limit-case12}
 \end{align}
\end{subequations}
The term $\boldsymbol{H}^{(0)}$ in equation \eqref{limit-case12} equals to 0 due to equation \eqref{limit-case11}. Multiplying both sides of (\ref{eqn:case11004}) by $\boldsymbol{n}$, taking its integral with respect to $\boldsymbol{n}$ and combining the obtained equation with  (\ref{eqn:case13004}), one can get the following momentum conservation equation:
\begin{equation}\label{m-e}
\partial_t\left(\f{\rho \boldsymbol{v}}{\PP_0}+\f{4\pi\eps}{c}\average{\boldsymbol{n}I}\right)+\nabla\cdot\left(\f{\rho \boldsymbol{v}\boldsymbol{v}-\boldsymbol{BB}+{\rm P^*}}{\PP_0}+4\pi\average{ \boldsymbol{n}\boldsymbol{n}I}\right)=0.
\end{equation}
By using (\ref{limit-case11}) and $D_d=\average{\boldsymbol{n}\boldsymbol{n}}=\f{1}{3}I_d$ (where $I_d$ denotes the 3 by 3 identity matrix), when $\eps\to0$, \eqref{m-e} gives
\begin{equation}\label{momen1}
\partial_t(\rho^{(0)} \boldsymbol{v}^{(0)})+\nabla\cdot\left(\rho^{(0)} \boldsymbol{v}^{(0)}\boldsymbol{v}^{(0)}-\boldsymbol{B}^{(0)}\boldsymbol{B}^{(0)}+({\rm P^*})^{(0)}\right)=-4\pi\PP_0D_d\nabla J^{(0)}.
\end{equation}
 By taking the integral with respect to $\boldsymbol{n}$  on  both sides of (\ref{eqn:case11004}), and combining it with  (\ref{eqn:case14004}), one can obtain the energy conservation equation
\begin{equation}\label{e-e}
\partial_t\left(\f{E}{\PP_0}+4\pi\average{I}\right)+\nabla\cdot\left[\f{(E+P^*)\boldsymbol{v}-\boldsymbol{B}(\boldsymbol{B}\cdot \boldsymbol{v})}{\PP_0}+\f{ 4 \pi c\average{\boldsymbol{n}I}}{\eps}\right]=0.
\end{equation}
By using (\ref{limit-case11}), when $\eps\to0$, \eqref{e-e} gives
\begin{equation*}
\partial_t\left(\f{E^{(0)}}{\PP_0}+ 4\pi\average{I^{(0)}}\right)+\nabla\cdot\left[\f{(E^{(0)}+(P^*)^{(0)})\boldsymbol{v}^{(0)}-\boldsymbol{B}^{(0)}(\boldsymbol{B}^{(0)}\cdot \boldsymbol{v}^{(0)})}{\PP_0}+ 4\pi c\average{\boldsymbol{n}I^{(1)}}\right]=0.
\end{equation*}
Noting (\ref{limit-case12}), we find
\begin{equation}\label{en-e}
\partial_t\left(E^{(0)}+4\pi\PP_0J^{(0)}\right)+\nabla\cdot\left[(E^{(0)}+(P^*)^{(0)})\boldsymbol{v}^{(0)}-\boldsymbol{B}^{(0)}(\boldsymbol{B}^{(0)}\cdot \boldsymbol{v}^{(0)})+16\pi\PP_0D_d\boldsymbol{v}^{(0)}J^{(0)}-\f{ 4\pi c\PP_0D_d}{\sigma_s}\nabla J^{(0)}\right]=0.
\end{equation}
\eqref{eqn:case12004}, \eqref{momen1}, \eqref{en-e} and \eqref{eqn:case15004} give four equations for $\rho,\ \rho \boldsymbol{v},\ E,\ J,\ \boldsymbol{B}$, to obtain a closed  system, we need one more equation for $J$. Taking $\average{\cdot}$ on both sides of equation (\ref{eqn:case11004}), yields:
\begin{equation*}
\partial_t\average{I}+\frac{c}{\varepsilon}\nabla\cdot \average{\boldsymbol{n}I}=\frac{c}{4\pi\varepsilon}S_{re}.
\end{equation*}
By using the expansion in \eqref{expan} and (\ref{limit-case11}),  the leading order terms are
\begin{equation*}
4\pi\partial_t\average{I^{(0)}}+4\pi c\nabla\cdot \average{\boldsymbol{n}I^{(1)}}=c\sigma_a\left((T^{(0)})^4-4\pi\average{I^{(0)}}\right)-4\pi\sigma_s\f{\boldsymbol{v}^{(0)}}{c}\cdot\left[c\average{\boldsymbol{n}I^{(1)}}-\boldsymbol{v}^{(0)}\left(\average{I^{(0)}}+\average{\boldsymbol{n}\boldsymbol{n}I^{(0)}}\right)\right],
\end{equation*}
then from (\ref{limit-case12}), one  gets 
\begin{equation}
4\pi\partial_tJ^{(0)}+\nabla\cdot\left(16\pi D_d\boldsymbol{v}^{(0)}J^{(0)}-\f{ 4\pi c D_d}{\sigma_s}\nabla J^{(0)}\right)=c\sigma_a\left((T^{(0)})^4-4\pi J^{(0)}\right)+4\pi D_d\boldsymbol{v}^{(0)}\cdot\nabla J^{(0)}.
\end{equation}
Therefore, when $\eps\to 0$ in (\ref{eqn:004}), the solution can be approximated by the solution of the following non-equilibrium system:
\begin{subequations}
\begin{numcases}{}
\partial_t\rho +\nabla\cdot(\rho \boldsymbol{v})=0, \\
 \partial_t(\rho \boldsymbol{v})+\nabla\cdot\left(\rho \boldsymbol{v}\boldsymbol{v}-\boldsymbol{B}\boldsymbol{B}+{\rm P^*}\right)=-4\pi\PP_0D_d\nabla J,\\
 \partial_t\left(E+4\pi\PP_0J\right)+\nabla\cdot\left[(E+P^*)\boldsymbol{v}-\boldsymbol{B}(\boldsymbol{B}\cdot \boldsymbol{v})+16\pi\PP_0D_d\boldsymbol{v}J\right]=\nabla\cdot\left(\f{4\pi c\PP_0D_d}{\sigma_s}\nabla J\right),\\
 4\pi\partial_tJ+\nabla\cdot\left(16\pi D_d\boldsymbol{v}J-\f{ 4\pi cD_d}{\sigma_s}\nabla J\right)=c\sigma_a\left(T^4-4\pi J\right)+4\pi D_d\boldsymbol{v}\cdot\nabla J, \\
   \partial_t \boldsymbol{B}+\nabla\times(\boldsymbol{v}\times \boldsymbol{B})=0. \,
 \end{numcases}
\end{subequations}
In the equilibrium regime, the  RMHD system \eqref{eqn:non} becomes:
\begin{subequations}
\begin{numcases}{}
\frac{\partial I}{\partial t}+\frac{c}{\varepsilon}\boldsymbol{n}\cdot\nabla I=\frac{c\sigma_a}{\varepsilon^2}\left( \f{T^4}{4\pi}-I\right)+c\sigma_s(J-I)+3\boldsymbol{n}\cdot \boldsymbol{v}\f{\sigma_a}{\eps}\left(\f{T^4}{4\pi}-J\right)+\boldsymbol{n}\cdot \boldsymbol{v}\left(\f{\sigma_a}{\eps}+\eps\sigma_s\right)(I+3J)\nonumber\\
\hspace{5cm}-2\eps\sigma_s \boldsymbol{v}\cdot \boldsymbol{H}-(\sigma_a-\eps^2\sigma_s)\f{\boldsymbol{v}\cdot \boldsymbol{v}}{c}J-(\sigma_a-\eps^2\sigma_s)\f{\boldsymbol{v}\cdot (\boldsymbol{v}\cdot {\rm K})}{c}, \\
 \frac{\partial\rho}{\partial t}+\nabla\cdot(\rho \boldsymbol{v})=0,\\
 \frac{\partial(\rho \boldsymbol{v})}{\partial t}+\nabla\cdot(\rho \boldsymbol{v}\boldsymbol{v}-\boldsymbol{B}\boldsymbol{B}+{\rm P^*})=-\PP_0\boldsymbol{S_{rp}},\\
 \frac{\partial E}{\partial t}+\nabla\cdot[(E+P^*)\boldsymbol{v}-\boldsymbol{B}(\boldsymbol{B}\cdot \boldsymbol{v})]=-\frac{c}{\varepsilon}\PP_0S_{re}, \\
   \frac{\partial \boldsymbol{B}}{\partial t}+\nabla\times(\boldsymbol{v}\times \boldsymbol{B})=0 \,,
 \end{numcases}
\end{subequations}
where
\begin{equation*}
S_{re}=\frac{\sigma_a}{\varepsilon}\left( T^4-4\pi J\right)+(\eps\sigma_a-\eps^3\sigma_s)\f{4\pi \boldsymbol{v}}{c^2}\cdot\left[\f{c}{\eps}\boldsymbol{H}-\left(\boldsymbol{v}J+\boldsymbol{v} \cdot {\rm K}\right)\right],
\end{equation*}
\begin{equation*}
\boldsymbol{S_{rp}}=-\frac{4\pi(\eps^2\sigma_s+\sigma_a)}{c}\left[\f{c}{\eps}\boldsymbol{H}-\left(\boldsymbol{v}J+\boldsymbol{v} \cdot {\rm K}\right)\right]+\f{\boldsymbol{v}}{c}\sigma_a\left(T^4-4\pi J\right).
\end{equation*}
Using the same expansion as in \eqref{expan},  one has
\begin{equation}\label{t4}
\begin{aligned}
T^4&= (T^{(0)}+\eps T^{(1)}+\eps ^2 T^{(2)}+\cdots)^4\\
&=(T^{(0)})^4+4\eps(T^{(0)})^3T^{(1)}+\cdots.
\end{aligned}
\end{equation}
Then  substituting \eqref{t4} into  (\ref{eqn:case11eqb}) and collecting the terms of the same order of $\eps$, one gets 
\begin{subequations} \label{limit-case21}
	\begin{align}
	&\order \left(\f{1}{\eps^2}\right): I^{(0)}=\f{(T^{(0)})^4}{4\pi},\label{limit-case211}\\
	& \order \left(\f{1}{\eps}\right): c\boldsymbol{n}\cdot\nabla I^{(0)}=c\sigma_a\left( \f{4(T^{(0)})^3T^{(1)}}{4\pi} -I^{(1)}\right)+\boldsymbol{n}\cdot \boldsymbol{v}\sigma_a\left(I^{(0)}+\f{3(T^{(0)})^4}{4\pi}\right),\label{limit-case212}
	\end{align}
\end{subequations}
by similar calculations  in section \ref{apen01}, one can get the following equilibrium diffusion limit:
\begin{subequations}
\begin{numcases}{}
\partial_t\rho +\nabla\cdot(\rho \boldsymbol{v})=0, \\
\partial_t(\rho \boldsymbol{v})+\nabla\cdot\left(\rho \boldsymbol{v}\boldsymbol{v}-\boldsymbol{B}\boldsymbol{B}+{\rm P^*}\right)=-\PP_0D_d\nabla T^4,\\
 \partial_t\left(E+\PP_0T^4\right)+\nabla\cdot\left[(E+P^*)\boldsymbol{v}-\boldsymbol{B}(\boldsymbol{B}\cdot \boldsymbol{v})+4\PP_0D_d\boldsymbol{v}T^4\right]=\nabla\cdot\left(\f{ c\PP_0D_d}{\sigma_a}\nabla T^4\right),\\
   \partial_t \boldsymbol{B}+\nabla\times(\boldsymbol{v}\times \boldsymbol{B})=0. \,
 \end{numcases}
\end{subequations}

	\section{The processes of UGKS}\label{apen1}
Recall the RTE \eqref{orig1} in the slab geometry. Integrating it over $[t_s,t_{s+1}]$ and $[x_{i-\f{1}{2}},x_{1+\f{1}{2}}]$ gives
\begin{equation*}
	\f{I_i^{s+1}-I_i^{s}}{\Delta t}+\f{1}{\Delta x}(\zeta_{i+\f{1}{2}}-\zeta_{i-\f{1}{2}})=\mathcal{C}\mathscr{L}_a\sigma_{a,i}\left( \f{(T_i^{s+1})^4}{4\pi}-I_i^{s+1}\right)+\mathcal{C}\mathscr{L}_s\sigma_{s,i}(J_i^{s+1}-I_i^{s+1})+G.\\
\end{equation*}
The microscopic flux $\zeta$ is computed by solving the following initial value problem at the cell boundary $x=x_{i+\f{1}{2}}$:
	\begin{equation}\label{ibdy}
	\left\{
	\begin{aligned}
	&\f{I_i^{s+1}-I_i^{s}}{\Delta t}+\f{1}{\Delta x}(\zeta_{i+\f{1}{2}}-\zeta_{i-\f{1}{2}})=\mathcal{C}\mathscr{L}_a\sigma_{a,i}\left( \f{(T_i^{s+1})^4}{4\pi}-I_i^{s+1}\right)+\mathcal{C}\mathscr{L}_s\sigma_{s,i}(J_i^{s+1}-I_i^{s+1})+G,\\
	&I(t,x,n)|_{t=t_s}=I(t_s,x,n).
	\end{aligned}
	\right.
	\end{equation}
	In the UGKS scheme, we firstly assume the coefficients $\sigma_a$ and $\sigma_s$ in space and time are piecewise constant, thus the radiative transfer equation in  RMHD system \eqref{orig} is equivalent to
	\begin{equation}\label{equal}
	\f{d}{dt}e^{\mu t}I\left(t,x+\mathcal{C}nt,n\right)=e^{\mu t}\left(\mathcal{C}\mathscr{L}_s\sigma_s J\left(t,x+\mathcal{C}nt,n\right)+\f{\mathcal{C}\mathscr{L}_a\sigma_a}{4\pi}T^4\left(t,x+\mathcal{C}nt,n\right)+G\right),
	\end{equation}
	with $\mu=\mathcal{C}\mathscr{L}_a\sigma_a+ \mathcal{C}\mathscr{L}_s\sigma_s$. Solving the equivalent equation \eqref{equal} with the initial value in the system \eqref{ibdy}, which gives the solution
	\begin{equation}\label{isol}
	\begin{aligned}
	I&(t,x_{i+\f{1}{2}},n)\approx e^{-\mu_{i+\f{1}{2}}(t-t_s)}I\left(t_s,x_{i+\f{1}{2}}-\mathcal{C}n(t-t_s)\right)+\f{1-e^{-\mu_{i+\f{1}{2}}(t-t_s)}}{\mu_{i+\f{1}{2}}}\hat{G}_{i+\f{1}{2}}\\
	&+\int\limits^t_{t_s}e^{-\mu_{i+\f{1}{2}}(t-t_s)}\left(\mathcal{C}\mathscr{L}_s\sigma_{s,i+\f{1}{2}}J\left(z,x_{i+\f{1}{2}}-\mathcal{C}n(t-z)\right)+\f{\mathcal{C}\mathscr{L}_a\sigma_{a,i+\f{1}{2}}}{4\pi}T^4\left(z,x_{i+\f{1}{2}}-\mathcal{C}n(t-z)\right)\right)dz,
	\end{aligned}
	\end{equation}
	with $\mu_{i+\f{1}{2}}$, $\sigma_{a,i+\f{1}{2}}$ and $\sigma_{s,i+\f{1}{2}}$ being the constant value at the corresponding cell boundary for $\mu$, $\sigma_a$ and $\sigma_s$ respectively. In solution \eqref{isol}, there remains three terms that need to approximate: the first one is the initial condition $I(t_s)$ around $x_{i+\f{1}{2}}$, namely the function $I\left(t_s,x_{i+\f{1}{2}}+\mathcal{C}n(t-t_s)\right)$;  the second one is the two functions $J, T^4$ localized  in the time interval $[t_s,t_{s+1}]$ and around the boundary $x_{i+\f{1}{2}}$, i.e., the functions $J\left(z,x_{i+\f{1}{2}}-\mathcal{C}n(t-z)\right)$ and $T^4\left(z,x_{i+\f{1}{2}}-\mathcal{C}n(t-z)\right)$;  the last one is the term $\hat{G}_{i+\f{1}{2}}$.
	
	\subsection*{The first term:}For the first term, a  piecewise constant  reconstruction function is used to approximate the initial function $I(t_s)$ around $x_{i+\f{1}{2}}$:
	\begin{equation}\label{equ:i0}
	I(t_s,x,n)=\left\{
	\begin{aligned}
	&I^s_{i}, \quad {\rm if}\quad x<x_{i+\f{1}{2}}, \\
	&I^s_{i+1}, \quad {\rm if}\quad  x>x_{i+\f{1}{2}}.
	\end{aligned}
	\right.
	\end{equation}
	
	\subsection*{The second term:}For the second term, the approximations of two functions $J$ and $T^4$ between $t_s$ and $t_{s+1}$ and around $x_{i+\f{1}{2}}$ are constructed in the same manner, so only the approximation of  $J$ is introduced in detail. The function $J$ is treated implicitly in time by using piecewise linear polynomials, so the reconstruction for $J$ reads:
	\begin{equation}\label{equ:j}
	J(t,x)=\left\{
	\begin{aligned}
	&J^{s+1}_{i+\f{1}{2}}+\delta_xJ^{s+1,+}_{i+\f{1}{2}}(x-x_{i+\f{1}{2}}), \quad {\rm if}\quad x<x_{i+\f{1}{2}}, \\
	&J^{s+1}_{i+\f{1}{2}}+\delta_xJ^{s+1,-}_{i+\f{1}{2}}(x-x_{i+\f{1}{2}}), \quad {\rm if}\quad  x>x_{i+\f{1}{2}},
	\end{aligned}
	\right.
	\end{equation}
	with the interface value of $J$  defined by
	$$
	J^{s+1}_{i+\f{1}{2}}=\average{I^{s+1}_{i+\f{1}{2}}(n)}=\f{1}{2}\average{I^{s+1}_i+I^{s+1}_{i+1}},
	$$
	and with the spatial derivatives given by
	$$
	\delta_xJ^{s+1,+}_{i+\f{1}{2}}=\f{J^{s+1}_{i+\f{1}{2}}-J^{s+1}_{i}}{\Delta x/2}, \quad \delta_xJ^{s+1,-}_{i+\f{1}{2}}=\f{J^{s+1}_{i+1}-J^{s+1}_{i+\f{1}{2}}}{\Delta x/2}.
	$$
	\subsection*{The third term:}At last, the approximation for the term $\hat{G}$ is discussed. Recall the definition of $G_i$:
$$
\begin{aligned}
G_i=&3\mathscr{L}_a\sigma_{a,i}nv_{x,i}^{s+1}\left(\f{(T^{s+1}_i)^4}{4\pi}-J_{i}^{s+1}\right)+n v_{x,i}^{s+1}\left(\mathscr{L}_a\sigma_{a,i}+\mathscr{L}_s\sigma_{s,i}\right)\left(4J^{s+1}_{i}+nR_{i}^{s+1}+Q_{i}^{s}\right)\\
&\hspace{2cm}-\f{2}{3}\mathscr{L}_s\sigma_{s,i} v_{x,i}^{s+1} R_{i}^{s+1}-\f{(\mathscr{L}_a\sigma_{a,i}-\mathscr{L}_s\sigma_{s,i})(v_{x,i}^{s+1})^2}{\mathcal{C}}\left(\f{4}{3}J^{s+1}_{i}+ K_{Q,i}^{s}\right),
\end{aligned}
$$
	The upwind scheme is used to determine $\hat{G}$, i.e., the  value of $\hat{G}$ at the center of the interval $[x_i,x_{i+1}]$ is determined by the  velocity field $v$ at the both sides boundary, 
	\begin{equation}\label{equ:g}
	\left\{
	\begin{aligned}
	&	\text { If } v^{s+1}_{x,i} \geq v^{s+1}_{x,i+1}: \\
	&\hat{G}_{i+\f{1}{2}}=G_i,\quad {\rm if} \ \left(G_{i+1}-G_{i}\right) /\left(v^{s+1}_{x,i+1}-v^{s+1}_{x,i}\right)>0,\\
	&\hat{G}_{i+\f{1}{2}}=G_{i+1},\quad {\rm if} \ \left(G_{i+1}-G_{i}\right) /\left(v^{s+1}_{x,i+1}-v^{s+1}_{x,i}\right)\leq 0; \\
	&	\text { If } v^{s+1}_{x,i} < v^{s+1}_{x,i+1}: \\
	&\hat{G}_{i+\f{1}{2}}=G_i,\quad {\rm if} \ v^{s+1}_{x,i}>0,\\
	&\hat{G}_{i+\f{1}{2}}=G_{i+1},\quad {\rm if} \ v^{s+1}_{x,i+1}<0,\\
	&\hat{G}_{i+\f{1}{2}}=0,\quad {\rm if} \ v^{s+1}_{x,i} \leq 0 \leq v^{s+1}_{x,i+1}.
	\end{aligned}
	\right.
	\end{equation}

\bibliography{RTMHD_reference}

\begin{bibdiv}
\begin{biblist}

\bib{BOLDING2017511}{article}{
      author={Bolding, Simon},
      author={Hansel, Joshua},
      author={Edwards, Jarrod~D.},
      author={Morel, Jim~E.},
      author={Lowrie, Robert~B.},
       title={Second-order discretization in space and time for
  radiation-hydrodynamics},
        date={2017},
        ISSN={0021-9991},
     journal={Journal of Computational Physics},
      volume={338},
       pages={511 \ndash  526},
  url={http://www.sciencedirect.com/science/article/pii/S0021999117301663},
}

\bib{brio1988upwind}{article}{
      author={Brio, Moysey},
      author={Wu, Cheng~Chin},
       title={An upwind differencing scheme for the equations of ideal
  magnetohydrodynamics},
        date={1988},
     journal={Journal of computational physics},
      volume={75},
      number={2},
       pages={400\ndash 422},
}

\bib{BCDB1}{article}{
      author={Buet, Christophe},
      author={Despres, Bruno},
       title={Asymptotic analysis of fluid models for the coupling of radiation
  and hydrodynamics},
        date={2004},
     journal={Journal of Quantitative Spectroscopy and Radiative Transfer},
      volume={85},
      number={3-4},
       pages={385\ndash 418},
}

\bib{BCDB2}{article}{
      author={Buet, Christophe},
      author={Despres, Bruno},
       title={Asymptotic preserving and positive schemes for radiation
  hydrodynamics},
        date={2006},
     journal={Journal of Computational Physics},
      volume={215},
      number={2},
       pages={717\ndash 740},
}

\bib{2018Tang}{article}{
      author={Chen, Hongfei},
      author={Chen, Gaoyu},
      author={Hong, Xiang},
      author={Gao, Hao},
      author={Tang, Min},
       title={Uniformly convergent scheme for rte with anisotropic scattering
  up to the boundary and interface layers},
        date={2018},
     journal={Communication in Computational Physics},
      volume={24},
       pages={1021\ndash 1048},
}

\bib{filbet2010class}{article}{
      author={Filbet, Francis},
      author={Jin, Shi},
       title={A class of asymptotic-preserving schemes for kinetic equations
  and related problems with stiff sources},
        date={2010},
     journal={Journal of Computational Physics},
      volume={229},
      number={20},
       pages={7625\ndash 7648},
}

\bib{godillon2005coupled}{article}{
      author={Godillon-Lafitte, Pauline},
      author={Goudon, Thierry},
       title={A coupled model for radiative transfer: Doppler effects,
  equilibrium, and nonequilibrium diffusion asymptotics},
        date={2005},
     journal={Multiscale Modeling \& Simulation},
      volume={4},
      number={4},
       pages={1245\ndash 1279},
}

\bib{golse1999convergence}{article}{
      author={Golse, Fran{\c{c}}ois},
      author={Jin, Shi},
      author={Levermore, C~David},
       title={The convergence of numerical transfer schemes in diffusive
  regimes i: Discrete-ordinate method},
        date={1999},
     journal={SIAM journal on numerical analysis},
      volume={36},
      number={5},
       pages={1333\ndash 1369},
}

\bib{hayes}{article}{
      author={Hayes, John~C},
      author={Norman, Michael~L},
       title={Beyond flux-limited diffusion: parallel algorithms for
  multidimensional radiation hydrodynamics},
        date={2003},
     journal={The Astrophysical Journal Supplement Series},
      volume={147},
      number={1},
       pages={197},
}

\bib{pub.1005319601}{article}{
      author={Jiang, Yan-Fei},
      author={Stone, James~M.},
      author={Davis, Shane~W.},
       title={A godunov method for multidimensional radiation
  magnetohydrodynamics based on a variable eddington tensor},
        date={2012},
     journal={The Astrophysical Journal Supplement Series},
      volume={199},
      number={1},
       pages={14},
         url={https://app.dimensions.ai/details/publication/pub.1005319601 and
  http://arxiv.org/pdf/1201.2223},
}

\bib{pub.1041264084}{article}{
      author={Jiang, Yan-Fei},
      author={Stone, James~M.},
      author={Davis, Shane~W.},
       title={An algorithm for radiation magnetohydrodynamics based on solving
  the time-dependent transfer equation},
        date={2014},
     journal={The Astrophysical Journal Supplement Series},
      volume={213},
      number={1},
       pages={7},
         url={https://app.dimensions.ai/details/publication/pub.1041264084 and
  https://iopscience.iop.org/article/10.1088/0067-0049/213/1/7/pdf},
}

\bib{jin1999efficient}{article}{
      author={Jin, Shi},
       title={Efficient asymptotic-preserving (ap) schemes for some multiscale
  kinetic equations},
        date={1999},
     journal={SIAM Journal on Scientific Computing},
      volume={21},
      number={2},
       pages={441\ndash 454},
}

\bib{JinReview}{article}{
      author={Jin, Shi},
       title={Asymptotic preserving (ap) schemes for multiscale kinetic and
  hyperbolic equations},
        date={2010},
}

\bib{jin1998diffusive}{article}{
      author={Jin, Shi},
      author={Pareschi, Lorenzo},
      author={Toscani, Giuseppe},
       title={Diffusive relaxation schemes for multiscale discrete-velocity
  kinetic equations},
        date={1998},
     journal={SIAM Journal on Numerical Analysis},
      volume={35},
      number={6},
       pages={2405\ndash 2439},
}

\bib{jin2009uniformly}{article}{
      author={Jin, Shi},
      author={Tang, Min},
      author={Han, Houde},
       title={A uniformly second order numerical method for the one-dimensional
  discrete-ordinate transport equation and its diffusion limit with interface},
        date={2009},
     journal={Networks \& Heterogeneous Media},
      volume={4},
      number={1},
       pages={35},
}

\bib{KADIOGLU20108313}{article}{
      author={Kadioglu, Samet~Y.},
      author={Knoll, Dana~A.},
      author={Lowrie, Robert~B.},
      author={Rauenzahn, Rick~M.},
       title={A second order self-consistent imex method for radiation
  hydrodynamics},
        date={2010},
        ISSN={0021-9991},
     journal={Journal of Computational Physics},
      volume={229},
      number={22},
       pages={8313 \ndash  8332},
  url={http://www.sciencedirect.com/science/article/pii/S0021999110004122},
}

\bib{klar1998asymptotic}{article}{
      author={Klar, Axel},
       title={An asymptotic-induced scheme for nonstationary transport
  equations in the diffusive limit},
        date={1998},
     journal={SIAM journal on numerical analysis},
      volume={35},
      number={3},
       pages={1073\ndash 1094},
}

\bib{2001Klar}{article}{
      author={Klar, Axel},
      author={Schmeiser, Christian},
       title={Numerical passage from radiative heat transfer to nonlinear
  diffusion models},
        date={2001},
     journal={Math. Meth. Mod. Appl. Sci.},
      volume={11},
      number={5},
       pages={749\ndash 767},
}

\bib{KFJ16}{article}{
      author={Küpper, Kerstin},
      author={Frank, Martin},
      author={Jin, Shi},
       title={An asymptotic preserving two-dimensional staggered grid method
  for multiscale transport equations},
        date={2016},
     journal={SIAM Journal on Numerical Analysis},
      volume={54},
      number={1},
       pages={440\ndash 461},
}

\bib{larsen1989asymptotic}{article}{
      author={Larsen, Edward~W},
      author={Morel, Jim~E},
       title={Asymptotic solutions of numerical transport problems in optically
  thick, diffusive regimes ii},
        date={1989},
}

\bib{larsen1987asymptotic}{article}{
      author={Larsen, Edward~W},
      author={Morel, Jim~E},
      author={Miller~Jr, Warren~F},
       title={Asymptotic solutions of numerical transport problems in optically
  thick, diffusive regimes},
        date={1987},
     journal={Journal of Computational Physics},
      volume={69},
      number={2},
       pages={283\ndash 324},
}

\bib{LM}{article}{
      author={Lemou, Mohammed},
      author={Mieussens, Luc},
       title={A new asymptotic preserving scheme based on micro-macro
  formulation for linear kinetic equations in the diffusion limit},
        date={2008},
     journal={SIAM Journal on Scientific Computing},
      volume={31},
      number={1},
       pages={334\ndash 368},
}

\bib{li2011central}{article}{
      author={Li, Fengyan},
      author={Xu, Liwei},
      author={Yakovlev, Sergey},
       title={Central discontinuous galerkin methods for ideal mhd equations
  with the exactly divergence-free magnetic field},
        date={2011},
     journal={Journal of Computational Physics},
      volume={230},
      number={12},
       pages={4828\ndash 4847},
}

\bib{lowrie1999coupling}{article}{
      author={Lowrie, RB},
      author={Morel, JE},
      author={Hittinger, JA},
       title={The coupling of radiation and hydrodynamics},
        date={1999},
     journal={The astrophysical journal},
      volume={521},
      number={1},
       pages={432},
}

\bib{lowrie2008radiative}{article}{
      author={Lowrie, Robert~B},
      author={Edwards, Jarrod~D},
       title={Radiative shock solutions with grey nonequilibrium diffusion},
        date={2008},
     journal={Shock Waves},
      volume={18},
      number={2},
       pages={129\ndash 143},
}

\bib{MCCLARREN20087561}{article}{
      author={McClarren, Ryan~G.},
      author={Evans, Thomas~M.},
      author={Lowrie, Robert~B.},
      author={Densmore, Jeffery~D.},
       title={Semi-implicit time integration for pn thermal radiative
  transfer},
        date={2008},
        ISSN={0021-9991},
     journal={Journal of Computational Physics},
      volume={227},
      number={16},
       pages={7561 \ndash  7586},
  url={http://www.sciencedirect.com/science/article/pii/S0021999108002489},
}

\bib{mieussens2013asymptotic}{article}{
      author={Mieussens, Luc},
       title={On the asymptotic preserving property of the unified gas kinetic
  scheme for the diffusion limit of linear kinetic models},
        date={2013},
     journal={Journal of Computational Physics},
      volume={253},
       pages={138\ndash 156},
}

\bib{mihalas1984foundations}{article}{
      author={Mihalas, D},
      author={Mihalas, BW},
       title={Foundations of radiation hydrodynamics oxford university press},
        date={1984},
     journal={New York},
}

\bib{MOREL1996445}{article}{
      author={Morel, J.E.},
      author={Wareing, Todd~A.},
      author={Smith, Kenneth},
       title={A linear-discontinuous spatial differencing scheme forsnradiative
  transfer calculations},
        date={1996},
        ISSN={0021-9991},
     journal={Journal of Computational Physics},
      volume={128},
      number={2},
       pages={445 \ndash  462},
  url={http://www.sciencedirect.com/science/article/pii/S0021999196902235},
}

\bib{pomraning2005equations}{book}{
      author={Pomraning, Gerald~C},
       title={The equations of radiation hydrodynamics},
   publisher={Courier Corporation},
        date={2005},
}

\bib{sekora2009higher}{article}{
      author={Sekora, Michael},
      author={Stone, James},
       title={A higher-order godunov method for radiation hydrodynamics:
  Radiation subsystem},
        date={2009},
     journal={Communications in Applied Mathematics and Computational Science},
      volume={4},
      number={1},
       pages={135\ndash 152},
}

\bib{SEKORA20106819}{article}{
      author={Sekora, Michael~D.},
      author={Stone, James~M.},
       title={A hybrid godunov method for radiation hydrodynamics},
        date={2010},
        ISSN={0021-9991},
     journal={Journal of Computational Physics},
      volume={229},
      number={19},
       pages={6819 \ndash  6852},
  url={http://www.sciencedirect.com/science/article/pii/S0021999110002858},
}

\bib{SKTS}{article}{
      author={Stamnes, Knut},
      author={Tsay, S-Chee},
      author={Wiscombe, Warren},
      author={Jayaweera, Kolf},
       title={Numerically stable algorithm for discrete-ordinate-method
  radiative transfer in multiple scattering and emitting layered media},
        date={1988},
     journal={Applied optics},
      volume={27},
      number={12},
       pages={2502\ndash 2509},
}

\bib{stone2008athena}{article}{
      author={Stone, James~M},
      author={Gardiner, Thomas~A},
      author={Teuben, Peter},
      author={Hawley, John~F},
      author={Simon, Jacob~B},
       title={Athena: a new code for astrophysical mhd},
        date={2008},
     journal={The Astrophysical Journal Supplement Series},
      volume={178},
      number={1},
       pages={137},
}

\bib{sun2015asymptotic}{article}{
      author={Sun, Wenjun},
      author={Jiang, Song},
      author={Xu, Kun},
       title={An asymptotic preserving unified gas kinetic scheme for gray
  radiative transfer equations},
        date={2015},
     journal={Journal of Computational Physics},
      volume={285},
       pages={265\ndash 279},
}

\bib{sun2020multiscale}{article}{
      author={Sun, Wenjun},
      author={Jiang, Song},
      author={Xu, Kun},
      author={Cao, Guiyu},
       title={Multiscale simulation for the system of radiation hydrodynamics},
        date={2020},
     journal={Journal of Scientific Computing},
      volume={85},
      number={2},
       pages={1\ndash 24},
}

\bib{2021Tang}{article}{
      author={Tang, Min},
      author={Wang, Li},
      author={Zhang, Xiao-Jiang},
       title={Accurate front capturing asymptotic preserving method for
  nonlinear grey radiative transport equation},
        date={2021},
     journal={SIAM Journal of Scientific Computing},
      volume={43},
      number={3},
       pages={759–783},
}

\end{biblist}
\end{bibdiv}
\bibliographystyle{siam}

\end{document}